\documentclass[11pt]{article}

\usepackage{graphicx}
\usepackage{marvosym}
\usepackage{amsmath}
\usepackage{amssymb}
\usepackage{mathtools}
\usepackage{amsthm}
\usepackage{tikz}
\usepackage{tikz-cd}
\usepackage[left = 1.5in, right = 1.5in]{geometry}
\usepackage[utf8]{inputenc}
\usepackage[english]{babel}
\usepackage[maxbibnames=10]{biblatex}
\usepackage{hyperref}
\usepackage{todonotes}
\usepackage{xcolor}

\newcommand{\reals}{\mathbb{R}}
\newcommand{\intgr}{\mathbb{Z}}

\renewcommand{\tilde}{\widetilde}
\renewcommand{\subset}{\subseteq}
\renewcommand{\phi}{\varphi}
\renewcommand{\epsilon}{\varepsilon}

\renewcommand{\H}{\mathbb{H}}

\newcommand{\M}{\mathcal{M}}
\newcommand{\T}{\mathcal{T}}

\newcommand{\vol}{\textrm{Vol}}
\newcommand{\prob}{\mathbb{P}}

\newtheorem{theorem}{Theorem}[section]
\newtheorem{lemma}[theorem]{Lemma}
\newtheorem{corollary}[theorem]{Corollary}

\theoremstyle{definition}
\newtheorem{definition}[theorem]{Definition}
\theoremstyle{remark}
\newtheorem{remark}[theorem]{Remark}

\title{Critical Exponents on Hyperbolic Surfaces with Long Boundaries and the Asymptotic Weil-Petersson Form}

\author{Henry Talbott}

\date{}

\begin{document}

\maketitle

\tableofcontents

\begin{abstract}

We study the critical exponent random variable $\delta_X$ on moduli spaces of hyperbolic surfaces with boundary, using the normalized Weil-Petersson measures $d\mu_{WP}$ as probability measures. We use the spine graph construction of Bowditch and Epstein to compare this random variable to the corresponding critical exponent random variable $\delta_\Gamma$ on moduli spaces of metric ribbon graphs with the normalized Kontsevich measures $d\mu_K$, proving an asymptotic convergence-in-mean result in the long boundary length regime. In particular, we show that $d\mu_K$ approximately pulls back to $d\mu_{WP}$ with quantitative uniform estimates.

\end{abstract}

\section{Introduction}

The critical exponent of a hyperbolic surface $X$ (with or without boundary) coarsely controls the number of geodesics on $X$. Specifically:

\begin{theorem}[Selberg, \cite{AS56}, Lalley, \cite{SL89}, Cor. 11.2]
\label{critexpdef}

Let $X$ be a compact hyperbolic surface with or without boundary, and let $L_X(x)$ be the number of simple closed geodesics on $X$ with length at most $x$. Then there exists a constant $0<\delta_X\leq 1$, the \emph{critical exponent} of $X$, such that
\begin{equation}\label{ecritexp}L_X(x)\sim\frac{e^{\delta_Xx}}{\delta_Xx}\end{equation}
where $f(x)\sim g(x)$ denotes $\lim_{x\to\infty}\frac{f(x)}{g(x)}=1$.
	
\end{theorem}

With this definition, the Prime Geodesic Theorem becomes the statement that $\delta_X=1$ for all $X$ without boundary.

Denote by $\M_{g,n}(L_1,...,L_n)$ the moduli space of all genus-$g$, $n$-boundary-component hyperbolic surfaces with prescribed boundary lengths $L_1,...,L_n$, and denote by $\T_{g,n}(L_1,...,L_n)$ the corresponding Teichm{\"u}ller space. This moduli space is equipped with the \emph{Weil-Petersson volume form} $dV_{WP}$, and has finite volume with respect to this form. Therefore, we may interpret $\delta_X$ as a random variable on $\M_{g,n}(L_1,...,L_n)$, and ask about its probability distribution. 

In this paper, we show that this distribution is asymptotically equivalent to the critical exponent distribution on the corresponding moduli space $MRG_{g,n}(L_1,...,L_n)$ of \emph{metric ribbon graphs} with respect to the \emph{Kontsevich volume form} $dV_K$ as the $L_i$ are sent to $\infty$. Like $\M_{g,n}(L_1,...L_n)$, $MRG_{g,n}(L_1,...,L_n)$ possesses a Teichm{\"u}ller space $TRG_{g,n}(L_1,...,L_n)$. Critical exponents of metric graphs are much more amenable to existing combinatorial and analytic tools, and in particular can be easily calculated in computer algebra systems via the machinery of graph zeta functions - see section \ref{graphzeta} for more details.

Define $\delta_X:\M_{g,n}(L_1,...,L_n)\to\reals$ to be the function sending each surface to its critical exponent, and define $\delta_\Gamma:MRG_{g,n}(L_1,...,L_n)\to\reals$ analogously. There exists a homeomorphism
$$\Phi:\M_{g,n}(L_1,...,L_n)\to MRG_{g,n}(L_1,...,L_n)$$
defined via the \emph{spine graph construction}, allowing us to pull back $\delta_\Gamma$ and $dV_K$ to $\M_{g,n}(L_1,...,L_n)$ and compare them to $\delta_X$ and $dV_{WP}$. Since $\M_{g,n}(L_1,...,L_n)$ and $MRG_{g,n}(L_1,...,L_n)$ have finite volumes with respect to $dV_{WP}$ and $dV_{K}$, we may normalize these volume forms to obtain probability measures $d\mu_{WP}$ and $d\mu_K$. 

\begin{remark}
Both $\M_{g,n}(L_1,...,L_n)$ and $MRG_{g,n}(L_1,...,L_n)$ possess cell complex structures, and in fact $\Phi$ is a diffeomorphism between the unions of top-dimensional cells on the two complexes. For a proof of this nontrivial fact, see Appendix \ref{nontrivalent}.
\end{remark}

To simplify the statements of our results, it is useful to take advantage of a scaling property of $MRG_{g,n}(L_1,...,L_n)$: for any $\alpha>0$, we can define a map $\rho_{1/\alpha}$ on metric ribbon graphs that scales all edge lengths by $1/\alpha$, in which case

\begin{equation}
\label{erescale}
\rho_{1/\alpha}:MRG_{g,n}(\alpha L_1,...,\alpha L_n)\to MRG_{g,n}(L_1,...,L_n)
\end{equation}

is a diffeomorphism satisfying $\rho_{1/\alpha}^*dV_K=dV_K$ and $\rho^*_{1/\alpha}\delta_{\Gamma}=\alpha\delta_\Gamma$. 

\begin{remark}
Note that this scaling property does not hold for surfaces, as the corresponding scaling map $\Phi^{-1}\circ\rho_{1/\alpha}\circ\Phi$ for $\M_{g,n}(\alpha L_1,...,\alpha L_n)$ does not preserve $dV_{WP}$, and does not scale $\delta_X$ in any uniform way. 
\end{remark}

\subsection{Main Results}

We can now define
$$\Phi_{1/\alpha}=\rho_{1/\alpha}\circ\Phi:\M_{g,n}(\alpha L_1,...,\alpha L_n)\to MRG_{g,n}(L_1,...,L_n)$$
and would like to compare the random variables $(\Phi_{1/\alpha})_*\delta_X$ and $\delta_\Gamma$. Two issues immediately arise.

First, as $\alpha$ tends to $\infty$, a `typical' critical exponent for a surface in $\M_{g,n}(\alpha L_1,...,\alpha L_n)$ will tend to $0$, and we will need to normalize $\delta_X$ via scaling by $\alpha$ for the comparison to $\delta_\Gamma$ to make sense. So we will work with $(\Phi_{1/\alpha})_*(\alpha\delta_X)$ instead of $(\Phi_{1/\alpha})_*\delta_X$.

Second, we need to account for the fact that $\delta_\Gamma$ is naturally a random variable with respect to $d\mu_K$ on $MRG_{g,n}(L_1,...,L_n)$, while $(\Phi_{1/\alpha})_*(\alpha\delta_X)$ is naturally a random variable with respect to $(\Phi_{1/\alpha})_*d\mu_{WP}$, where $d\mu_{WP}$ is the Weil-Petersson probability measure on $\M_{g,n}(L_1,...,L_n)$. We resolve this issue by showing that $(\Phi_{1/\alpha})_*d\mu_{WP}$ is absolutely continuous with respect to $d\mu_K$, so we can apply the Radon-Nikodym Theorem and write
$$(\Phi_{1/\alpha})_*d\mu_{WP}=\frac{(\Phi_{1/\alpha})_*d\mu_{WP}}{d\mu_K}d\mu_K$$
where $\frac{(\Phi_{1/\alpha})_*d\mu_{WP}}{d\mu_K}$ is a measurable function defined $d\mu_K$-almost-everywhere. We then have the following imprecise statement. This statement is a consequence of Theorem \ref{maintheorem}, and is stated precisely as Corollary \ref{cor1}.

\begin{theorem}[Main result, non-technical version]
\label{mainnontechnical}
As $\alpha\to\infty$, the random variables
$$(\Phi_{1/\alpha})_*(\alpha\delta_X) \frac{(\Phi_{1/\alpha})_*d\mu_{WP}}{d\mu_K}$$
on the probability space $(MRG_{g,n}(L_1,...,L_n),d\mu_{K})$ converge in the $L^1$ norm (equivalently, converge in mean) to the random variable $\delta_\Gamma$.
	
\end{theorem}

This result significantly improves upon work of Mondello and Do \cite{GM09,ND10}, whose work implies only pointwise convergence of this random variable. Furthermore, our proof requires new tools to provide uniform error control on portions of the moduli space, and we believe these tools are likely of independent interest.

\begin{remark}

While Theorem \ref{mainnontechnical} is easier to state, for the rest of the paper we will work exclusively on moduli spaces of surfaces. This choice has the downside that the spaces will depend on $\alpha$ (or on another quantity being sent to $\infty$, depending on the setup), so $L^1$ convergence of random variables no longer makes sense. On the other hand, more is known about moduli spaces of surfaces, making them a more reasonable setting for us.
	
\end{remark}

\begin{figure}[h]

\centering
\includegraphics[scale=0.8]{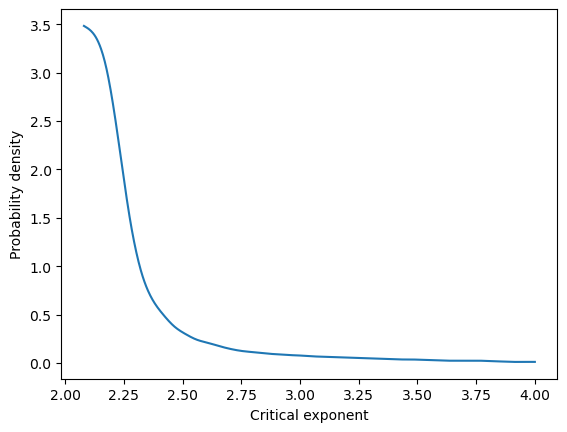}
\caption{An approximation for the distribution of critical exponent values on $MRG_{1,1}(1)$, based on smoothed data obtained from a computer simulation. Notice the lower bound of $3\log(2)\approx 2.08$.}

\label{critexpdist}

\end{figure}

We may also say something about critical exponent distributions as probability measures on $\reals$. Fixing $L_1,...,L_n$, and considering the probability spaces $(\M_{g,n}(\alpha L_1,...,\alpha L_n),d\mu_{WP})$ and $(MRG_{g,n}(L_1,...,L_n),d\mu_K)$, we would like to say that as $\alpha\to\infty$, the probability measure $(\alpha\delta_X)_*d\mu_{WP}$ is related to $(\delta_\Gamma)_*d\mu_K$ via increasingly minor perturbation. We show that $(\alpha\delta_X)_*d\mu_{WP}$ does in fact converge to $(\delta_\Gamma)_*d\mu_K$ in the space of probability measures on $\reals$ in the Wasserstein metric (also known as the Kantorovich metric) at an explicit rate as $\alpha\to\infty$. Recall that the Wasserstein or `sandpile' distance on probability measures $\mu$ and $\nu$ expresses the minimal amount of work (in the sense of physics) to create the distribution $\nu$ by moving the mass of $\mu$. Since $(\delta_\Gamma)_*d\mu_K$ can be approximately computed, we can therefore obtain rough pictures of $(\alpha\delta_X)_*d\mu_{WP}$ for large $\alpha$ (see Figure \ref{critexpdist}).

\

Our main technical result is as follows:

\begin{theorem}
\label{maintheorem}

Fix $\M_{g,n}(L_1,...,L_n)$, and constants $C$ and $L$ satisfying
$$\frac{1}{C}L\leq L_i\leq CL$$ for all $L_i$. Let $f:\reals\to\reals$ be an arbitrary $1$-Lipschitz function satisfying $f(0)=0$. Consider the probability spaces
$$(\M_{g,n}(L_1,...,L_n),d\mu_{WP}),\ (MRG_{g,n}(L_1,...,L_n),d\mu_K)$$
and the spine map
$$\Phi:\M_{g,n}(L_1,...,L_n)\to MRG_{g,n}(L_1,...,L_n)$$
Then the Radon-Nikodym derivative $\frac{\Phi^*d\mu_K}{d\mu_{WP}}$ is well-defined off a set of $\mu_{WP}$-measure $0$, and if $L$ is sufficiently large in terms of $g$, $n$, and $C$,
\begin{equation}\label{emaintheorem}\int_{\M_{g,n}(L_1,...,L_n)}\left|f(\delta_X)-f(\Phi^*\delta_\Gamma)\frac{\Phi^*d\mu_K}{d\mu_{WP}}\right|d\mu_{WP}=O_{g,n,C}\left(\frac{1}{L^{13/12}}\right)\end{equation}

\end{theorem}

This result has two useful consequences. The first is a technical restatement of Theorem \ref{mainnontechnical}.

\begin{corollary}
\label{cor1}

Fix $\M_{g,n}(L_1,...,L_n)$. For $\alpha>0$, define 
$$\Phi_{1/\alpha}:\M_{g,n}(\alpha L_1,...,\alpha L_n)\to MRG_{g,n}(L_1,...,L_n)$$
to be the composition of
$$\Phi:\M_{g,n}(\alpha L_1,...,\alpha L_n)\to MRG_{g,n}(\alpha L_1,...,\alpha L_n)$$
and $\rho_{1/\alpha}$ as defined in equation \ref{erescale}. Then
\begin{equation}\label{ecor1}\int_{\M_{g,n}(\alpha L_1,...,\alpha L_n)}\left|\alpha\delta_X-\Phi_\alpha^*\delta_\Gamma\frac{\Phi^*_\alpha d\mu_K}{d\mu_{WP}}\right|d\mu_{WP}=O_{g,n,L_1,...,L_n}\left(\frac{1}{\alpha^{1/12}}\right)\end{equation}
	
\end{corollary}

While the $L^1$ norm here is on $(\M_{g,n}(\alpha L_1,...,\alpha L_n),d\mu_{WP})$ rather than $(MRG_{g,n}(L_1,...,L_n),d\mu_K)$, the fact (see Lemma \ref{absolutecontinuity}) that $d\mu_{WP}$ is absolutely continuous with respect to $d\mu_K$ and vice-versa allows us to derive Theorem \ref{mainnontechnical} from this corollary.

The second concerns Wasserstein convergence. Formally, for two probability measures $\mu$ and $\nu$ on $\reals$, their Wasserstein distance is
$$d_W(\mu,\nu)=\inf_\Gamma \int_{\reals^2} |x-y|d\Gamma(x,y)$$
where the infimum is over all probability measures $\Gamma$ on $\reals^2$ with $\mu$ and $\nu$ as marginals, i.e.
$$\int_\reals d\Gamma(x_0,y)=\mu(x_0),\ \int_\reals d\Gamma(x,y_0)=\nu(y_0)$$
for all $x_0$ and $y_0$. As mentioned previously, the main advantage of working with the Wasserstein metric over, say, the stronger total variation metric is that it is insensitive to very small perturbations of the underlying measures. For example, if $\mu$ is a unit point mass at $1$ and $\nu$ is a unit point mass at $0.99$, then $d_{TV}(\mu,\nu)=2$ while $d_W(\mu,\nu)=0.01$.
\begin{corollary}
\label{cor2}

Fix $\M_{g,n}(L_1,...,L_n)$. For $\alpha>0$, let $\delta_X$ be the critical exponent function on the probability space $(\M_{g,n}(\alpha L_1,...,\alpha L_n),d\mu_{WP})$, and let $\delta_\Gamma$ be the critical exponent function on $(MRG_{g,n}(L_1,...,L_n),d\mu_K)$. Then
\begin{equation}\label{ecor2}d_W((\alpha\delta_X)_*(d\mu_{WP}),(\delta_\Gamma)_*(d\mu_K))=O\left(\frac{1}{\alpha^{1/12}}\right)\end{equation}
where $d_W(\cdot,\cdot)$ is the Wasserstein $1$-metric on probability measures on $\reals$.
	
\end{corollary}

With the help of a computer algebra system, we may estimate $(\delta_\Gamma)_*(d\mu_K)$ for the space $MRG_{1,1}(1)$; the result is depicted in Figure \ref{critexpdist}. Corollary \ref{cor2} then implies that, after appropriate rescaling, the critical exponent distribution on $\M_{1,1}(\alpha)$ should converge to the distribution depicted in Figure \ref{critexpdist} as $\alpha$ is sent to $\infty$.

We would also like to point out our primary technical result, since others may find it independently useful. 

\begin{lemma}
\label{earlyformratio}

Fix $X\in\M_{g,n}(L_1,...,L_n)$ with $\Gamma=\Phi(X)$, and assume that $\Gamma$ has no geodesics of length $\leq B$ for some constant $B$ depending only on $g$, $n$, and $C$. Then
$$\frac{\Phi^*dV_K}{dV_{WP}}|_X=O_{g,n,C}(1)$$
If in fact every edge of $\Gamma$ has length at least $\eta$ and $\eta$ is sufficiently large in terms of $g$, $n$, and $C$, then
$$\frac{\Phi^*dV_K}{dV_{WP}}|_X=1+O_{g,n,C}(e^{-\eta/2})$$
\end{lemma}

An important part of proving this lemma is that, in the no-short-geodesics case, we can find a `nice' set of geodesics with length functions that serve as local coordinates around $X$, and that only intersect with angles close to $0$ or $\pi$. This is handled in Lemma \ref{goodgeodframe}.

Lastly, note that the methods used to prove Lemma \ref{earlyformratio} recover the pointwise convergence results of Do and Mondello.

\subsection{Background}

The connection between hyperbolic surfaces with boundary and metric ribbon graphs (also called fatgraphs) was initially developed by Harer-Mumford-Thurston, and moduli spaces of metric ribbon graphs were used in Kontsevich's proof of Witten's conjecture \cite{JH86, MK92}. However, the geometric construction presented here is due to Penner and Bowditch-Epstein \cite{RP87, BE88}. Convergence of the Weil-Petersson form to the Kontsevich form in the fixed-topological-type, long-boundary-length regime was initially studied by Do and Mondello \cite{ND10, GM09}. Connections between the Kontsevich volume form and complementary subsurfaces of a given hyperbolic surface have also been studied by Arana-Herrara-Calderon \cite{AHC22}.

Furthermore, the asymptotic  counts of simple multicurves with bounded lengths, averaged over the respective moduli spaces, are the same in the hyperbolic and combinatorial settings, and in fact coincide with the Masur-Veech volumes of the principal stratum of the moduli space of quadratic differentials. This result follows from the works of Mirzakhani, Delecroix-Goujard-Zograf-Zorich, and Andersen-Borot-Charbonnier-Delecroix-Giacchetto-Lewa{\'n}ski-Wheeler \cite{MM08a, MM08b, DGZZ21,ABCDGLW19}. Giacchetto additionally shows that the respective geometric and combinatorial averages may be computed (before taking the asymptotic) via essentially identical topological recursion formulas \cite{AG21}.

We make significant use of the method of topological recursion to calculate volumes of subsets of $\M_{g,n}(L_1,...,L_n)$. This technique was first developed by Mirzakhani, and has since found numerous applications in problems relating to random surfaces \cite{MM13}. Recently, analogous recursion techniques were established for $MRG_{g,n}(L_1,...,L_n)$ by Andersen-Borot-Charbonnier-Giacchetto-Lewa\'nski-Wheeler \cite{ABCGLW21}. 

Note that the asymptotic properties of the length spectra of surfaces and metric graphs are well-understood in the large-genus regime; in both cases, the spectrum is asymptotically described by a Poisson point process, with identical intensities. This is due to Mirzakhani-Petri for the hyperbolic case and Jason-Louf and Barazer-Giacchetto-Liu for the combinatorial case \cite{MP19, JL23, BGL23}.

\subsection{Setup and Asymptotic Notation}

Throughout this paper, we will fix $g\in\intgr_{\geq 0}$, $n\in\intgr_{\geq 1}$, and $C\in\reals_{\geq 1}$, with the caveat that if $g=0$, then $n\geq 3$. We use $\omega$ to denote the Weil-Petersson symplectic 2-form on $\M_{g,n}(L_1,...,L_n)$ or $\T_{g,n}(L_1,...,L_n)$, while $\omega_K$ will denote the Kontsevich $2$-form on $MRG_{g,n}(L_1,...,L_n)$ and $TRG_{g,n}(L_1,...,L_n)$. The respective volume forms will be denoted $dV_{WP}$ and $dV_K$.

\begin{remark}
Note that geodesics are generally assumed to be closed, unless otherwise specified.
\end{remark}

For real-valued functions $f$ and $g$ on an arbitrary domain, $f=O(g)$ denotes that there exists some constant $A>0$ such that $|f|\leq A|g|$ for all inputs, and $f=\Omega(g)$ denotes that there exists $A>0$ such that $|f|\geq A|g|$. If constants are attached to $O$ or $\Omega$ as subscripts, then $A$ is allowed to depend on those constants. For example, $f=O_{a,b}(g)$ means that for any fixed choice of $a$ and $b$, there exists a $A=A(a,b)>0$ such that $|f|\geq A|g|$. The related notation $f\asymp g$ means that $f=O(g)$ and $g=O(f)$, or explicitly there exists some $A>0$ such that $\frac{1}{A}f\leq g\leq Af$. The notation $f\prec g$ is another way of writing $f=O(g)$.

We will also make the following definition:

\begin{definition}
Let $L,L_1,...,L_n,C$ be positive constants. We say the $L_i$ are $C$-\emph{controlled by} $L$ if for all $L_i$ we have
$$\frac{1}{C}L\leq L_i\leq CL$$
\end{definition}

\begin{remark}

Essentially all implicit constants in this paper will depend on the $g$ and $n$ defining $\M_{g,n}(L_1,...,L_n)$ and $MRG_{g,n}(L_1,...,L_n)$, and on the $C$ used to control the growth of the $L_i$ as above. Therefore, we will \textbf{we will generally drop the subscripts in the corresponding asymptotic notation}, so $O_{g,n,C}(...)$, $\Omega_{g,n,C}(...)$, and $\prec_{g,n,C}$ become $O(...)$, $\Omega(...)$, and $\prec$ throughout the text.
	
\end{remark}

\subsection{Organization}

Section 2 provides background material on metric ribbon graphs and spine graphs, and explains how critical exponent distributions may be computed. Section 3 provides volume bounds on various subsets of $\M_{g,n}(L_1,...,L_n)$ using Mirzakhani's topological recursion scheme. Section 4 establishes many preliminary geometric and trigonometric results, and provides quantitative estimates for $|\delta_X-\delta_{\Phi(X)}|$. Sections 5 and 6 are devoted to proving Lemma \ref{earlyformratio}, which appears as Lemma \ref{ratiocontrol}. Section 7 establishes the main results. Appendix A contains elementary trigonometric formulas, while appendix B contains some technical arguments related to certain loci of $\M_{g,n}(L_1,...,L_n)$ and $MRG_{g,n}(L_1,...,L_n)$ having full measure.

\subsection{Acknowledgements}

The author would like to thank Alex Wright for suggesting this problem, Richard Canary, Gabriele Mondello, and Alex Wright for helpful conversations, and Alessandro Giacchetto for alerting the author to several useful papers in the literature. The author would also like to thank the reviewer for their insightful feedback. This material is based upon work supported by the National Science Foundation Graduate Research Fellowship under Grant No. DGE 2241144.

\section{Background and Applications}

\subsection{Metric Ribbon Graphs}

A \emph{ribbon graph} $\Gamma$ is a finite topological graph such that each vertex carries a cyclic ordering of its incident edges; when ribbon graphs are drawn, edges are usually arranged around a vertex counterclockwise according to the order. Although they do not need to be planar, ribbon graphs have well-defined \emph{faces}. One can imagine thickening each edge into a ribbon and gluing the ribbons at each vertex according to the order; this produces an oriented topological surface with boundary components corresponding to the faces of the ribbon graph. Furthermore, this construction allows us to assign a well-defined topological type $(g,n)$ to any ribbon graph.

A \emph{metric ribbon graph} is a ribbon graph with a length assigned to each edge; these lengths induce lengths on faces in the expected way. Likewise, an isometry of metric ribbon graphs is a metric graph isometry that preserves the vertex orderings. The \emph{combinatorial moduli space} $MRG_{g,n}(L_1,...,L_n)$ is the space of all metric ribbon graphs (up to isometry) with topological type $(g,n)$ and faces with lengths $L_1$ through $L_n$. This space is the mapping-class-group quotient of the corresponding \emph{combinatorial Teichm{\"u}ller space} $TRG_{g,n}(L_1,...,L_n)$, which is defined in terms of ribbon graphs with additional marking data analogously to $\T_{g,n}(L_1,...,L_n)$. 

The space $MRG_{g,n}(L_1,...,L_n)$ possesses a natural cell structure, and this cell structure may be pulled back to $\M_{g,n}(L_1,...,L_n)$ via the spine homeomorphism. Cells in this structure represent distinct topological isomorphism classes of ribbon graphs, and maximal-dimension cells correspond to trivalent ribbon graphs. For any $MRG_{g,n}(L_1,...,L_n)$, there is a finite number of cells, dependent only on $g$ and $n$.

Furthermore, Andersen-Borot-Charbonnier-Giacchetto-Lewa{\'n}ski-Wheeler in \cite{ABCGLW21} study the geometry of $MRG_{g,n}(L_1,...,L_n)$ and $TRG_{g,n}(L_1,...,L_n)$, and establish parallels to many results about $\M_{g,n}(L_1,...,L_n)$. In particular, $(TRG_{g,n}(L_1,...,L_n),\omega_K)$ is a symplectic manifold, and possesses analogs of Fenchel-Nielsen coordinate systems corresponding to topological pants decompositions. The length coordinates in these systems are defined as lengths of closed geodesics as one might expect, and every simple closed geodesic $\gamma$ gives rise to a twist flow $\frac{\partial}{\partial\tau_\gamma}$. However, twist coordinates and twist flows come with an important technical caveat; see Remark \ref{twistproblem}.

Closed geodesics on ribbon graphs are non-backtracking closed paths; notice this definition makes sense regardless of a choice of metric. Also, any closed path is homotopic to a unique geodesic representative. The geodesic length spectrum for a metric ribbon graph is defined as expected, and a theorem analogous to Theorem \ref{critexpdef} holds. In particular, every metric ribbon graph $\Gamma$ has a well-defined critical exponent $0<\delta_\Gamma$. As mentioned above, $MRG_{g,n}(L_1,...,L_n)$ carries a natural symplectic $2$-form $\omega_K$, the \emph{Kontsevich form}, and a corresponding volume form $dV_K$.

Note that defining intersections of geodesics on ribbon graphs is somewhat subtle; see Definition \ref{graphintersectiondef}.

Throughout the paper, we will make use of both cutting and gluing operations along geodesics on metric ribbon graphs. These are defined essentially how one would expect: for cutting, one thinks of slightly inflating a given graph $\Gamma$ into a union of ribbons, cutting the ribbons along some geodesic $\tilde{\gamma}$, and then deflating the ribbons to obtain the new graph(s). Gluing is analogous, with the caveat that gluing two geodesics with a certain twist offset $\tau$ may only be defined for $\tau$ in a dense open subset of $\reals$; see Remark \ref{twistproblem}.

\begin{remark}
The length spectrum and critical exponent of a metric ribbon graph is independent of the ribbon structure (i.e., the cyclic edge-ordering information attached to each vertex).
\end{remark}

\begin{remark}
Notice that metric graphs may have critical exponent greater than $1$, something that is not possible in the surface case. This is essentially due to the collar lemma; no equivalent constraint exists on the geometry of metric graphs. For example, consider a graph $\Gamma_\epsilon$ with a single vertex and two loops of length $\epsilon$. Then $\lim_{\epsilon\to 0}\delta_{\Gamma_\epsilon}=\infty$. 
\end{remark}

\begin{remark}
\label{twistproblem}
While Andersen-Borot-Charbonnier-Giacchetto-Lewa{\'n}ski-Wheeler construct in \cite{ABCGLW21} a combinatorial analogue of Fenchel-Nielsen coordinates for metric ribbon graphs, their notion of twisting comes with a technical caveat. Specifically, for a given geodesic $\gamma$ on a metric ribbon graph $\Gamma$, the graph $\Gamma_t$ produced by twisting $\Gamma$ along $\gamma$ for distance $t$ (via a cut-and-glue operation analogous to twisting on surfaces) may fail to be well-defined for $t$ lying in some subset $B\subset\reals$. 

However, Andersen-Borot-Charbonnier-Giacchetto-Lewa{\'n}ski-Wheeler prove (Proposition $2.34$ in \cite{ABCGLW21}) that $B$ is countable, and one may extract from their proof that if $\Gamma$ is trivalent, then there exists is some $\epsilon>0$ such that $\Gamma_t$ is well-defined for all $t\in (-\epsilon,\epsilon)$. In particular, infinitesimal twists of trivalent graphs are always well-defined. For more information on gluing and twisting for metric ribbon graphs, see subsection 2.4 in \cite{ABCGLW21}.
	
\end{remark}

\subsection{Spine Graphs}

For any $X\in \M_{g,n}(L_1,...,L_n)$, define its spine to be the set of points $x\in X$ such that $d(x,\partial X)$ is realized by multiple points on $\partial X$. This spine has a nice geometric structure.

\begin{theorem}[Bowditch-Epstein, \cite{BE88}]
\label{spineconst}

Fix $X\in \M_{g,n}(L_1,...,L_n)$, let $E$ be the set of $x\in X$ such that $d(x,\partial X)$ is realized by exactly two boundary points, and let $V$ be the set of $x\in X$ such that $d(x,\partial X)$ is realized by at least three boundary points. Then:

\begin{enumerate}
\item $V$ is a finite set.
\item $E$ is a finite set of disjoint open geodesic arcs with endpoints in $V$.
\item $E\cup V=S_X$ is a metric ribbon graph embedded in $X$, with the orientation on $X$ determining the ribbon structure.
\item $X$ deformation retracts onto $E$, so $X$ and $E$ have the same topological type.
\end{enumerate}
	
\end{theorem}

We can now construct a surface-to-spine map
$$\Phi:\M_{g,n}(L_1,...,L_n)\to MRG_{g,n}(L_1,...,L_n)$$
Topologically, $\Phi(X)$ will be the same as the spine graph $S_X$ defined above. However, defining the lengths of its edges is more subtle. First draw for each $x\in V$ and $y\in \partial X$ realizing $d(x,\partial X)$ the geodesic arc from $x$ to $y$. These arcs are the \emph{ribs} of the spine, and divide $X$ into a union of hexagons we will call \emph{corridors}. Each corridor contains one edge of $S_X$ and two boundary arcs of $\partial X$ with equal length; we define the length of the corresponding edge in $\Phi(X)$ to be the length of either boundary arc. The key property of this construction is that the lengths of the faces of $\Phi(X)$ are identical to the boundary lengths of $X$, so $\Phi$ has the correct codomain.

It will also be handy to decompose a surface into polygons in a slightly different way: for each corridor, one may draw an arc contained in the corridor between the two boundary components of that corridor, and then find the minimal-length representative of the homotopy class of that arc relative to the boundary. These representatives are called \emph{intercostals}, and cut a given surface into polygons we will call \emph{sectors}. In the case where the spine graph is trivalent, all sectors are hexagonal, and it is this case that will appear frequently throughout the paper.

\begin{figure}[h]

\centering
\includegraphics[scale=0.16]{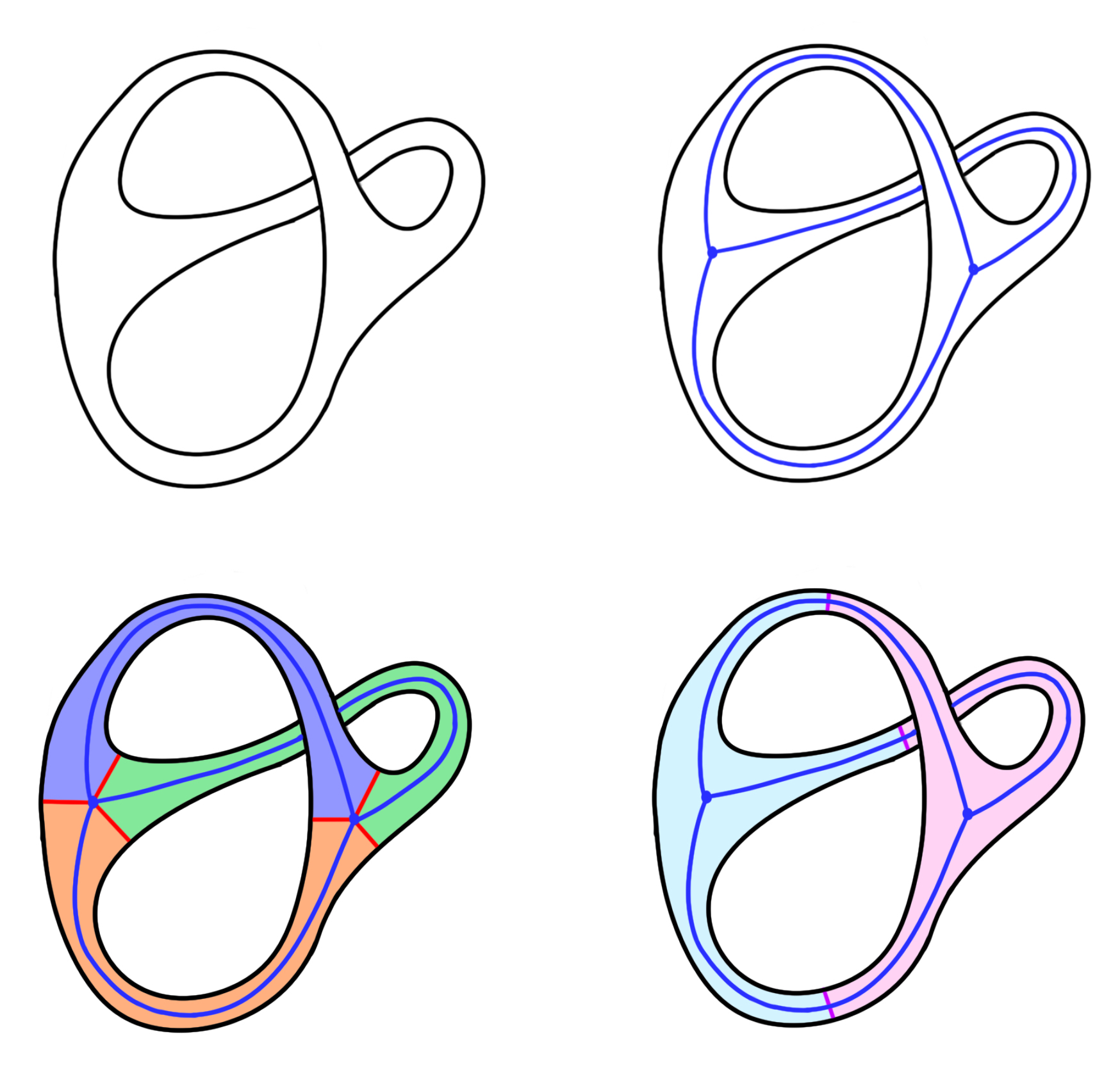}
\caption{Examples of constructions used in this paper. Top left: a surface $X$ with $g=1$ and $n=1$ (i.e., a one-holed torus). Top right: the same surface with its spine graph marked. Notice that this metric ribbon graph has a single face. Bottom left: the same surface, with its ribs added in red and the three corridors shaded. Bottom right: the same surface, with its intercostals added in magenta and its two sectors shaded.}

\label{spineexample}

\end{figure}

For an example of the spine graph, corridor decomposition, and sector decomposition constructions, see Figure \ref{spineexample}. 

\begin{theorem}[Bowditch-Epstein, \cite{BE88}, Thm. 9.5]
\label{spinemap}

The function $\Phi:\M_{g,n}(L_1,...,L_n)\to MRG_{g,n}(L_1,...,L_n)$ defined above is a homeomorphism.
	
\end{theorem}

Later, we will see that when the $L_i$ are large and roughly of the same order, $\Phi$ both approximately preserves critical exponents and approximately preserves the volume form.

\subsection{Graph Zeta Functions and Critical Exponent Distribution}

This section outlines one motivation for relating surface and graph critical exponents. Let $\Gamma$ be a metric graph (or a metric ribbon graph, although the ribbon structure is irrelevant here). Double each edge, and label each resulting pair with opposite orientations. Then construct a square matrix $M_{\Gamma}=(a_{ij})$ with a row and column for each edge in this new directed graph. Set $a_{ij}=z^{l(i)}$ if the initial vertex of edge $j$ is the terminal vertex of edge $i$ (i.e. edge $i$ feeds into edge $j$) and edge $j$ is not the reverse of edge $i$. Otherwise, set $a_{ij}=0$. The \emph{Ihara zeta function} $\zeta_{\Gamma}(z)$ of $\Gamma$ is defined as the resolvent
$$\zeta_\Gamma(z)=\frac{1}{\det(M_\Gamma-I)}$$
This function is the reciprocal of a complex function $p_\Gamma(z)=\det(M_\Gamma-I)$; furthermore, if all edge lengths are integers, then $p_\Gamma(z)$ is a polynomial. Moreover, we have the following relation between the poles of $\zeta_\Gamma$ and the critical exponent $\delta_\Gamma$ (see chapter 10 of \cite{AT10} for a proof\footnote{Terras only covers the case when all edge lengths are integers, but the full result is easily obtained by rational approximation and variable substitution.}):

\begin{theorem}
\label{graphzeta}

Let $\Gamma$ be a finite metric graph, and let $\delta_\Gamma$ be its critical exponent. With $p_{\Gamma}(z)$ defined as above,
\begin{equation}\label{ezetazero}\delta_\Gamma=\min_{p_\Gamma(z)=0}\log(1/|z|)\end{equation}
	
\end{theorem}

Formula \ref{ezetazero} is easily implemented in computer algebra systems such as Sagemath, allowing for concrete calculations. For example, consider $MRG_{1,1}(1)$. In general, $MRG_{g,n}(L_1,...,L_n)$ will be a cell complex; here, a trivalent ribbon graph satisfying the topological information $g=n=1$ must be topologically isomorphic to the spine graph in figure \ref{spineexample}, so $MRG_{1,1}(1)$ has a single top-dimensional cell with dimension $2$. The non-trivalent graphs in this space lie on the boundary of the cell and can be ignored; see section 4. We may parameterize these trivalent graphs by labeling the edges $x$, $y$, and $z$, in which case the boundary length condition becomes $x+y+z=1$.  This yields a finite cover of the top-dimensional cell of $MRG_{1,1}(1)$:
$$F_{x,y}=\{(x,y)\in\reals^2:0<x<1,0<y<1,x+y<1\}$$

\begin{figure}[h]

\centering
\includegraphics[scale=0.3]{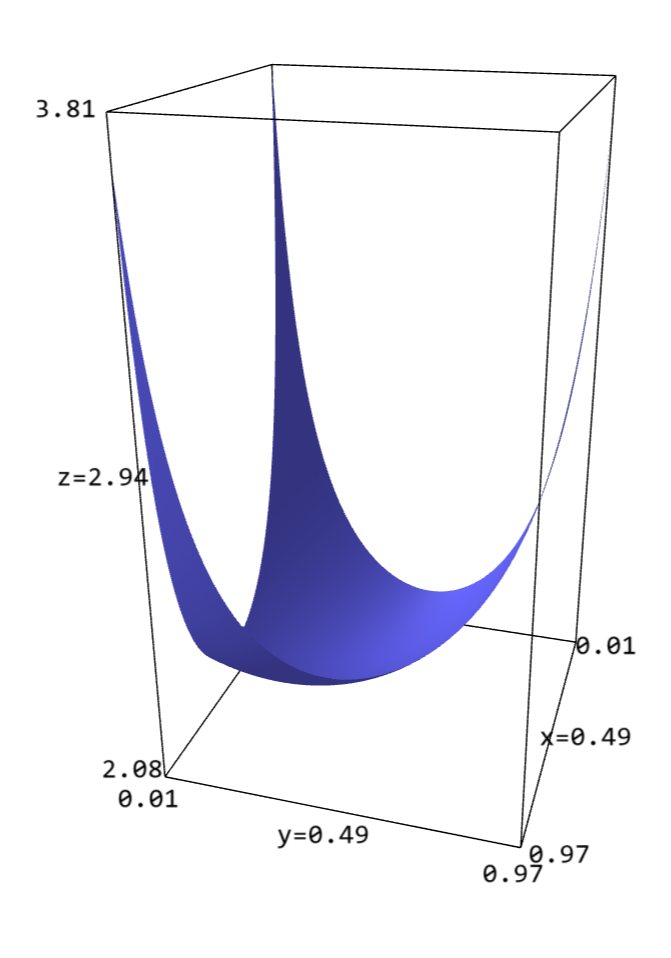}
\caption{A graph of $\delta_\Gamma$ as $\Gamma$ varies over the domain $F_{x,y}$ defined above.}

\label{critexpgraph}

\end{figure}

In this case, the Kontsevich $2$-form is also the volume form, and in the $(x,y)$-coordinates becomes a scalar multiple of the Euclidean volume form $dx\wedge dy$. We can now compute an approximate graph for $\delta_\Gamma$ as a function of $x$ and $y$, as well as the probability density function of $(\delta_\Gamma)_*(dV_K)$; see Figures \ref{critexpgraph} and \ref{critexpdist}.\footnote{All computations were done in Sagemath, and the author is happy to share the code used upon request.}

Corollary \ref{cor2} tell us that as $L\to\infty$ and under suitable normalization, the distribution of critical exponents on $\M_{1,1}(L)$ approaches the distribution shown in Figure \ref{critexpdist}.

\section{Pants Defects and Integration Formulas}

\subsection{Hyperbolic Surfaces and Topological Recursion}

All volumes in this paper will be respect to either $dV_{WP}$ (for volumes on $\M_{g,n}(L_1,...,L_n)$) or $dV_K$ (for volumes on $MRG_{g,n}(L_1,...,L_n)$).

Given measurable sets $U\subset \M_{g,n}(L_1,...,L_n)$ and $V\subset MRG_{g,n}(L_1,...,L_n)$, make the following definitions:
$$V_{g,n}(L_1,...,L_n)=\vol(\M_{g,n}(L_1,...,L_n))$$
$$V^K_{g,n}(L_1,...,L_n)=\vol(MRG_{g,n}(L_1,...,L_n))$$
$$\prob_{\M_{g,n}(L_1,...,L_n)}(U)=\frac{\vol_{\M_{g,n}(L_1,...,L_n)}(U)}{V_{g,n}(L_1,...,L_n)}$$
$$\prob_{MRG_{g,n}(L_1,...,L_n)}(V)=\frac{\vol_{MRG_{g,n}(L_1,...,L_n)}(V)}{V^K_{g,n}(L_1,...,L_n)}$$

Mirzakhani proved that the functions $V_{g,n}(L_1,...,L_n)$ are explicit polynomials: 

\begin{theorem}[Mirzakhani, \cite{MM07}, Thm. 6.1]
\label{volpoly}

We have
\begin{equation}\label{evolpoly}V_{g,n}(L_1,...,L_n)=\sum_{\substack{d\in\intgr^n_{\geq 0}\\ |d|\leq 3g-3+n}} C_g(d)L_1^{2d_1}\cdots L_n^{2d_n}\end{equation}
where the $C_g(d)$ are positive constants depending only on $g$ and $d$.
	
\end{theorem}

One can give formulas for the $C_g(d)$ in terms of intersection numbers of divisors of the compactification $\overline{\M}_{g,n}(L_1,...,L_n)$, viewed as an algebraic variety, but here these coefficients will simply be absorbed into the various implicit constants. 

We now need several probability computations. The first concerns the probability of finding a short geodesic, and is the analogous version of Theorem 4.2 in \cite{MM13} for the long-boundary-length regime.

\begin{lemma}	
\label{shortloopvol}
Consider $\M_{g,n}(L_1,...,L_n)$ with $g$ and $n$ fixed, and choose a constant $\epsilon>0$. Assume that the $L_i$ are $C$-controlled by $L$. Let $U(\epsilon)$ be the subset of $\M_{g,n}(L_1,...,L_n)$ containing surfaces $X$ with at least one closed geodesic with length at most $\epsilon$. Then for $L$ sufficiently large in terms of $\epsilon$,
$$\prob_{\M_{g,n}(L_1,...,L_n)}(U(\epsilon))=O\left(\frac{\epsilon^2}{L^2}\right)$$
\end{lemma}

\textbf{Proof.} First, if a surface has a closed geodesic $\gamma$ with length at most $\epsilon$, a standard argument shows that the surface has a \emph{simple} closed geodesic with length at most $\epsilon$. 

Let $c_\epsilon:\M_{g,n}(L_1,...,L_n)\to\intgr_{\geq 0}$ count the number of simple closed geodesics with length at most $\epsilon$ on a given $X$; we ultimately want to bound
$$\frac{1}{V_{g,n}(L_1,...,L_n)}\int_{\M_{g,n}(L_1,...,L_n)}1_{c_\epsilon(X)\neq 0}dV_{WP}$$
Since $c_\epsilon$ is $\intgr_{\geq 0}$-valued, $1_{c_\epsilon(X)\neq 0}\leq c_\epsilon(X)$ for all $X$, and our first step will be to replace $1_{c_\epsilon(X)\neq 0}$ with $c_\epsilon(X)$ above. The resulting expression is the average number of simple closed geodesics of length at most $\epsilon$ on a random surface in $\M_{g,n}(L_1,...,L_n)$ with respect to the Weil-Petersson probability measure.

The set of all geodesics $\gamma$ on any $X\in\M_{g,n}(L_1,...,L_n)$ can be organized into a disjoint union of finitely many mapping class group orbits; the orbit to which $\gamma$ is assigned is determined by the topology of the surface produced by cutting along $\gamma$. We will handle the contributions of the orbits separately; this is sufficient since the total number of orbits is bounded in terms of $g$ and $n$. First, consider the single orbit containing all non-separating geodesics, and fix a representative $\gamma_0$. Define $f(x)=1_{0\leq x\leq \epsilon}$; as mentioned above, it suffices to bound the following average:
\begin{equation}\label{einitcount}\frac{1}{V_{g,n}(L_1,...,L_n)}\int_{\M_{g,n}(L_1,...,L_n)}\sum_{\gamma\in MCG\cdot\gamma_0}f(l_X(\gamma))dV_{WP}\end{equation}
Since $\sum_{\gamma\in MCG\cdot\gamma_0}f(l(\gamma))$ is a geometric function in the sense defined by Mirzakhani, we may `unfold' the integral in Equation \ref{einitcount} and cut along $\gamma$ to obtain a more tractable integral. This unfolding is a standard technique in Mirzakhani's integration machinery, and is stated as Theorem 8.1 in \cite{MM07}. The result is a bound of
$$\prec\frac{1}{V_{g,n}(L_1,...,L_n)}\int_0^\infty xf(x)V_{g-1,n+2}(x,x,L_1,...,L_n)dx$$
Applying Theorem \ref{volpoly}, this is
$$\prec\frac{1}{V_{g,n}(L_1,...,L_n)}\sum_{|d|\leq 3g-6+n+2}L_1^{2d_1}\cdots L_n^{2d_n}\int_0^\epsilon x^{2d_x+1}dx$$
$$\prec\frac{1}{V_{g,n}(L_1,...,L_n)}\sum_{|d|\leq 3g-6+n+2} \epsilon^{2d_x+2} L_1^{2d_1}\cdots L_n^{2d_n}$$
where $d=(d_x,d_1,...,d_n)$ is the degree tuple for coefficients of the volume polynomial corresponding to the cut surface. Applying $C$-control yields
$$\prec\sum_{0\leq d_x\leq 3g-6+2n+2}\frac{\epsilon^{2d_x+2}L^{6g-12+2n+4-2d_x}}{L^{6g-6+2n}}$$
$$=\sum_{0\leq d_x\leq 3g-6+2n+2}\frac{\epsilon^{2d_x+2}}{L^{2d_x+2}}\prec\frac{\epsilon^2}{L^2}$$

Next, consider the orbit that cuts $X$ into a surface with genus $g_1$ and $k$ boundary components, and a surface with genus $g_2=g-g_1$ and $n-k$ boundary components; this orbit is entirely determined by $g_1$ and $k$. Assume without loss of generality that the first surface has boundary components $\beta_1,...,\beta_k$ and the second surface has boundary components $\beta_{k+1},...,\beta_n$. Fixing a representative $\gamma_0$, we employ the same unfolding technique:
$$\frac{1}{V_{g,n}(L_1,...,L_n)}\int_{\M_{g,n}(L_1,...,L_n)}\sum_{\gamma\in MCG\cdot\gamma_0} f(l_X(\gamma))dV_{WP}$$
$$\prec\frac{1}{V_{g,n}(L_1,...,L_n)}\int_0^\infty xf(x)V_{g_1,k+1}(x,L_1,...,L_k)V_{g_2,n-k+1}(x,L_{k+1},...,L_n)dx$$
$$\prec\frac{1}{V_{g,n}(L_1,...,L_n)}\sum_{\substack{|(d_{x,1},d_1,...,d_k)|\leq3g_1-3+k+1\\|(d_{x,2},d_{k+1},...,d_n)|\leq3g_2-3+n-k+1}}L_1^{2d_1}\cdots L_n^{2d_n}\int_0^\epsilon x^{2(d_{x,1}+d_{x,2})+1}dx$$
$$\prec \frac{1}{V_{g,n}(L_1,...,L_n)}\sum_{\substack{|(d_{x,1},d_1,...,d_k)|\leq3g_1-3+k+1\\|(d_{x,2},d_{k+1},...,d_n)|\leq3g_2-3+n-k+1}}\epsilon^{2(d_{x,1}+d_{x,2})+2}L_1^{2d_1}\cdots L_n^{2d_n}$$
$$\prec \sum_{\substack{0\leq d_{x,1}\leq 3g_1-3+k+1\\ 0\leq d_{x,2}\leq 3g_2-3+n-k+1}}\frac{\epsilon^{2(d_{x,1}+d_{x,2})+2}L^{(6g_1-6+2k+2)+(6g_2-6+2n-2k+2)-2(d_{x,1}+d_{x,2})}}{L^{6g-6+2n}}$$
This simplifies to
$$ \sum_{\substack{0\leq d_{x,1}\leq 3g_1-3+k+1\\ 0\leq d_{x,2}\leq 3g_2-3+n-k+1}}\frac{\epsilon^{2(d_{x,1}+d_{x,2})+2}}{L^{2(d_{x,1}+d_{x,2})+2}}\prec\frac{\epsilon^2}{L^2}\text{$\blacksquare$}$$

\

Denote the boundary components of surfaces in $\M_{g,n}(L_1,...,L_n)$ by $\beta_1,...,\beta_n$. For a constant $\eta>0$, define $U_1(\eta)$ as the set of all surfaces $X\in \M_{g,n}(L_1,...,L_n)$ such that for at least one pants (subsurface homeomorphic to a thrice-punctured sphere) $P$ of $X$:

\begin{enumerate}
\item The boundary components of $P$ are two boundary components $\beta_i$ and $\beta_j$ of $X$ and a non-boundary geodesic $\xi$ of $X$.
\item We have 
\begin{equation}\label{epantsdefect1}|l(\beta_i)+l(\beta_j)-l(\xi)|\leq \eta\end{equation}
\end{enumerate}

Likewise, let $U_2(\eta)$ be the set of all $X\in\M_{g,n}(L_1,...,L_n)$ such that for at least one pants $P$ of $X$:
\begin{enumerate}
\item The boundary components of $P$ are one boundary component $\beta$ of $X$ and two non-boundary geodesics $\xi_1$ and $\xi_2$ of $X$.
\item We have
\begin{equation}\label{epantsdefect2}|l(\beta)-l(\xi_1)-l(\xi_2)|\leq \eta\end{equation}	
\end{enumerate}

\begin{figure}[h]

\centering
\includegraphics[scale=0.15]{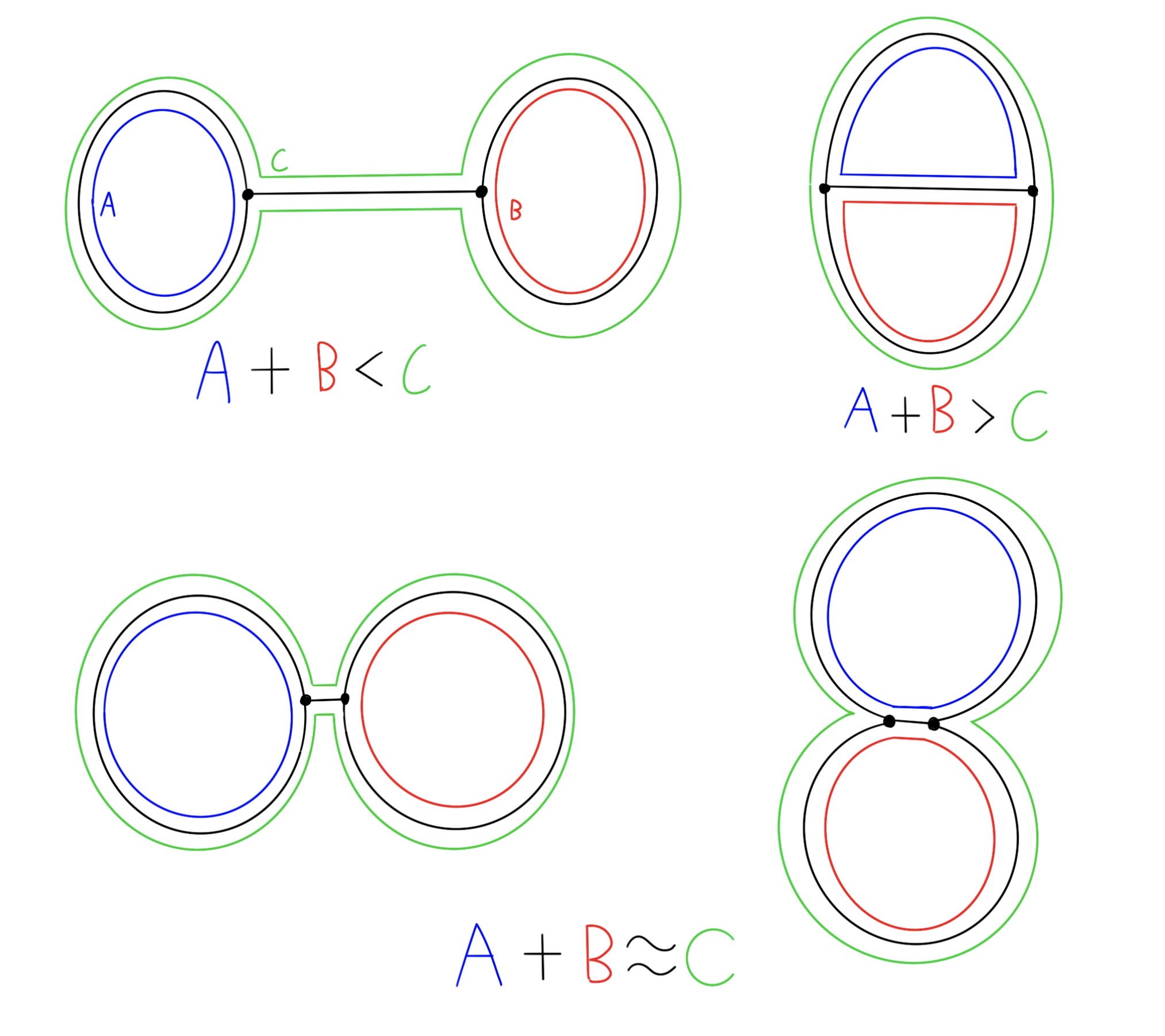}
\caption{A schematic illustration of four pants with inscribed spine graphs, representing the possible shapes when all boundaries of the pants are long. The shape of the pants is determined by the relationship between its boundary lengths. We have a `dumbbell' in the top left position, a `theta' in the top right position, and two approximate `figure eights' in the bottom positions.}

\label{figureeight}

\end{figure}

We will call the quantities $|l(\beta_i)+l(\beta_j)-l(\xi)|$ and $|l(\beta)-l(\xi_1)-l(\xi_2)|$ \emph{pants defects} - when all three lengths are long, they loosely indicate how close a pants is to a `figure-eight shape', as opposed to a `theta' or `dumbbell.' See Figure \ref{figureeight} for an illustration of these possibilities. 

Then:

\begin{theorem}
\label{shortedgevol}

Consider $\M_{g,n}(L_1,...,L_n)$ with $g$ and $n$ fixed, and assume that the $L_i$ are $C$-controlled by $L$. Also fix $\eta>0$ and assume $L$ is sufficiently large in terms of $\eta$. Letting $U_1(\eta)$ and $U_2(\eta)$ be as above,
$$\prob_{\M_{g,n}(L_1,...,L_n)}(U_1(\eta)\cup U_2(\eta))=O\left(\frac{\eta}{L}\right)$$
	
\end{theorem}

\textbf{Proof.} We prove the bound only for $U_2$, since the proof for $U_1$ is similar.

We begin by choosing a single boundary component $\beta$ and considering all pairs $(\gamma_1,\gamma_2)$ of geodesics that form a pants with $\beta$. Since there are $n$ boundary components and our implicit constants are allowed to depend on $n$, it suffices to prove a bound for the contribution of $\beta$. Furthermore, we can assume without loss of generality that $\beta=\beta_1$. The set of such pairs is a disjoint union of finitely many mapping class group orbits, corresponding to the possible topological types of surface(s) obtained by cutting along $\gamma_1$ and $\gamma_2$. Since the number of such orbits is bounded in terms of $g$ and $n$, it suffices to bound the contribution of each orbit separately. We will present only the case where cutting along $\gamma_1$ and $\gamma_2$ decomposes our surface into a pants and a single connected surface with genus $g-1$ and $n+1$ boundary components, as the others are analogous. Fix a representative pair $(\gamma_{1,0},\gamma_{2,0})$. 

As in Lemma \ref{shortloopvol}, we bound the volume of $U_2(\eta)$ by calculating the \emph{average} number of pants on a surface in $\M_{g,n}(L_1,...,L_n)$ with pants defect at most $\eta$. This average is calculated by the following integral, which we must bound:

\begin{equation}\label{einitbound2}\frac{1}{V_{g,n}(L_1,...,L_n)}\int_{\M_{g,n}(L_1,...,L_n)}\sum_{(\gamma_1,\gamma_2)\in MCG\cdot(\gamma_{1,0},\gamma_{2,0})}f(l_{\gamma_1}(X),l_{\gamma_2}(X))dV_{WP}\end{equation}
where $f(x,y)=1_{|x+y-L_1|\leq\eta}$.

Equation \ref{einitbound2} is an integral of a geometric function defined by summing over a mapping class group orbit, so we can unfold to the space $(X,\gamma_1,\gamma_2)$ of hyperbolic surfaces in $\M_{g,n}(L_1,...,L_n)$ with two marked simple closed geodesics of a fixed topological type. Using the Fenchel-Nielsen formula for $V_{WP}$, and calculating the topological type of the surface produced by cutting along $\gamma_1$ and $\gamma_2$, we can bound the expression above by
$$\frac{1}{V_{g,n}(L_1,...,L_n)}\int_0^\infty\int_0^\infty xyf(x,y)V_{g-1,n+1}(x,y,L_2,...,L_n)dxdy$$
\begin{equation}\label{edoubleint}\prec \frac{1}{V_{g,n}(L_1,...,L_n)}\sum_{|d|\leq 3g-6+n+1}L_2^{2d_2}\cdots L_n^{2d_n}\int_{x,y\geq 0,x+y\in[L_1-\eta,L_1+\eta]}x^{2d_x+1}y^{2d_y+1}dxdy\end{equation}
Estimating the integral contained in Equation \ref{edoubleint},
$$\int_{x,y\geq 0,x+y\in[A-\eta,A+\eta]}x^{2d_x+1}y^{2d_y+1}dxdy$$
$$=\int_0^{A+\eta}x^{2d_x+1}\int_{\max(A-\eta-x,0)}^{A+\eta-x}y^{2d_y+1}dydx$$
$$\leq 2\eta\int_0^{A+\eta}x^{2d_x+1}(A+\eta-x)^{2d_y+1}dx$$
Treating $((A+\eta)-x)^{2d_y+1}$ as a binomial and expanding results in
$$2\eta\sum_{\zeta=0}^{2d_y+1}(-1)^{\zeta}\frac{\binom{2d_y+1}{\zeta}}{2d_x+2}(A+\eta)^{2d_y+2d_x+3}$$
$$\prec_{d_x,d_y}\eta(A+\eta)^{2d_y+2d_y+3}$$
Since in our calculation the inputs $d_x$ and $d_y$ are both bounded in terms of $g$ and $n$, we can absorb the implied constants. Returning to Equation \ref{edoubleint}, we have
$$\prec\frac{1}{V_{g,n}(L_1,...,L_n)}\sum_{|d|\leq 3g-6+n+1}L_2^{2d_2}\cdots L_n^{2d_n}\eta(L_1+\eta)^{2d_x+2d_y+3}$$
Performing another binomial expansion (resulting in a new index variable $\beta$) leads to
$$\prec\frac{1}{V_{g,n}(L_1,...,L_n)}\sum_{\substack{|d|\leq 3g-6+n+1\\ 0\leq \beta\leq 2d_x+2d_y+3}}\eta^{\beta+1}L_1^{2d_x+2d_y+3-\beta}L_2^{2d_2}\cdots L_n^{2d_n}$$
We now use $C$-control. For $L$ large, this implies $V_{g,n}(L_1,...,L_n)\asymp_{C,g,n}L^{6g-6+2n}$. Moreover, we can ignore all terms except the leading terms. We obtain.
$$\prec \sum_{0\leq \beta\leq 6g-12+2n+2}\frac{\eta^{\beta+1}L^{6g-12+2n+2+3-\beta}}{L^{6g-6+2n}}$$
$$=\sum_{0\leq \beta\leq 6g-18+2n+8}\frac{\eta^{\beta+1}}{L^{\beta+1}}\prec\frac{\eta}{L}\text{$\blacksquare$}$$

\subsection{Metric Ribbon Graphs}

It turns out that the topological recursion methods employed in the previous section generalize naturally to $MRG_{g,n}(L_1,...,L_n)$. 
\begin{lemma}
\label{shortloopgraphvol}

Consider $MRG_{g,n}(L_1,...,L_n)$ with $g$ and $n$ fixed, and choose a constant $\epsilon>0$. Assume that the $L_i$ are $C$-controlled by $L$. Let $U(\epsilon)$ be the subset of $MRG_{g,n}(L_1,...,L_n)$ containing metric ribbon graphs $\Gamma$ with at least one closed geodesic with length at most $\epsilon$. Then
$$\prob_{MRG_{g,n}(L_1,...,L_n)}(U(\epsilon))\prec O\left(\frac{\epsilon^2}{L^2}\right)$$
	
\end{lemma}

The ribbon-graph analogue of Theorem \ref{shortedgevol} is also easily obtainable, but is not needed for our results.

\textbf{Proof Sketch.} This follows from a computation essentially identical to that performed in Lemma \ref{shortloopvol}, once the following prerequisites are known:
\begin{enumerate}
\item $TRG_{g,n}(L_1,...,L_n)$ possesses a version of Fenchel-Nielsen coordinates. Just as for $\T_{g,n}(L_1,...,L_n)$, a choice of Fenchel-Nielsen coordinates is highly non-canonical, but in any such coordinate system a direct analogue of Wolpert's magic formula holds.
\item There is a direct analogue of topological recursion for geometric functions integrated over $MRG_{g,n}(L_1,...,L_n)$. This allows for the 'unfolding step' to be performed. 
\item $\vol(MRG_{g,n}(L_1,...,L_n))$ satisfies a polynomial formula identical to \ref{evolpoly}, except with different coefficients. 
\end{enumerate}

For rigorous statements and proofs of these facts, see \cite{ABCGLW21}.

\section{Geometry Control}

\subsection{Classifying Surfaces}

We will use the results of the previous sections to subdivide $\M_{g,n}(L_1,...,L_n)$ and $MRG_{g,n}(L_1,...,L_n)$. 

\begin{definition} 
\label{surfaceclassification}
We define the following subsets of $\M_{g,n}(L_1,...,L_n)$.

\begin{enumerate}
\item $\M_I$ is the locus of all surfaces $X$ such that $\Phi(X)$ is not trivalent. We say surfaces $X\in \M_I$ are \emph{irrelevant}.
\item $\M_{NG}(\epsilon)$ is the locus of all surfaces $X$ such that $\Phi(X)$ is trivalent and $\Phi(X)$ contains a geodesic of length at most $\epsilon$. We say such surfaces are \emph{not good}.
\item $\M_{G}(\epsilon,\eta)$ is the locus of all surfaces $X$ such that $\Phi(X)$ is trivalent, all geodesics of $\Phi(X)$ have length at least $\epsilon$, but some edge of $\Phi(X)$ has length less than $\eta$. We say such surfaces are $(\epsilon,\eta)$\emph{-good}.
\item $\M_{GG}(\epsilon,\eta)$ is the locus of all surfaces $X$ such that $\Phi(X)$ is trivalent, all geodesics of $\Phi(X)$ have length at least $\epsilon$, and all edges of $\Phi(X)$ have length at least $\eta$. We say such surfaces are $(\epsilon,\eta)$\emph{-great}. We also say that a surface $X$ is $\eta$\emph{-great} if it is $(\eta,\eta)$\emph{-great}, i.e. all edges and all geodesics of $\Phi(X)$ have length at least $\eta$ - however, note that by the upcoming Lemma \ref{pathlengthbound}, the bound on edges implies the bound on geodesics.

\end{enumerate}
\end{definition}

We will also need to define absolute constants $A$ and $B$ in connection to these loci. We put the following conditions on $A$ and $B$, and on all $\epsilon$ and $\eta$ used to partition $\M_{g,n}(L_1,...,L_n)$ as above:

\begin{enumerate}
\item All $\epsilon$ satisfy $B\leq\epsilon$, and all $\eta$ satisfy $A\leq\eta$.
\item $A$ and $B$ are both sufficiently large in terms of $g$ and $n$.
\item 	$B$ is a sufficiently large multiple of $A$ in terms of $g$ and $n$, and $\epsilon$ is a sufficiently large multiple of $\eta$ in terms of $g$ and $n$.
\end{enumerate}

Note that the second and third conditions may always be simultaneously satisfied.

\subsection{Corridors} For any $X\in\M_{g,n}(L_1,...,L_n)$ with spine graph $S_X$, the ribs corresponding to $S_X$ cut $X$ into a finite number of geodesic hexagons we will call corridors. Each corridor has four right angles and two angles that may vary, and its edges consist of two boundary arcs of $X$ and four ribs. 

We first need some basic geometric properties of corridors, all of which follow from elementary trigonometric arguments.

\begin{lemma} 
\label{wedgeprops}

Consider a corridor $W$ embedded in a surface $X$ as above. Then:
\begin{enumerate}
\item The six edges of $W$ are geodesic arcs, and all angles between the boundary arcs and the ribs are right angles.
\item $W$ is mirror-symmetric around the spine edge it contains.
\item If $\theta$ and $\theta'$ are the two non-right angles in $W$, then $\theta$ and $\theta'$ cannot both be at least $\pi$.
\item The shortest paths between one pair of ribs and the other are given by the boundary arcs.
\end{enumerate}
	
\end{lemma}

Next, we need to know that the spine graph cannot be too complex. This result follows from a straightforward Euler characteristic computation.

\begin{lemma}
\label{graphnum}

Fix $X\in\M_{g,n}(L_1,...,L_n)$ with spine graph $S_X$, and assume $S_X$ is trivalent. Then $S_X$ has $6g-6+3n$ edges and $4g-4+2n$ vertices. In particular, both quantities are $O(1)$.
	
\end{lemma}

\

Recall that each corridor $W$ has a unique associated intercostal $\iota(W)$, the unique minimal-length arc homotopic to a path in $W$ between the boundary arcs rel $\partial X$. Again, this result is elementary.

\begin{lemma}
\label{intercostal}

Fix $X\in\M_{g,n}(L_1,...,L_n)$ such that its spine $S_X$ is trivalent, and set $\Gamma=\Phi(X)$. For each corridor $W\subset X$, define its intercostal $\iota(W)$ as above. Then:
\begin{enumerate}
\item Each $\iota(W)$ is a simple geodesic arc that lies in the interior of $X$ except at its endpoints.
\item Cutting along the union of the $\iota(W)$ decomposes $X$ into hexagons. 
\item There is a natural correspondence between these hexagons and vertices of $\Gamma$,  and the hexagons are connected precisely as indicated by the edges of $X$.
\end{enumerate}
	
\end{lemma}

\

We will refer to the hexagons obtained by cutting along intercostals as \emph{sectors}.

\subsection{Geometric Correspondence}

Fix $X\in\M_{g,n}(L_1,...,L_n)$ with spine $S_X$ homeomorphic to the metric ribbon graph $\Gamma=\Phi(X)$. By the spine construction, $X$ deformation retracts onto $S_X$, and so there is a natural homotopy equivalence between $X$ and $\Gamma$. Since both $X$ and $\Gamma$ have the property that every homotopy class of closed path contains a unique geodesic representative, this homotopy equivalence allows us to define a correspondence between geodesics on $X$ and geodesics on $\Gamma$. We will usually use a tilde to represent this correspondence, so the geodesic $\gamma$ on $X$ corresponds to a geodesic $\tilde{\gamma}$ on $\Gamma$.

For a geodesic $\gamma$ on a metric graph $\Gamma$, we define its combinatorial length $c(\gamma)$ to be the number of edges $\gamma$ traverses, counted with multiplicity. For $X\in \M_{g,n}(L_1,...,L_n)$ and $\gamma$ a geodesic on $X$, we define $c(\gamma)$ to be the number of corridors $\gamma$ traverses, again counted with multiplicity. If $\tilde{\gamma}$ on $\Phi(X)$ corresponds to $\gamma$, then $c(\gamma)=c(\tilde{\gamma})$.

Correspondence of geodesics extends to intersections, but we first need to deal with subtleties related to how intersections on $\Gamma$ are defined. 

\begin{figure}

\centering
\includegraphics[scale=0.15]{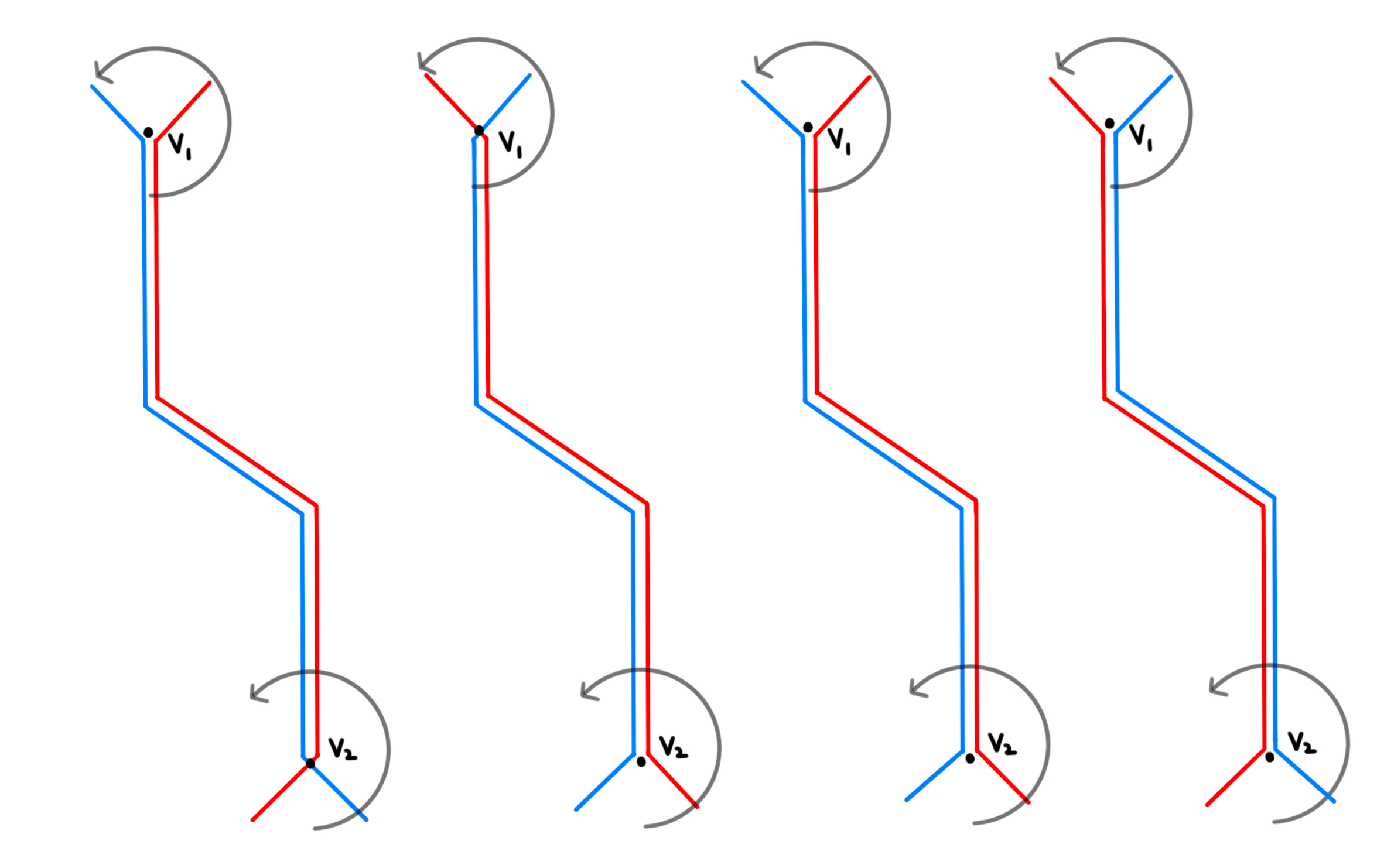}
\caption{Four possible ways two arcs may overlap on a trivalent ribbon graph. We consider the left two to represent one intersection each, and the right two to not represent intersections. Note that the cases may be distinguished by examining the cyclic edge orderings at the vertices $v_1$ and $v_2$.}

\label{intersectionongraph}
	
\end{figure}

\begin{definition} 
\label{graphintersectiondef}
Given two closed geodesics $\tilde{\gamma}$ and $\tilde{\eta}$ on a \emph{trivalent} ribbon graph $\Gamma$, we say $\tilde{\gamma}$ and $\tilde{\eta}$ \emph{intersect} at a pair of subarcs $\tilde{\gamma}'\subset\tilde{\gamma}$ and $\tilde{\eta}'\subset\tilde{\eta}$ if the following properties hold:
\begin{enumerate}
\item Both $\tilde{\gamma}'$ and $\tilde{\eta}'$ consist of a series of edges with a half-edge attached at each end. The edges will be referred to as the \emph{bodies}, while the half-edges will be referred to as the \emph{tails}.

\item As arcs in $\Gamma$, $\tilde{\gamma}'$ and $\tilde{\eta}'$ have equal bodies ending at vertices $v_1$ and $v_2$.
\item $\tilde{\gamma}'$ and $\tilde{\eta}'$ `diverge' in the sense that, at either endpoint $v_i$, $\tilde{\gamma}'$ and $\tilde{\eta}'$ have distinct tails.
\item By the previous property, at $v_1$ we may label the three incident edges as the `body' edge shared by the bodies of $\tilde{\gamma}'$ and $\tilde{\eta}'$ , the tail edge of $\tilde{\gamma}$, and the tail edge of $\tilde{\eta}$, and make the same labeling at $v_2$. Then the cyclic ordering induced on these three labels at both $v_1$ and $v_2$ agree, as in the left half of Figure \ref{intersectionongraph}.
\end{enumerate}
	
\end{definition}

Conceptually, the fourth condition above requires that $\tilde{\eta}$ passes from one side of $\tilde{\gamma}$ to the other, or vice-versa. 
\begin{definition}

Given a closed geodesic $\tilde{\gamma}$ on a trivalent ribbon graph $\Gamma$, we say $\tilde{\gamma}$ is \emph{simple} if there exists no pair of distinct subarcs (i.e. subarcs coming from distinct domains under some parameterization) $\tilde{\gamma}_1,\tilde{\gamma}_2$ of $\tilde{\gamma}$ satisfying the four properties given in Definition \ref{graphintersectiondef}.
	
\end{definition}

\begin{remark}

When working with intersections of graph geodesics, instead of referring to a pair of subarcs $\tilde{\gamma}'$ and $\tilde{\eta}
'$ as above, we will typically refer to a single `object' $\tilde{p}=(\tilde{\gamma}',\tilde{\eta}')$. We will also write $\tilde{p}\in\tilde{\gamma}\cap\tilde{\eta}$, thinking of $\tilde{\gamma}\cap\tilde{\eta}$ not as a literal intersection set but as a label for the set of all such $\tilde{p}$. With this setup, $i(\tilde{\gamma},\tilde{\eta})$ is defined as $|\tilde{\gamma}\cap\tilde{\eta}|$.

\end{remark}

We record the following result, which essentially follows from the definitions and the fact that a surface deformation retracts onto its spine.

\begin{lemma}
\label{intersectioncorr}

Fix $X\in\M_{g,n}(L_1,...,L_n)$ with $\Gamma=\Phi(X)$, and assume $\Gamma$ is trivalent. Given simple closed geodesics $\gamma$ and $\eta$ on $X$, there is a natural correspondence between points $p\in \gamma\cap \eta$ and arcs $\tilde{p}\in\tilde{\gamma}\cap\tilde{\eta}$. 
	
\end{lemma}

\

An analogous argument works for self-intersections, leading to:

\begin{lemma}
\label{simplecorr}

With the setup above and given a closed geodesic $\gamma$ on $X$, $\gamma$ is simple if and only if $\tilde{\gamma}$ on $\Gamma$ is simple.
	
\end{lemma}

Lastly, we state a simple result about how the geodesic correspondence relates to lengths.

\begin{lemma}
\label{pathlengthbound}

Let $X\in\M_{g,n}(L_1,...,L_n)$ be a surface with spine $S_X$, and let $\gamma$ be a geodesic on $X$. Let $\tilde{\gamma}$ be the corresponding geodesic on $\Phi(X)$ as above. Then
$$l_X(\gamma)\geq l_\Gamma(\tilde{\gamma})$$
	
\end{lemma}

\textbf{Proof.} The geodesic $\gamma$ may be subdivided into a series of geodesic arcs through corridors, and each such arc corresponds to an edge of $\Phi(X)$ traversed by $\tilde{\gamma}$. Therefore, it suffices to prove that any arc through a corridor in $X$ is at least as long as the edge in $\Phi(X)$ corresponding to that corridor. Since the length of an edge in $\Phi(X)$ is by definition the length of either boundary arc in the corresponding corridor, we may immediately apply item 4 of Lemma \ref{wedgeprops}.

\begin{corollary}
\label{goodbridge}

If $X$ is $(B,A)$-good, then any geodesic $\gamma$ of $X$ has length $l(\gamma)\geq B$.
	
\end{corollary}

\begin{corollary}

For any surface $X$, $\delta_X\leq\delta_{\Phi(X)}$.
	
\end{corollary}

Combining Corollary \ref{goodbridge} with Lemma \ref{shortloopvol} yields:

\begin{corollary}
\label{shortgeodcontrol}

Fix $\epsilon\geq B$ and assume $L_1,...,L_n$ are $C$-controlled by a parameter $L$. Then for $L$ sufficiently large in terms of $\epsilon$,
$$\prob_{\M_{g,n}(L_1,...,L_n)}(\M_{NG}(\epsilon))=O\left(\frac{\epsilon^2}{L^2}\right)$$
	
\end{corollary}

\subsection{Rib Control}

Recall that all $X\in\M_{g,n}(L_1,...,L_n)$ have the same area $\alpha_{g,n}=(4g-4+2n)\pi=O_{g,n}(1)$.

\begin{lemma}
\label{ribbound}

As long as $B$ is sufficiently large in terms of $g$ and $n$, there exists a constant $R=O(1)$ such that all $(B,A)$-good surfaces have maximal rib length bounded above by $R$.
	
\end{lemma}

\textbf{Proof.} Choose a constant $r_0>0$ such that a disk in $\H$ with radius $r_0$ has area $\alpha_{g,n}+1$. Assume $B\geq r_0$; this is allowed since $r_0$ depends only on $g$ and $n$. For any $X\in \M_{g,n}(L_1,...,L_n)$ with spine graph $S_X$, and any $v\in S_X$, assume the rib lengths at $v$ are greater than $r_0$. Therefore, if $B(r_0)$ is the ball of radius $r_0$ in $X$ centered at $v$, $B(r_0)$ does not intersect the boundary. Since furthermore $X$ has no geodesics of length $r_0$ or less, $B(r_0)$ in fact embeds into $X$, and the area of $X$ is too large. $\blacksquare$

\

As an immediate consequence:

\begin{corollary}
\label{ribcontrol}

Let $X\in\M_{g,n}(L_1,...,L_n)$ be $(B,A)$-good, with $B$ taken sufficiently large in terms of $g$ and $n$. Set $\Gamma=\Phi(X)$, let $\gamma$ be a closed geodesic on $X$, and let $\tilde{\gamma}$ be the corresponding geodesic on $\Gamma$. Then
$$l_\Gamma(\tilde{\gamma})\leq l_X(\gamma)\leq l_\Gamma(\tilde{\gamma})+O(c(\gamma))$$
where $c(\gamma)$ is the combinatorial length of $\gamma$.
\end{corollary}

\textbf{Proof.} Apply the same subdivision scheme as in Lemma \ref{pathlengthbound}. Any geodesic arc passing through a corridor has length bounded above by the length of either boundary of that corridor, plus the sum of the lengths of the ribs of that corridor. But by Lemma \ref{ribbound}, the lengths of these ribs add to $O(1)$. $\blacksquare$

\begin{lemma}
\label{shortedgecontrol}

Fix $g$, $n$, and $C$, and consider $\M_{g,n}(L_1,...,L_n)$ with the $L_i$ $C$-controlled by a parameter $L$. Assume that $A$ and $B$ are sufficiently large in terms of $g$, $n$, and $C$, and choose parameters $\epsilon$ and $\eta$ satisfying $B\leq\epsilon$, $A\leq\eta$. Then
$$\prob_{\M_{g,n}(L_1,...,L_n)}(\M_G(\epsilon,\eta))=O\left(\frac{\eta}{L}\right)$$
	
\end{lemma}

\textbf{Proof.} First replace $\epsilon$ with $B$; this can only increase the volume.

\begin{figure}[h]

\centering
\includegraphics[scale=0.18]{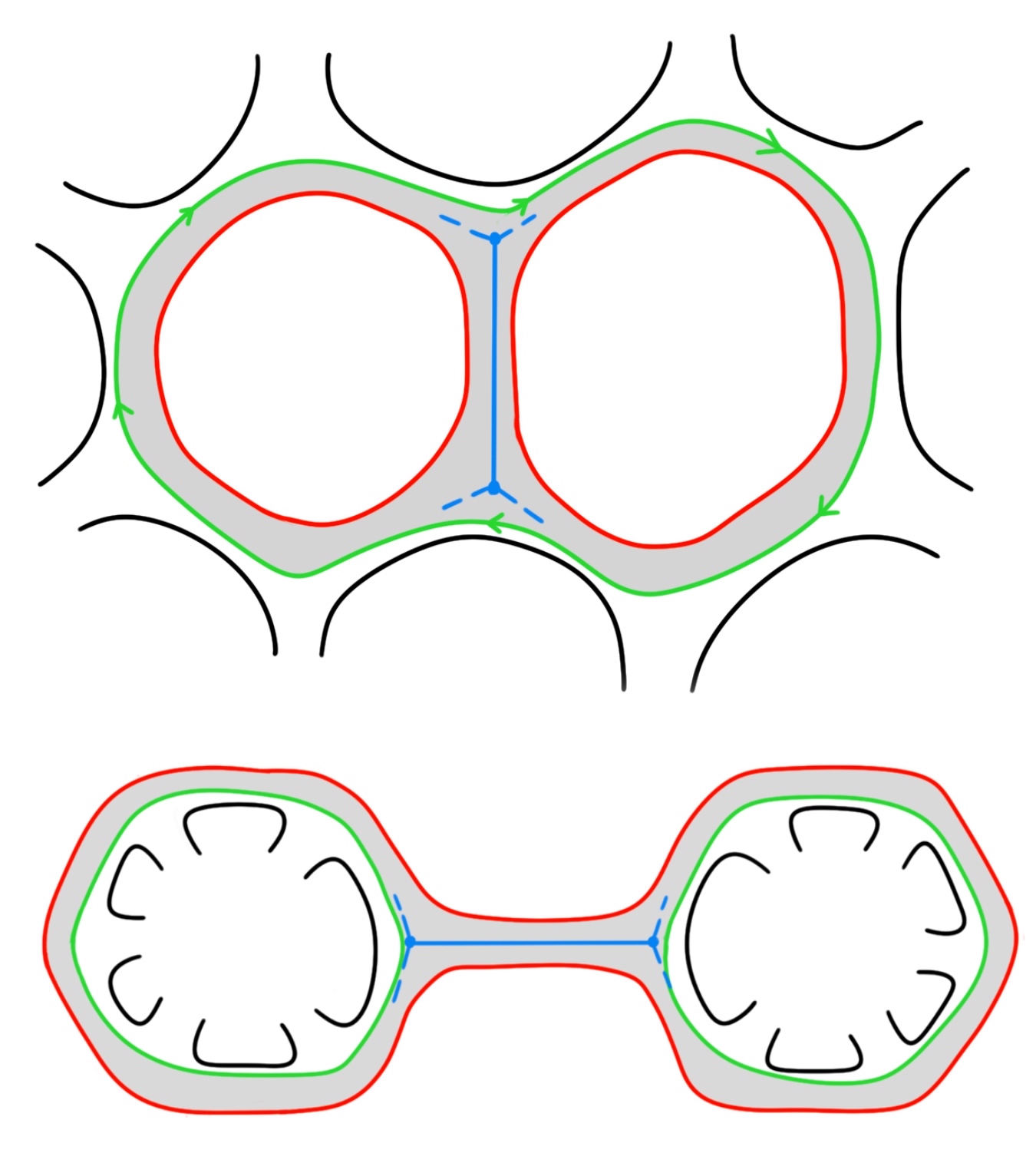}
\caption{Two possible cases of constructing a pair of pants from a given edge (blue). In the first case, the pants has two boundaries (red) that are also boundaries of the surface, and one boundary that is a geodesic on the surface (green). In the second case, the pants has one boundary that is also a boundary of the surface (red) and two boundaries that are geodesics on the surface (green).}

\label{buildpants}

\end{figure}

The key idea is that a short edge always corresponds to a small pants defect. Specifically, for each spine edge $e\in S_X\subset X$, we construct a pants $P_e$ via one of two templates, depending on whether the boundary arcs on either side of $e$ lie in the same or different boundary components (see Figure \ref{buildpants}). This pants construction is due to Do \cite{ND10}.

Corollary \ref{ribcontrol} now allows us to verify from this figure that an edge of length $l$ corresponds to a pants defect of length $l+O(1)$, in the sense of equations \ref{epantsdefect1} and \ref{epantsdefect2} (equation \ref{epantsdefect1} corresponds to the top part of figure \ref{buildpants}, and equation \ref{epantsdefect2} corresponds to the bottom part of figure \ref{buildpants}). For $A$ sufficiently large in terms of this $O(1)$ error, we may apply Theorem \ref{shortedgevol} to conclude the desired result. $\blacksquare$

\subsection{Intercostal Bounds}

We now focus on $\eta$-great surfaces. Our goal will be to show that if $X\in \M_{g,n}(L_1,...,L_n)$ is $\eta$-great and $L$ is sufficiently large, then $\Phi(X)$ is roughly `geometically equivalent' to $X$. First, corridors must contain their intercostals:

\begin{lemma}
\label{containintercostal}

Let $X\in\M_{g,n}(L_1,...,L_n)$ be $(B,A)$-good, and let $W$ be a corridor of $X$. As long as the lengths of the boundary arcs of $W$ are at least $A$ for $A$ sufficiently large in terms of $g$ and $n$, then $W$ contains its intercostal. Equivalently, both non-right angles of $W$ are less than $\pi$.
	
\end{lemma}

\begin{figure}[h]

\centering
\includegraphics[scale=0.1]{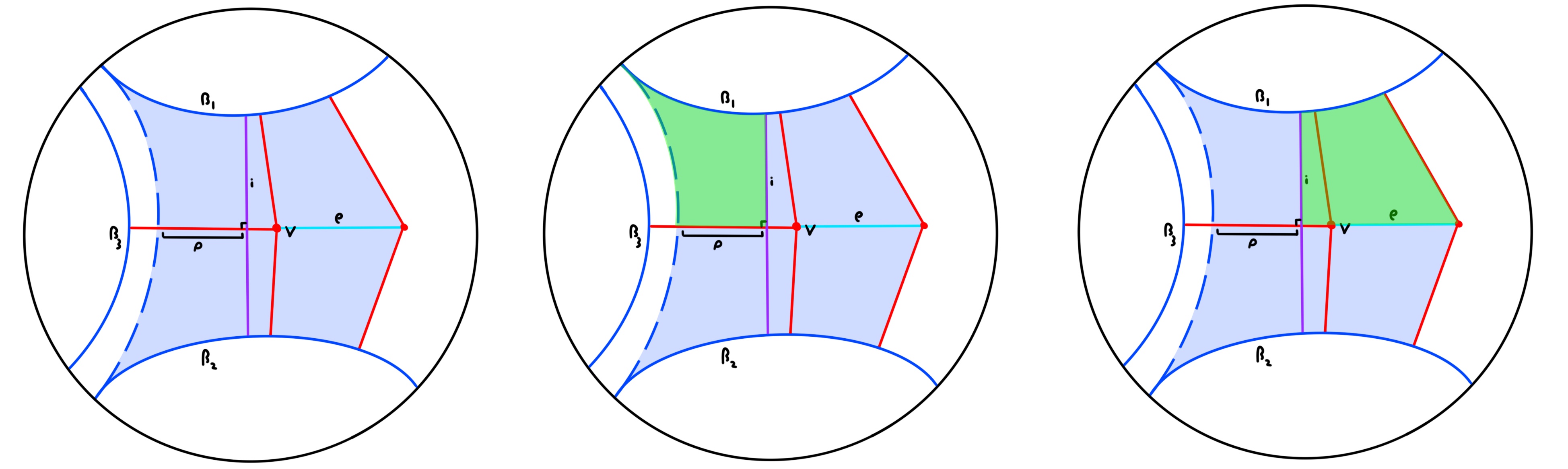}
\caption{Left: the construction presented in Lemma \ref{containintercostal}. Center: the same construction, with the first quadrilateral mentioned in the proof of Lemma \ref{containintercostal} shaded. Right: the same construction, with the second quadrilateral mentioned in the proof of Lemma \ref{containintercostal} shaded.}

\label{intercostalfigure}

\end{figure}

\textbf{Proof.} Take the two corridors meeting $W$ at its obtuse angle and glue them to $W$, then embed the resulting geodesic polygon in $\H$ (see Figure \ref{intercostalfigure}). Let $v$ be the spine vertex corresponding to the ribs along which this gluing occurs. We obtain a nine-sided polygon with three sides corresponding to boundary arcs in $X$; extend these arcs into disjoint geodesics $\beta_1$, $\beta_2$, and $\beta_3$, with $\beta_1$ and $\beta_2$ incident to $W$. There exists some intercostal arc $\iota$ between $\beta_1$ and $\beta_2$.

The two symmetric regions bounded in $\H$ by $\beta_1$, $\beta_2$, and $\iota$ are convex, so if $\beta_3$ has both endpoints in one such region, it cannot intersect $\iota$. Furthermore, the distances from $\beta_3$ to $\iota$ and from $\beta_3$ to $v$ can only decrease if $\beta_3$ is moved so that it meets $\beta_1$ and $\beta_2$ on $\partial\H$. After this move is made, we may form a quadrilateral with one ideal angle and one side that is exactly half of $\iota$ (see Figure \ref{intercostalfigure}). By Corollary \ref{ribbound} applied to the other finite-length side of this quadrilateral, the length $\rho$ from $\beta_3$ to $\iota$ is $O(1)$, so by the first formula of Lemma \ref{quadrilateral} we have $\iota=\Omega(1)$. Now consider the quadrilateral indicated in the third part of Figure $\ref{intercostalfigure}$. By the second formula of Lemma \ref{quadrilateral} and since $l(\iota)=\Omega(1)$, the ribs of $W$ attached to the acute angles have length at least $\Omega(A)$, which is a contradiction for $A$ sufficiently large in terms of $g$ and $n$ by Lemma \ref{ribbound}. $\blacksquare$

\

From here onwards, we will assume that $A$ is sufficiently large such that $\eta\geq A$ implies all corridors contain their intercostals.

\begin{lemma}
\label{intercostalcontrol}

Let $X\in\M_{g,n}(L_1,...,L_n)$ be $\eta$-great, and let $W$ be a corridor of $X$ with intercostal $\iota$. Then
$$l(\iota)=O(e^{-\eta/2})$$
	
\end{lemma}

\textbf{Proof.} The intercostal $\iota$ and the spine edge $e$ contained in $W$ split $W$ into two pairs of symmetric quadrilaterals. Choose a quadrilateral such that if $a$ is the corresponding boundary subarc, $l(a)\geq \eta/2$. Letting $\theta$ be the sole non-right angle of the quadrilateral, since $\theta$ is acute, the first formula of Lemma \ref{quadrilateral} tells us that
$$0\leq \sinh(l(a))\sinh(l(\iota)/2)\leq 1$$
Applying $\sinh(x)=\frac{1}{2}e^x+O(e^{-x})$ and $\sinh^{-1}(x)=x+O(x^2)$ yields the desired bound. $\blacksquare$

\

\begin{lemma}
\label{intercostalbalance}

Let $X\in\M_{g,n}(L_1,...,L_n)$ be $\eta$-great, and let $W$ be a corridor of $X$ with intercostal $\iota$ such that $\iota$ is contained in $W$.  The intercostal $\iota$ cuts each boundary arc into respective subarcs $\beta_1$ and $\beta_2$. Then for $i=1,2$,
$$l(\beta_i)=\Omega(\eta)$$
	
\end{lemma}

\textbf{Proof.} The intercostal $\iota$ and the spine edge $e$ split $W$ into two pairs of symmetric quadrilaterals; choose any one such quadrilateral. One of its sides has length $l(\iota)/2$, and another one of its sides is a rib with length $O(1)$. The first formula of \ref{quadrilateral} and the previous lemma yield the desired asymptotic. $\blacksquare$

\

\subsection{Critical Exponent Convergence for Great Surfaces}

\begin{lemma}
\label{greatlengthbound}

Fix $X\in\M_{g,n}(L_1,...,L_n)$ with spine graph $S_X$, and assume that $X$ is $\eta$-great. Choose a geodesic $\gamma$ on $X$, and let $\tilde{\gamma}$ be the corresponding geodesic on $\Phi(X)$. Then
$$l(\tilde{\gamma})\leq l(\gamma)\leq l(\tilde{\gamma})+O(e^{-\eta/2}c(\gamma))$$
	
\end{lemma}

\textbf{Proof.} Assume that $\tilde{\gamma}$ traverses edges $e_{i_1},...,e_{i_{c(\gamma)}}$ corresponding to corridors $W_{i_1},...,W_{i_{c(\gamma)}}$. Since $X$ is great, each corridor contains its intercostal, so we may cut each corridor into half-corridors and attach the second component of each $W_{i_k}$ (in order of traversal by $\gamma$) to the first component of $W_{i_{k+1}}$ to form polygons $V_{i_k}$, wrapping around at the end of the list. $\gamma$ traverses these $V_{i_k}$ in order, forming a geodesic arc between the two intercostals in each $V_{i_k}$. The boundary arc between these two intercostals is the shortest path between them, and adding up the lengths of these arcs gives $l(\tilde{\gamma})$. On the other hand, any arc between the intercostals must have length at most the length of the boundary arc plus the lengths of the two intercostals, and adding up these bounds using Lemma \ref{intercostalcontrol} gives $l(\tilde{\gamma})+O(e^{-\eta/2}c(\gamma))$. $\blacksquare$

\

Furthermore, the combinatorial length $c(\gamma)$ of a geodesic $\gamma$ can be controlled in terms of its length.

\begin{lemma}
\label{wordcontrol}

Fix $X\in\M_{g,n}(L_1,...,L_n)$ with the $L_i$ $C$-controlled by $L$, and assume $X$ is $\eta$-great. Also fix a geodesic $\gamma$ on $X$. Then
$$\frac{l(\gamma)}{L}\prec c(\gamma)\prec \frac{l(\gamma)}{\eta}$$
	
\end{lemma}

\textbf{Proof.} Every edge of $\Phi(X)$ (and hence every corridor boundary arc of $X$) has length at least $\eta$ and at most $O(L)$, since every edge is part of a face. Corollary \ref{ribcontrol} now allows us to obtain the desired bounds. $\blacksquare$

\

We can now show that $\Phi$ roughly preserves critical exponents of great surfaces.

\begin{lemma}
\label{critexpcontrol}

Fix $X\in\M_{g,n}(L_1,...,L_n)$ and set $\Gamma=\Phi(X)$. Let $\delta_X$ be the critical exponent of $X$, and let $\delta_\Gamma$ be the critical exponent of $\Gamma$. If $X$ is $\eta$-great, then
$$\delta_X\leq\delta_\Gamma\leq\delta_X\left(1+O\left(\frac{e^{-\eta/2}}{\eta}\right)\right)$$
	
\end{lemma}

\textbf{Proof.} Let $\pi_X(x)$ count the number of closed geodesics on $X$ with length at most $x$, and let $\pi_\Gamma(x)$ count the number of closed geodesics on $\Gamma$ with length at most $x$.

We have already established a bijection $\gamma\leftrightarrow\tilde{\gamma}$ between geodesics on $X$ and geodesics on $\Gamma$, and we know by Lemma \ref{pathlengthbound} that $l(\tilde{\gamma})\leq l(\gamma)$. Therefore $\pi_X(x)\leq\pi_{\Gamma}(x)$ for all $x$, and $\delta_X\leq\delta_\Gamma$. On the other hand, by Lemmas \ref{greatlengthbound} and \ref{wordcontrol}, for any $x>0$, every $\tilde{\gamma}$ counted by $\pi_{\Gamma}(x)$ must have its corresponding $\gamma$ counted by
$$\pi_X\left(l(\tilde{\gamma})+O(e^{-\eta/2}c(\gamma))\right)\leq\pi_X\left(x\left(1+O\left(\frac{e^{-\eta/2}}{\eta}\right)\right)\right)$$
Therefore
$$\pi_\Gamma(x)\leq \pi_X\left(x\left(1+O\left(\frac{e^{-\eta/2}}{\eta}\right)\right)\right)$$
We may now apply Theorem \ref{critexpdef} and its metric-graph counterpart, send $x\to\infty$ in the resulting formulas for $\pi_\Gamma(x)$ and $\pi_X(x)$, and obtain the desired inequalities. $\blacksquare$

\section{Quantitative Convergence: Forms}

In this section and afterwards, we will often work with quantities such as tangent vectors, cotangent vectors, and differential forms that are technically only defined at a point. Since the point of definition will generally be clear, we will often leave off referencing the point in notation - for example, given $X\in \M_{g,n}(L_1,...,L_n)$ and a smooth function $f:\M_{g,n}(L_1,...,L_n)\to\reals$, we will write $df$ for the corresponding cotangent vector in $T^*_X\M_{g,n}(L_1,...,L_n)$, rather than $df|_X$. 

Also, all local coordinates will be smooth. To justify this assertion, we will use a few results from Appendix \ref{nontrivalent}. Recall that $\M_I$ is the locus of $\M_{g,n}(L_1,...,L_n)$ consisting of surfaces $X$ such that $\Phi(X)$ is not trivalent. Briefly, both $\M_{g,n}(L_1,...,L_n)$ and $MRG_{g,n}(L_1,...,L_n)$ have cell complex structures, such that $\M_I$ is the complement of the union of top-dimensional cells in $\M_{g,n}(L_1,...,L_n)$ and $\Phi(\M_I)$ is the complement of the union of top-dimensional cells in $MRG_{g,n}(L_1,...,L_n)$. Moreover, the $dV_{WP}$-volume of $\M_I$ and the $dV_K$-volume of $\Phi(\M_I)$ are both zero, and $\Phi$ restricted to the union of top-dimensional cells of $\M_{g,n}(L_1,...,L_n)$ is a diffeomorphism onto the union of top-dimensional cells of $MRG_{g,n}(L_1,...,L_n)$ (with respect to Fenchel-Nielsen coordinates on both spaces).

\subsection{The Weil-Petersson Form}

To work with the Weil-Petersson form, we will use Wolpert's cosine formula. Let $X\in\M_{g,n}(L_1,...,L_n)$ be a surface, and let $\gamma$ be an oriented simple closed multi-geodesic on $X$. Then one may define a twist function $\tau:\reals\to\M_{g,n}(L_1,...,L_n)$, where $\tau(t)$ is the surface obtained by cutting along $\gamma$, twisting by $t$, and gluing the surface back together (see the remark below for a comment on how we define positive twists). The tangent vector $\frac{\partial}{\partial \tau}\in T_X$ is then defined as $\frac{\partial}{\partial t}\tau(t)|_{t=0}$. We stretch notation by writing $\frac{\partial}{\partial\tau}$ instead of $\frac{\partial}{\partial\tau}|_X$, but for our purposes this specificity will not be needed. Crucially, we can define twist derivatives without the need to appeal to a coordinate system where the twists appear as coordinates. 

\begin{remark}
We define positive twists/earthquakes to be such that, if a geodesic is oriented vertically, the right `side' is shifted up, or equivalently the left side is shifted down. See Figure \ref{twistchange} for an illustration of this convention.
\end{remark}

Given two distinct geodesics $\gamma_1$ and $\gamma_2$ on a surface $X\in \M_{g,n}(L_1,...,L_n)$, they intersect at most finitely many times. At each intersection point $p$, we measure the intersection angle $\theta_p$ counterclockwise from $\gamma_1$ to $\gamma_2$. We then define the \emph{Wolpert number} $w(\gamma_1,\gamma_2)$ as

$$w(\gamma_1,\gamma_2)=\sum_{p\in\gamma_1\cap\gamma_2}\cos(\theta_p)$$

This quantity is clearly bounded in absolute value by $i(\gamma_1,\gamma_2)=|\gamma_1\cap\gamma_2|$, the usual intersection number.

\begin{lemma}[Wolpert, \cite{SW83}, Lem. 4.2]
\label{wolpertformula}

Fix $X\in\M_{g,n}(L_1,...,L_n)$, and choose two simple closed geodesics $\gamma_1$ and $\gamma_2$ on $X$. Let $\frac{\partial}{\partial\tau_2}$ be the twist derivative corresponding to $\gamma_2$. Then
$$\frac{\partial}{\partial \tau_2}l(\gamma_1)=w(\gamma_1,\gamma_2)$$
	
\end{lemma}

Wolpert also proved the following length-twist duality result.

\begin{lemma}[Wolpert, \cite{SW83}, Lem. 4.1]
\label{duality}

Fix $X\in\M_{g,n}(L_1,...,L_n)$, and let $\gamma$ be a simple closed geodesic on $X$. Let $l$ be the length of $\gamma$, and let $\frac{\partial}{\partial \tau}$ be the corresponding twist derivative. Then
$$dl(\cdot)=\omega_{WP}\left(\frac{\partial}{\partial\tau},\cdot\right)$$
	
\end{lemma}

A useful consequence of this lemma, due to Wolpert, is as follows. Since we will copy it more or less verbatim for the case of metric ribbon graphs, we provide the proof for ease of referral.

\begin{lemma}[Wolpert, \cite{SWP83}, Lem. 4.5]
\label{symplectic}

Fix $X\in\M_{g,n}(L_1,...,L_n)$ such that the spine $S_X$ of $X$ is trivalent, and choose simple closed geodesics $\gamma_1,...,\gamma_{6g-6+2n}$ on $X$. For each $\gamma_i$, define $l_i=l(\gamma_i)$, and let $\frac{\partial}{\partial\tau_i}$ be the twist derivative corresponding to $\gamma_i$. Define a $(6g-6+2n)\times (6g-6+2n)$ skew-symmetric matrix $W$ by $W_{ij}=w(\gamma_i,\gamma_j)$. Then $l_1,...,l_{6g-6+2n}$ form local coordinates for $\M_{g,n}(L_1,...,L_n)$ in a neighborhood of $X$ if and only if $W$ is invertible. Furthermore, if $l_1,...,l_{6g-6+2n}$ do form local coordinates, then
$$\omega_{WP}\left(\frac{\partial}{\partial l_i},\frac{\partial}{\partial l_j}\right)=(W^{-1})_{ij}$$
	
\end{lemma}

\textbf{Proof.} It is known that $(\M_{g,n}(L_1,...,L_n),\omega_{WP})$ is a symplectic manifold, and the length functions $l_i$ are real analytic with respect to any choice of Fenchel-Nielsen coordinates. Set $\M_{g,n}(L_1,...,L_n)=M$. For a point $X\in\M_{g,n}(L_1,...,L_n)$, the symplectic structure gives both a linear isomorphism $\Omega:T_X^*M\to T_XM$ between the tangent and cotangent spaces, and a nondegenerate bilinear pairing $\omega_{WP}:T_XM\otimes T_XM\to\reals$. Furthermore, $\Omega$ is defined via Hamiltonian vector fields: for a cotangent vector $df\in T^*_XM$, $\Omega(df)$ is the unique tangent vector satisfying $df(\cdot)=\omega_{WP}(\Omega(df),\cdot)$. 

As a consequence of Lemma \ref{duality}, $\Omega(dl_i)=\frac{\partial}{\partial\tau_i}$ for all $dl_i$, i.e. $\frac{\partial}{\partial \tau_i}$ is the Hamiltonian vector field associated to the function $l_i$. Furthermore, by Lemma \ref{wolpertformula},
$$W_{ij}=\frac{\partial l_i}{\partial \tau_j}=\omega_{WP}\left(\frac{\partial}{\partial\tau_i},\frac{\partial}{\partial\tau_j}\right)$$
Since $\omega_{WP}$ is nondegenerate, we may conclude that $W$ is invertible at $X$ if and only if the tangent vectors $\left\{\frac{\partial}{\partial\tau_1},...,\frac{\partial}{\partial\tau_{6g-6+2n}}\right\}$ form a basis for $T_XM$. Applying $\Omega$ to this prospective basis, $W$ is invertible if and only if $\{dl_1,...,dl_{6g-6+2n}\}$ form a basis for $T^*_XM$, and this is equivalent to the $l_i$ forming a real analytic coordinate system by the Inverse Function Theorem.

For the second part of the lemma, we first apply the chain rule to $\frac{\partial}{\partial \tau_k}$:
$$\frac{\partial}{\partial\tau_k}=\sum_{i=1}^{6g-6+2n}\frac{\partial l_i}{\partial \tau_k}\frac{\partial}{\partial l_i}=\sum_{i=1}^{6g-6+2n} W_{ik}\frac{\partial}{\partial l_i}$$
Since $W$ is invertible,
$$\frac{\partial}{\partial l_i}=\sum_{k=1}^{6g-6+2n}(W^{-1})_{ik}\frac{\partial}{\partial \tau_k}$$
and so, applying Lemma \ref{duality},
$$\omega_{WP}\left(\frac{\partial}{\partial l_i},\frac{\partial}{\partial l_j}\right)=\sum_{k=1}^{6g-6+2n}(W^{-1})_{ik}\omega_{WP}\left(\frac{\partial}{\partial \tau_k},\frac{\partial}{\partial l_j}\right)$$
$$=\sum_{k=1}^{6g-6+2n}(W^{-1})_{ik}dl_k\left(\frac{\partial}{\partial l_j}\right)=\sum_{k=1}^{6g-6+2n}(W^{-1})_{ik}\delta_{kj}=(W^{-1})_{ij}\text{ $\blacksquare$ }$$

\

As a corollary, we obtain the following theorem.

\begin{theorem}[Wolpert, \cite{SW83}]
\label{cosineformula}

Fix a hyperbolic surface $X\in\M_{g,n}(L_1,...,L_n)$ and choose $6g-6+2n$ distinct simple closed multi-geodesics $\gamma_1,...,\gamma_{6g-6+2n}$. Define a $(6g-6+2n)\times(6g-6+2n)$ skew-symmetric matrix $W$ by $W_{ij}=w(\gamma_1,\gamma_2)$. If $W$ is invertible, then the length functions $l_i=l(\gamma_i)$ form local coordinates for $\M_{g,n}(L_1,...,L_n)$ in a neighborhood of $X$, and we have the following formulas on this neighborhood:
\begin{equation}\label{einitwpform}\omega_{WP}=\sum_{i<j}(W^{-1})_{ij}dl_i\wedge dl_j\end{equation}
\begin{equation}\label{ewpform}dV_{WP}=\sqrt{\det(W^{-1})}dl_1\wedge\cdots\wedge dl_{6g-6+2n}\end{equation}
	
\end{theorem}

Wolpert does not include Formula \ref{ewpform}, but it may be derived from Formula \ref{einitwpform} via a standard computation. 

\subsection{The Kontsevich Form}

In order to work with the Kontsevich form and obtain a combinatorial analogue for Theorem \ref{cosineformula}, we will need to define the \emph{ideal Wolpert number} of two geodesics on a surface or ribbon graph, $iw(\gamma_1,\gamma_2)$. This number can be thought of as an asymptotic version of the Wolpert number in which all angles become either $0$ or $\pi$.

Fix $X\in\M_{g,n}(L_1,...,L_n)$ such that its spine $S_X$ is trivalent. Recall that every corridor has an intercostal representing the local shortest path between the two boundary components incident to that corridor, and that the intercostals cut the surface into hexagons. We will call these hexagons \emph{sectors}; there is a bijection between sectors and vertices of $S_X$. Note that this definition makes sense even if some corridors do not contain their intercostals. 

We will now give two definitions of ideal Wolpert numbers. The first definition is more symmetric and closer to the definition of Wolpert numbers, while the second definition is more concrete and will be useful for proofs. Note that while ideal Wolpert numbers are defined for unoriented curves, the orientation-independence is somewhat subtle.

\begin{figure}[h]

\centering
\includegraphics[scale=0.23]{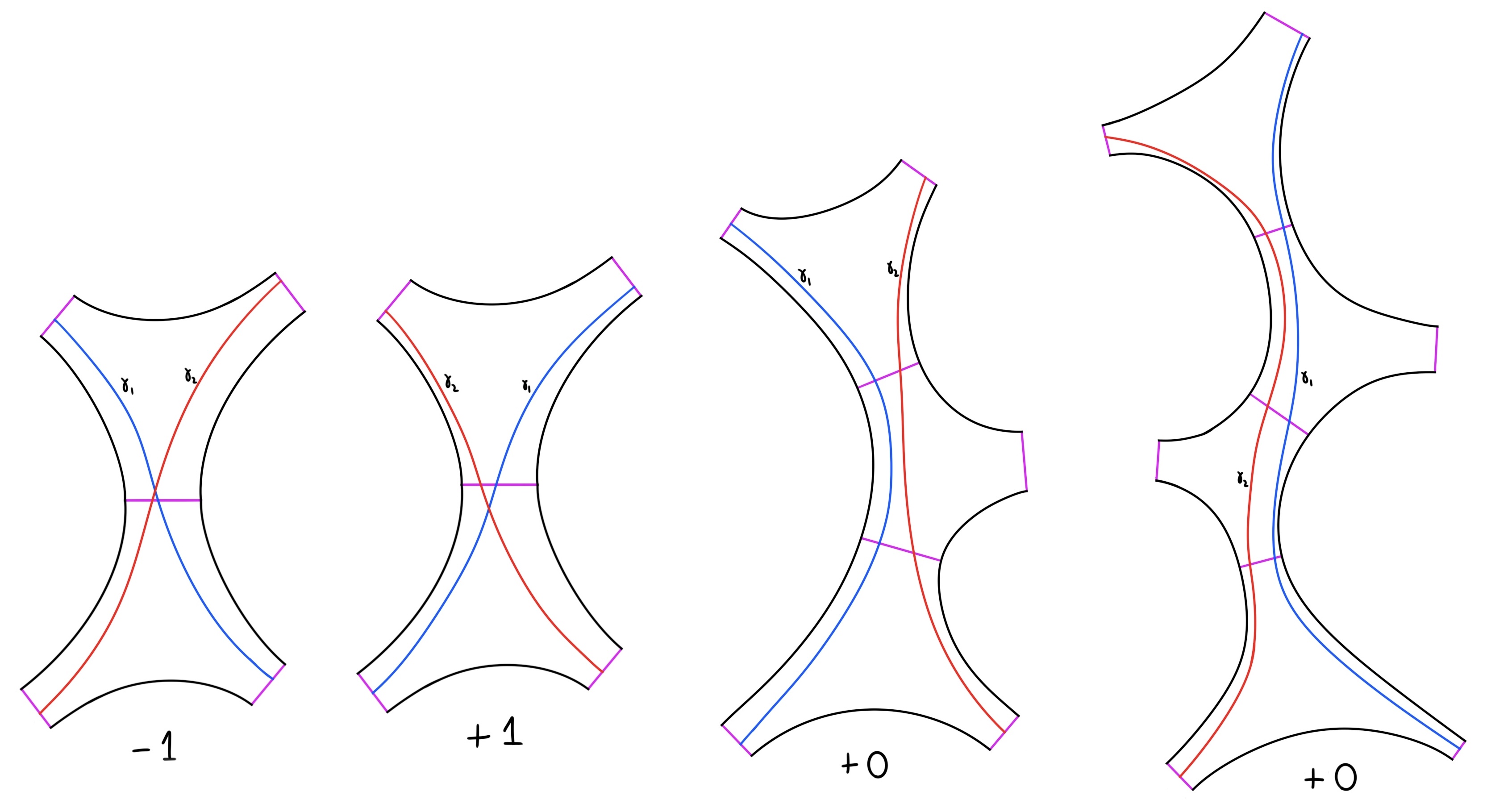}
\caption{Types of maximal intersection pairs, showing varying numbers of sectors between the end sectors of the pair. Far left: a type-$(-1)$ pair, where in both end sectors the endpoint of $\gamma_1'$ lies counterclockwise of the endpoint of $\gamma_2'$. Center-left: a type-$(1)$ pair, where in both end sectors the endpoint of $\gamma_1'$ lies clockwise of the endpoint of $\gamma_2'$. Right: Two type-$(0)$ pairs, where in one end sector the endpoint of $\gamma_1'$ lies counterclockwise from the endpoint of $\gamma_2'$, and in the other end sector the reverse holds.}

\label{maxpairdef}

\end{figure}

Take two geodesics $\gamma_1$ and $\gamma_2$ on a surface $X$. A \emph{maximal intersection pair} is a pair of subarcs $\gamma_1'\subset\gamma_1$ and $\gamma_2'\subset\gamma_2$ with the following property: $\gamma_1'$ and $\gamma_2'$ may be simultaneously oriented so that $\gamma_1'$ and $\gamma_2'$ start at two distinct intercostals in the same sector, pass through that sector and potentially multiple further sectors via the same intercostal(s), and then end at distinct intercostals (see Figure \ref{maxpairdef}). These are the surface equivalent of the pairs of arcs constructed in Definition \ref{graphintersectiondef}, although we do not yet require a version of property three in that definition. Each maximal intersection pair can then be labeled as type-$(0)$, type-$(1)$, or type-$(-1)$ as in Figure \ref{maxpairdef}, based on the configurations of the endpoints of $\gamma_1'$ and $\gamma_2'$ in the two end sectors. For example, if the endpoint of $\gamma_1'$ lies counterclockwise of the endpoint of $\gamma_2'$ in both end sectors, the pair is type-$(-1)$.

The assigned number is the \emph{signed crossing number} of the pair, and the ideal Wolpert number $iw(\gamma_1,\gamma_2)$ is the sum of signed crossing numbers for all maximal intersection pairs. Crucially, the signed crossing number is independent of the simultaneous orientation chosen. Furthermore, by following geodesics, the intersection points $p\in\gamma_1\cap\gamma_2$ of $\gamma_1$ and $\gamma_2$ are in bijection with the maximal intersection pairs that are type-$(1)$ and type-$(-1)$. Ideal Wolpert numbers satisfy the property $iw(\gamma_1,\gamma_2)=-iw(\gamma_2,\gamma_1)$. 

\begin{figure}[h]

\centering
\includegraphics[scale=0.15]{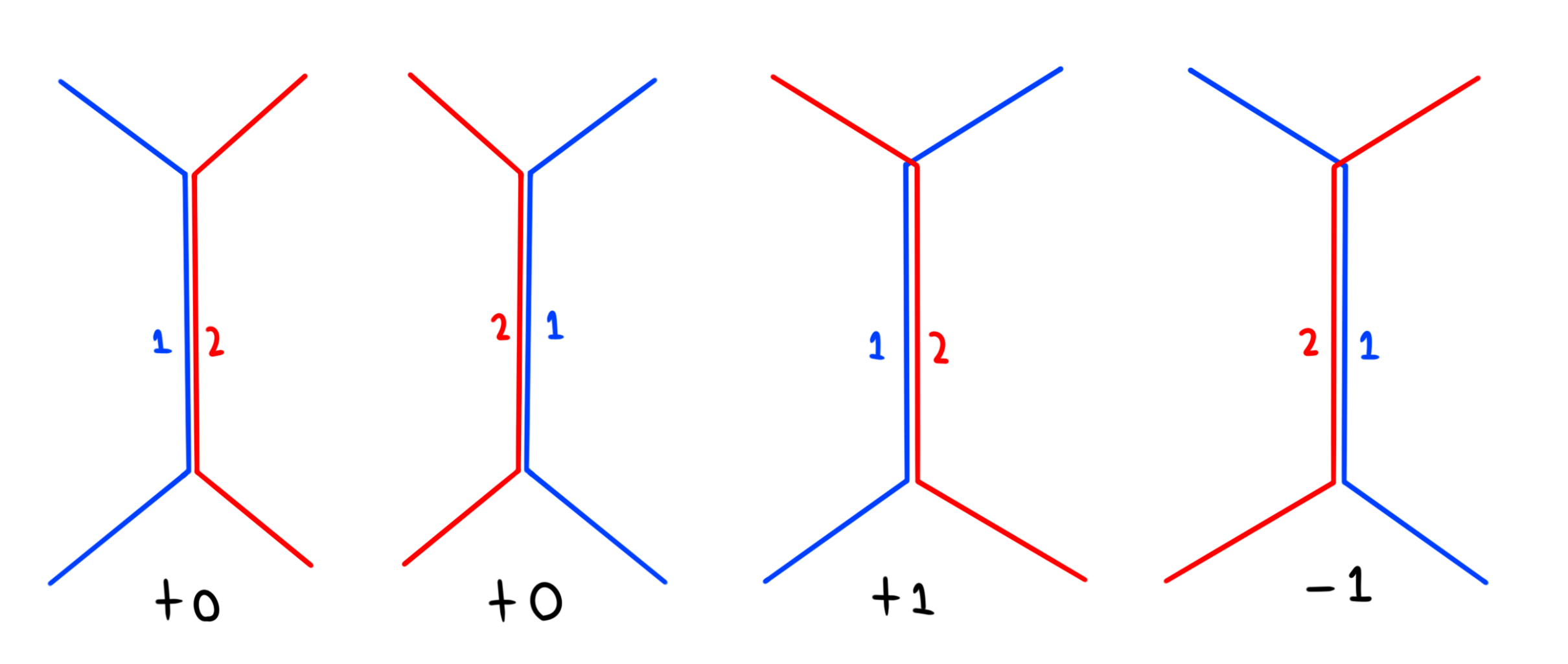}
\caption{A schematic drawing illustrating signed crossing numbers on graphs. Note that the arc where the two geodesics overlap may contain multiple edges and vertices, and so be incident to other edges `branching off' that are not shown here.}

\label{graphcrossing}

\end{figure}

If instead $\tilde{\gamma}_1$ and $\tilde{\gamma}_2$ lie on a ribbon graph $\Gamma$, a maximal intersection pair is a pair of subarcs satisfying properties one and two in Definition \ref{graphintersectiondef}, and we define the signed crossing number as in Figure \ref{graphcrossing}. Adding contributions from all pairs gives the ideal Wolpert number. In particular, every $\tilde{p}\in\tilde{\gamma}_1\cap\tilde{\gamma}_2$ has a well-defined crossing number $cr(\tilde{\gamma}_1,\tilde{\gamma}_2,\tilde{p})$ that equals $\pm 1$.

\begin{figure}[h]

\centering
\includegraphics[scale=0.2]{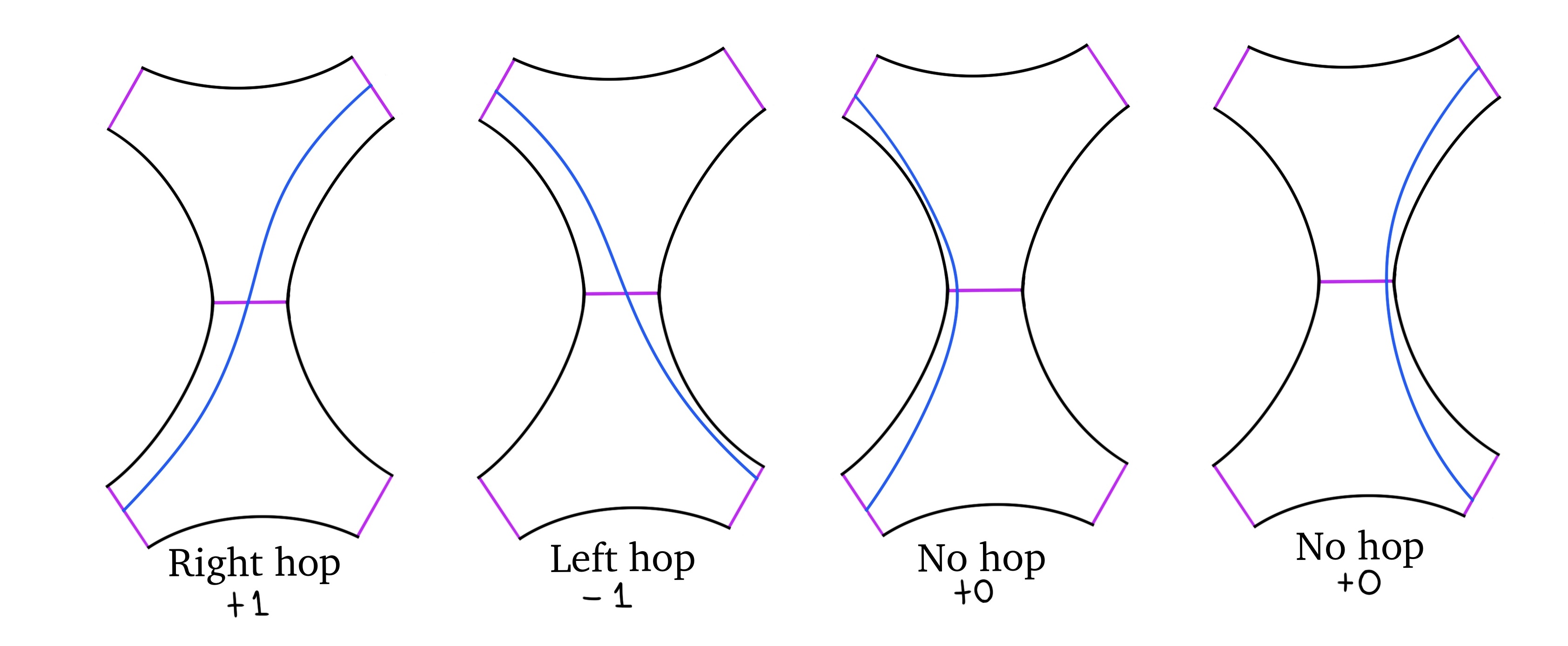}
\caption{The four possible ways a geodesic arc can pass through two adjacent sectors, with labels describing the different possible `hops' and the corresponding contributions to Formula \ref{egeodderiv}.}

\label{hopdef}

\end{figure}

We now give the second definition, and prove it is equivalent to the first. Given $X\in\M_{g,n}(L_1,...,L_n)$ and simple closed geodesics $\gamma_1$ and $\gamma_2$ on $X$, cut $X$ into sectors as above. Let $I$ be the collection of all points where $\gamma_1$ intersects an intercostal, and for all $i\in I$, let $\iota_i$ be the corresponding intercostal. We say that $\gamma_1$ makes a `left hop', `right hop', or no hop at $\iota_{i}$ as in Figure \ref{hopdef}, and assign a corresponding number $h_{i}$ from $\{-1,0,1\}$. We then let $c_i$ be the number of times $\gamma_2$ intersects $\iota_{i}$.

\begin{lemma}
\label{intersectiondef}

Let $X$, $\gamma_1$, $\gamma_2$, the $\iota_{i}$, and the constants $h_{i}$ and $c_{i}$. Then
\begin{equation}\label{egeodderiv}iw(\gamma_1,\gamma_2)=\sum_{j=1}^m h_{j}c_{j}\end{equation}

\end{lemma}

\textbf{Proof.} Choose a maximal intersection pair $\gamma_1'\subset\gamma_1,\gamma_2'\subset\gamma_2$, and assume it is type-$(1)$. Arbitrarily orient $\gamma_1$, and say the intercostals intersected by both $\gamma_1'$ and $\gamma_2'$ are $\iota_{i_1},...,\iota_{i_m}$, with corresponding hop numbers $h_{i_1},...,h_{i_m}$ for $\gamma_1$. We now claim that $h_{i_1}+...+h_{i_m}$ is equal to $1$. One may verify that, ignoring all `no hops', two right hops cannot occur sequentially, and two left hops cannot occur sequentially. Furthermore, based on the type-$(1)$ assumption, if all `no hops' are ignored, the sequence of right and left hops must start and end with a right hop. So this sequence looks like $RLR...LR$, and $h_{i_1}+...+h_{i_m}=1$ must hold. To finish, $\gamma_2'$ intersects every intercostal under consideration exactly once, so the pair contributes $+1$ to $iw(\gamma_1,\gamma_2)$ under the new definition. By symmetry, type-$(-1)$ pairs contribute $-1$, and an analogous argument works for type-$(0)$ pairs. $\blacksquare$

\

Note that both $|iw(\gamma_1,\gamma_2)|$ and $|w(\gamma_1,\gamma_2)|$ are bounded by the usual intersection number $i(\gamma_1,\gamma_2)$ - this is immediate for Wolpert intersection numbers, and follows for ideal Wolpert numbers due to the fact that maximal intersection pairs not corresponding to intersection points contribute zero to $iw(\gamma_1,\gamma_2)$.

The following lemma is a consequence of the fact that the topological structure of a spine graph reflects the sector gluing structure of the corresponding surface, and can be verified through straightforward casework.

\begin{lemma}
\label{preserveintersection}

Let $\gamma_1$ and $\gamma_2$ be geodesics on a surface $X\in \M_{g,n}(L_1,...,L_n)$, let $\tilde{\gamma}_1$ and $\tilde{\gamma}_2$ be the corresponding geodesics on $\Phi(X)$, and assume $\Phi(X)$ is trivalent. Then
$$iw(\gamma_1,\gamma_2)=iw(\tilde{\gamma}_1,\tilde{\gamma}_2)$$
and
$$i(\gamma_1,\gamma_2)=i(\tilde{\gamma}_1,\tilde{\gamma}_2).$$
	
\end{lemma}

We will now pursue a formula for the Kontsevich form. Recall that Andersen-Borot-Charbonnier-Giacchetto-Lewa{\'n}ski-Wheeler give in \cite{ABCGLW21} a method for constructing twist flows associated to geodesics in $MRG_{g,n}(L_1,...,L_n)$ (however, some care must be taken; see Remark \ref{twistproblem}).

\begin{lemma}[Andersen-Borot-Charbonnier-Giacchetto-Lewa{\'n}ski-Wheeler, \cite{ABCGLW21}, Prop. 3.11]
\label{graphduality}

Fix $\Gamma\in MRG_{g,n}(L_1,...,L_n)$, and assume $\Gamma$ is trivalent. Let $\gamma$ be a simple closed geodesic on $\Gamma$, let $l$ be the length of $\gamma$, and let $\frac{\partial}{\partial \tau}$ be the corresponding twist derivative. Then there is an open neighborhood $U$ of $\Gamma$ on which
$$dl(\cdot)=\omega_K\left(\frac{\partial}{\partial\tau},\cdot\right)$$
	
\end{lemma}

The following lemma extends the methods of Andersen-Borot-Charbonnier-Giacchetto-Lewa{\'n}ski-Wheeler, and is the ribbon-graph analogue of Wolpert's cosine formula. This lemma is the reason for defining ideal Wolpert numbers.

\begin{lemma}
\label{graphtwist}

Fix $\Gamma\in MRG_{g,n}(L_1,...,L_n)$, and assume $\Gamma$ is trivalent. Fix two simple closed geodesics $\tilde{\gamma}_1$ and $\tilde{\gamma}_2$ on $\Gamma$. Let $l_1=l(\tilde{\gamma}_1)$, and let $\frac{\partial}{\partial\tau_2}$ be the twist derivative corresponding to $\tilde{\gamma}_2$. Then
\begin{equation}\label{egraphtwist}\frac{\partial l_1}{\partial\tau_2}=iw(\tilde{\gamma}_1,\tilde{\gamma}_2)\end{equation}

\end{lemma}

\begin{figure}[h]

\centering
\includegraphics[scale=0.27]{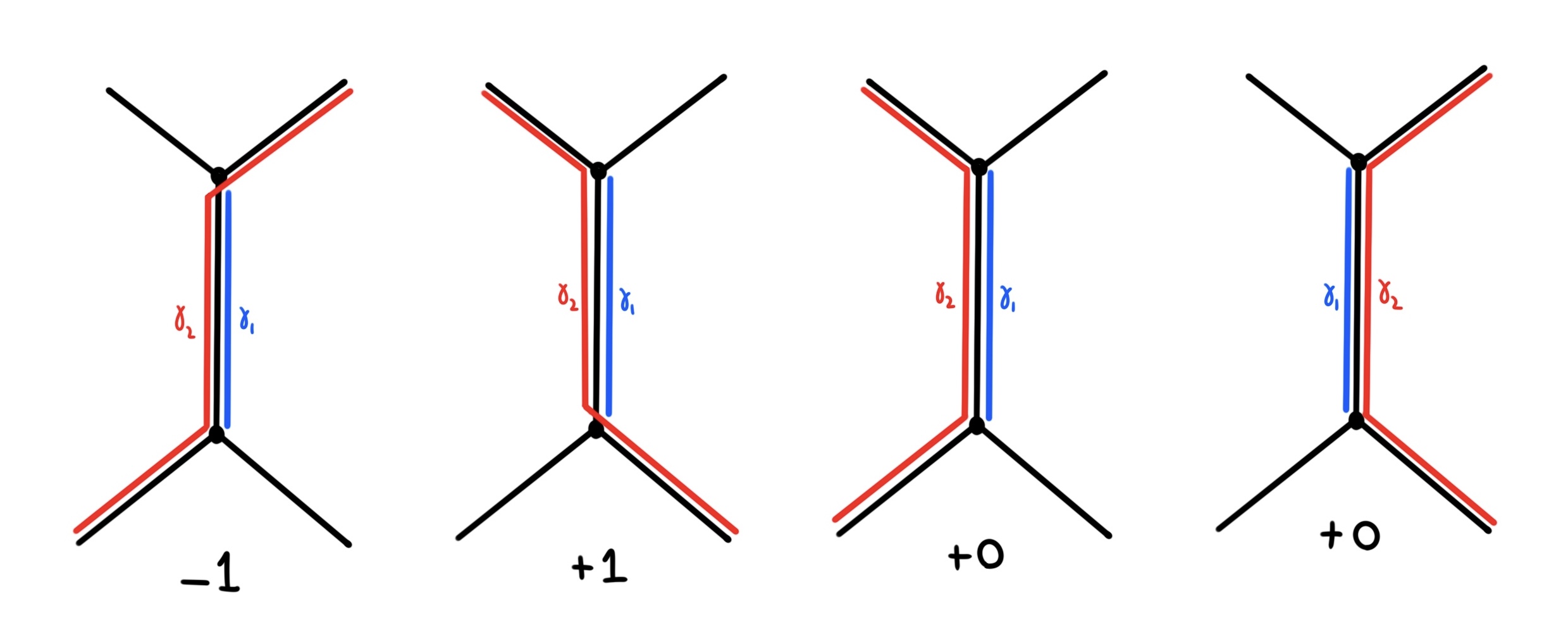}
\caption{Possible configurations for an edge in a metric ribbon graph, and two paths $\gamma_1$ and $\gamma_2$ passing through it. The adjacent edges traversed by $\gamma_2$ are marked, and determine the derivative of the edge length with respect to twisting along $\gamma_2$. These derivatives are indicated below each configuration.}

\label{edgelengthderiv}

\end{figure}

\textbf{Proof.} Arbitrarily orient $\tilde{\gamma}_1$ and $\tilde{\gamma}_2$, and say they traverse respective edges $e_{1,1},...,e_{m_1,1}$ and $e_{1,2},...,e_{m_2,2}$ in order. In the third paragraph of Section 3.2.1 in \cite{ABCGLW21}, Andersen-Borot-Charbonnier-Giacchetto-Lewa{\'n}ski-Wheeler argue that for a given edge $e$, $\frac{\partial l(e)}{\partial\tau_2}$ may be computed by identifying all $e_{j,2}$ equal to $e$ and then summing contributions from each $e_{j,2}$ according to the rule presented in Figure \ref{edgelengthderiv} (paying attention only to $\gamma_2$). To calculate $\frac{\partial l_1}{\partial\tau_2}$, we therefore sum over the $e_{j,2}$, and for each $e_{j,2}$ multiply its contribution by the number of times $\tilde{\gamma}_1$ traverses that edge. Say $e_{j,2}$ contributes $h_j'$, and $\tilde{\gamma}_1$ traverses that edge $c_j'$ times, so
\begin{equation}\label{eedgederiv}\frac{\partial l_1}{\partial\tau_2}=\sum_{j=1}^{m_2}h_j'c_j'\end{equation}

Say $X=\Phi^{-1}(\Gamma)$ and the simple closed geodesics $\gamma_1$ and $\gamma_2$ respectively correspond to $\tilde{\gamma}_1$ and $\tilde{\gamma}_2$. Edges of $\Gamma$ correspond to intercostals on $X$, and any traversal of an edge by $\tilde{\gamma}_i$ corresponds to an intersection of $\gamma_i$ with the corresponding intercostal. Furthermore, since vertices of $\Gamma$ correspond to sectors of $X$ and $\Gamma$ determines how these sectors are glued, we may observe that the first possibility in Figure \ref{edgelengthderiv} corresponds to a right hop, the second corresponds to a left hop, and the final two correspond to no hop (see the discussion before Lemma \ref{intersectiondef} for definitions of hops). However, here a right hop contributes $-1$ and a left hop contributes $+1$. 

Therefore, Equation \ref{eedgederiv} recovers \ref{egeodderiv}, only with $\gamma_1$ and $\gamma_2$ swapped and all signs flipped. Using Lemmas \ref{intersectiondef} and \ref{preserveintersection}, and the antisymmetry of ideal Wolpert numbers,
$$\frac{\partial l_1}{\partial\tau_2}=-iw(\gamma_2,\gamma_1)=-iw(\tilde{\gamma}_2,\tilde{\gamma}_1)=iw(\tilde{\gamma}_1,\tilde{\gamma}_2)\text{ $\blacksquare$}$$

\

By the work of Kontsevich (see \cite{MK92}), $(TRG_{g,n}(L_1,...,L_n),\omega_K)$ is known to be a symplectic manifold. Applying Lemmas \ref{graphduality} and \ref{graphtwist}, we obtain the following analogue of Theorem \ref{cosineformula}.

\begin{theorem}
\label{kontsevichformula}

Fix a trivalent ribbon graph $\Gamma\in MRG_{g,n}(L_1,...,L_n)$, and choose simple closed geodesics $\tilde{\gamma}_1,...,\tilde{\gamma}_{6g-6+2n}$. For each $\tilde{\gamma}_i$, define $k_i=l(\tilde{\gamma}_i)$, and let $S$ be the $(6g-6+2n)\times (6g-6+2n)$ skew-symmetric matrix defined by $S_{ij}=iw(\tilde{\gamma}_1,\tilde{\gamma}_2)$. Then $S$ is invertible if and only if the functions $k_i$ form local coordinates for $MRG_{g,n}(L_1,...,L_n)$ in a neighborhood of $\Gamma$. Furthermore, if the $k_i$ form local coordinates, we have the following formulas on this neighborhood:
$$\omega_K=\sum_{i<j}(S^{-1})_{ij}dk_i\wedge dk_j$$
\begin{equation}\label{ekontform}dV_K=\sqrt{\det(S^{-1})}dk_1\wedge\cdots\wedge dk_{6g-6+2n}\end{equation}
	
\end{theorem}

\textbf{Proof.} Using Lemmas \ref{graphduality} and \ref{graphtwist}, this result may be obtained by following essentially the same argument as for Lemma \ref{symplectic}. $\blacksquare$

\

Note that $S$ as defined above is constant on the interiors of maximal-dimension cells of $MRG_{g,n}(L_1,...,L_n)$, so $dV_K$ is a constant multiple of the Lebesgue volume form with respect to $k_1,...,k_{6g-6+2n}$ on such interiors. Morever, the $k_i$ are analytic (in fact, linear) with respect to the edge length functions on any such interior, which in turn are analytic with respect to any choice of Fenchel-Nielsen coordinates by Corollary \ref{analytic}. Lastly, $dV_{WP}$ is the Lebesgue volume form on $\M_{g,n}(L_1,...,L_n)$ with respect to those Fenchel-Nielsen coordinates, so we obtain the following result.

\begin{corollary}
\label{absolutecontinuity}

$\Phi^*dV_K$ is absolutely continuous with respect to $dV_{WP}$ on $\M_{g,n}(L_1,...,L_n)$.
	
\end{corollary}

The usefulness of Theorem \ref{kontsevichformula} hinges on being able to find geodesics $\tilde{\gamma}_1,...,\tilde{\gamma}_{6g-6+2n}$ such that the corresponding matrix $S$ is invertible. Following Do, for every edge $e$ of a trivalent $\Gamma\in MRG_{g,n}(L_1,...,L_n)$, we will define a corresponding simple multi-geodesic $\gamma_e$ using the boundary component(s) adjacent to $e$. The definition splits into cases depending on whether $e$ is incident to one or two boundary components of $\Gamma$, and is illustrated in Figure \ref{buildpants}. We will say such a $\gamma_e$ is \emph{derived from} $e$. Furthermore, Lemma 14 in \cite{ND10} essentially states that the $(6g-6+3n)\times(6g-6+3n)$ matrix $S'$ obtained by arbitrarily numbering the edges $e_1,...,e_{6g-6+3n}$, letting $\tilde{\gamma}_1,...,\tilde{\gamma}_{6g-6+3n}$ be the corresponding derived geodesics, and setting $S'_{ij}=iw(\gamma_i,\gamma_j)$ is the same (up to a constant multiple) as the matrix Do calls $B$ in \cite{ND10}. Applying Lemma 15 in \cite{ND10}, which states that $S'$ has rank $6g-6+2n$, and Lemma \ref{kontsevichformula}, we obtain:

\begin{lemma}
\label{geodchoice}

Assume $\Gamma\in MRG_{g,n}(L_1,...,L_n)$ is trivalent. Then there exist geodesics $\tilde{\gamma}_1,...,\tilde{\gamma}_{6g-6+2n}$ such that if $S$ is the matrix defined by $S_{ij}=iw(\gamma_i,\gamma_j)$ and $k_i=l(\tilde{\gamma}_i)$, then $S$ is invertible and hence the $k_i$ form local coordinate functions. Furthermore, all $i(\gamma_i,\gamma_j)$ are uniformly bounded in terms of $g$ and $n$.
	
\end{lemma}

\subsection{Angle Control}

We first prove angle control for great surfaces.

\begin{lemma}
\label{greatanglecontrol}

Fix $X\in\M_{g,n}(L_1,...,L_n)$, and assume that $X$ is $\eta$-great. Fix two closed oriented geodesics $\gamma_1$ and $\gamma_2$ in $X$, and assume these geodesics intersect at a point $p$ with angle of incidence $\theta$ and signed crossing number $C$, measured from $\gamma_1$ to $\gamma_2$. Then
$$\cos(\theta)=C+O(e^{-2\eta})$$
	
\end{lemma}

\textbf{Proof.} Let $H$ be a sector containing $p$. Extending $\gamma_1$ and $\gamma_2$ from $p$, we obtain geodesic arcs passing through the sector. These arcs each intersect two distinct intercostals, and so intersect a common intercostal corresponding to a corridor $W$.

\begin{figure}[h]

\centering
\includegraphics[scale=0.2]{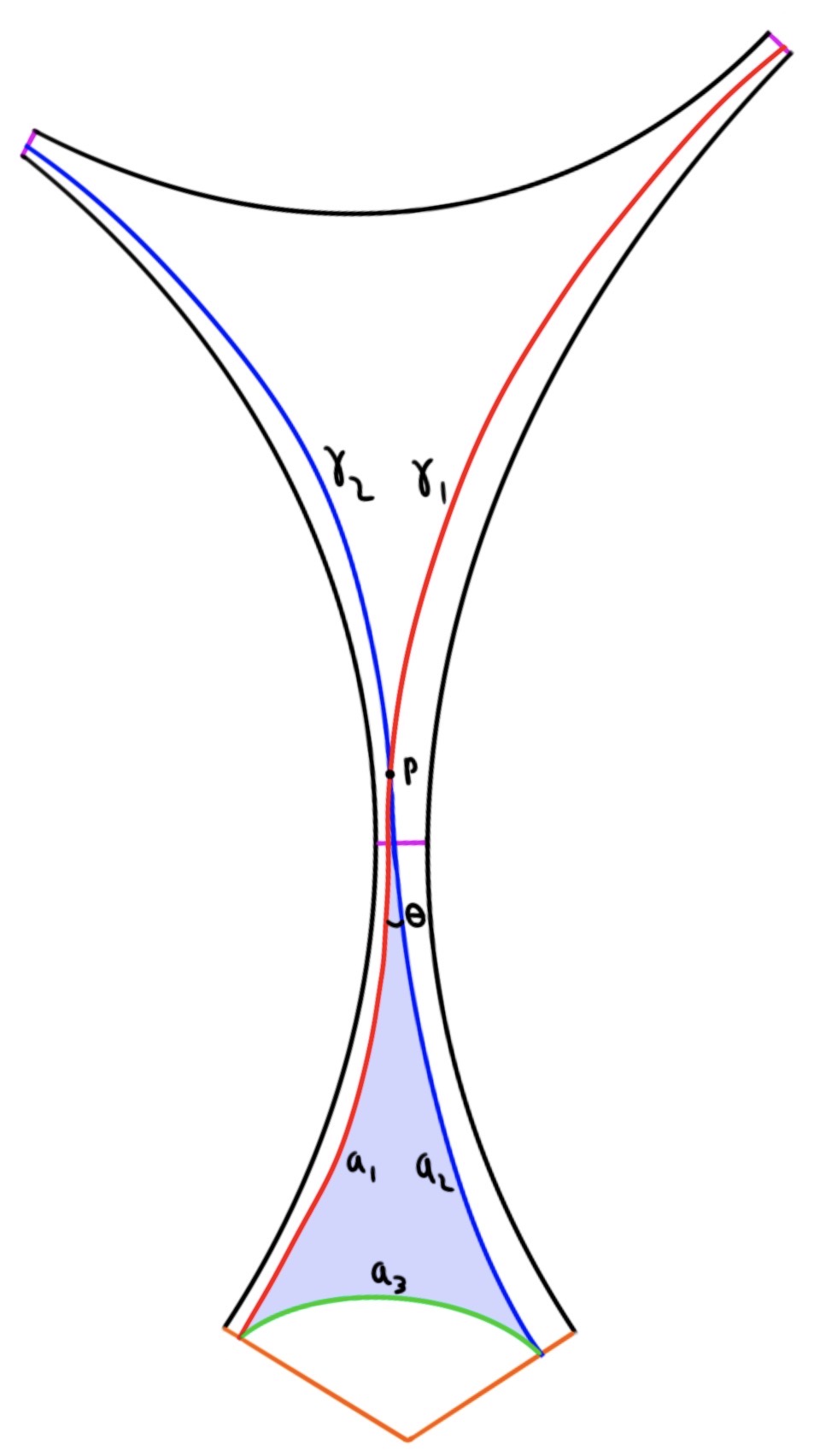}
\caption{The constructions performed in the proof of Lemma \ref{greatanglecontrol}. The triangle formed by $a_1$, $a_2$, and $a_3$ is shaded.}

\label{anglecontrol1}

\end{figure}

Now, swap labels as necessary so that $C=1$. Follow the geodesics from $p$ through the common intercostal until they exit $W$, labeling the resulting arcs $a_1$ and $a_2$. Let $a_3$ be the geodesic arc in $W$ between the exit points, so that the $a_i$ from a triangle. This entire construction is shown in Lemma \ref{anglecontrol1}. By Lemma \ref{triangle} we now have

$$\cos(\theta)=\coth(l(a_1))\coth(l(a_2))-\frac{\cosh(l(a_3))}{\sinh(l(a_1))\sinh(l(a_2))}$$
By Corollary \ref{ribcontrol} and the fact that all half-corridors in $X$ have boundary arcs of length at least $\Omega(\eta)$ by Lemma \ref{intercostalbalance}, $l(a_1)=\Omega(\eta)$, $l(a_2)=\Omega(\eta)$, and $l(a_3)=O(1)$. Standard asymptotics for $\coth$, $\cosh$, and $\sinh$ now yield the desired result. $\blacksquare$

\

\begin{lemma}
\label{intercostalintersectioncontrol}

Fix $X\in\M_{g,n}(L_1,...,L_n)$, and assume that $X$ is $\eta$-great. Let $\iota$ be an intercostal of a corridor $W$ in $X$, and let $\gamma$ be a geodesic intersecting $\iota$ at a point $p$. Let $\theta$ be the angle of incidence. Then
$$\sin(\theta)=1+O(e^{-\eta/2})$$
	
\end{lemma}

\begin{figure}[h]

\centering
\includegraphics[scale=0.14]{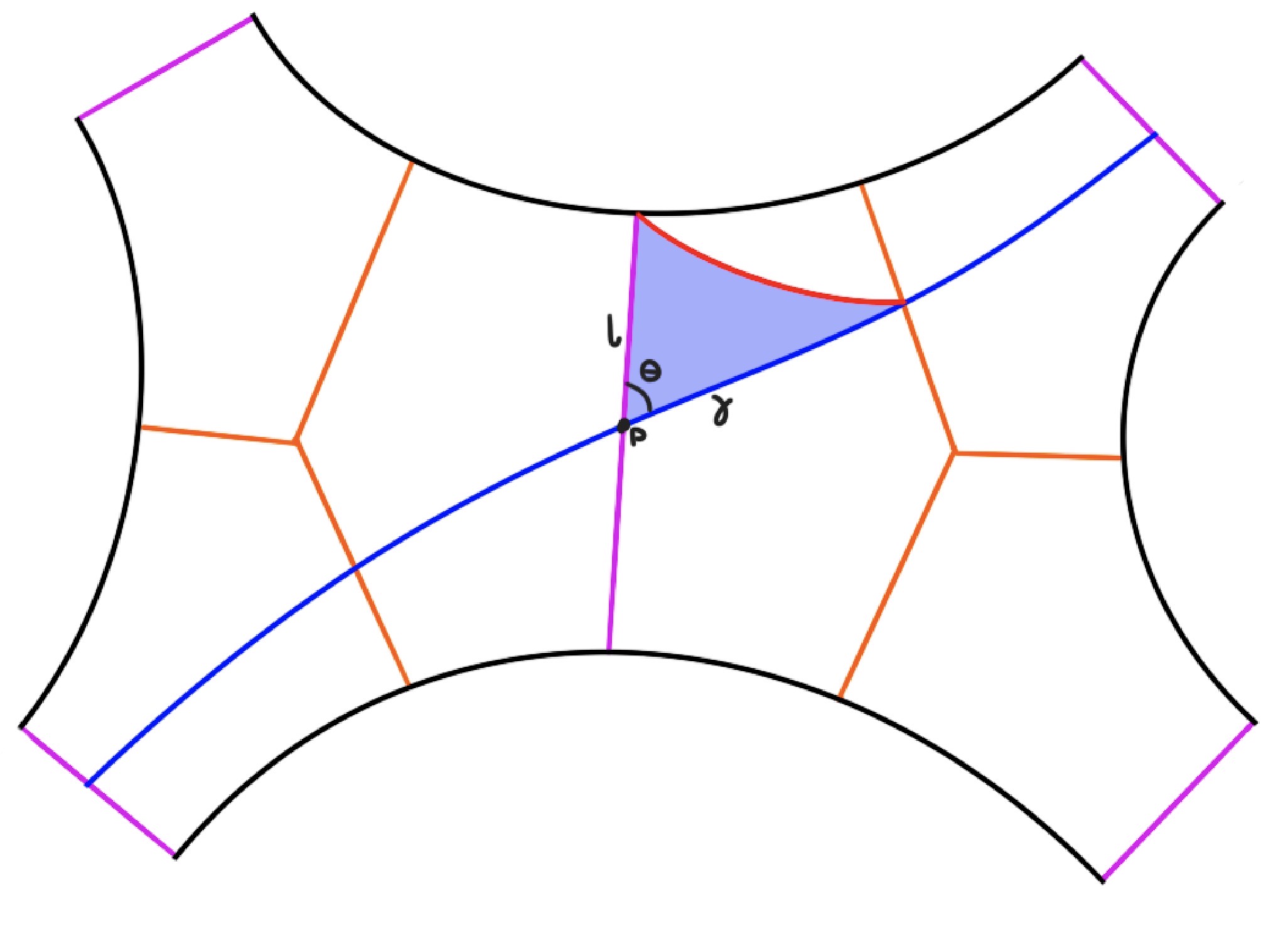}
\caption{The construction performed in Lemma \ref{intercostalintersectioncontrol}. The triangle is shaded.}

\label{anglecontrol2}

\end{figure}

\textbf{Proof.} Choose one pair of incident ribs of $W$, and a geodesic arc from the intersection of $\gamma$ with that pair to one endpoint of $\iota$, so that a triangle is formed with $\theta$ as one of the other angles - see Figure \ref{anglecontrol2}. For this triangle, the two sides not lying on $\iota$ are $\Omega(\eta)$ by Lemma \ref{intercostalbalance} and Corollary \ref{ribcontrol}, while the side lying on $\iota$ is $O(e^{-\eta/2})$. Therefore the two long sides differ in length by $O(e^{-\eta/2})$, and Lemma \ref{righttriangle} yields the desired result. $\blacksquare$

\

Angle control for good surfaces will require more work. Let $X\in\M_{g,n}(L_1,...,L_n)$ be a $(\epsilon,\eta)$-good surface with $B\leq\epsilon$ and $A\leq \eta$. So $\Phi(X)$ may have edges of length at most $\eta$, but no geodesics of length at most $\epsilon$. In $\Phi(X)$, cut at the midpoints of all edges with length at least $ \eta$. If any of the resulting components is not contractible, it must contain a simple closed geodesic of length at most $O(\eta)$. On $X$, this implies that if we cut along the intercostals corresponding to all corridors derived from edges with length at least $ \eta$, either all resulting components are simply connected or there is a simple closed geodesic on $X$ of length at most $O(\eta)+O(1)$, where the $O(1)$ comes from Corollary \ref{ribbound}. Therefore, as long as $\epsilon$ is chosen to be a sufficiently large multiple of $\eta$ in terms of $g$ and $n$, cutting $X$ in this way will decompose it into polygons with alternating intercostals and boundary arcs as sides. We will call these polygons \emph{supersectors}. 

\

\textbf{From now on, we will assume that $B$ is a sufficiently large multiple of $A$ as above, and that whenever a surface is assumed to be $(\epsilon,\eta)$-good or $(\epsilon,\eta)$-great, then $\epsilon$ is a sufficiently large multiple of $\eta$.}

\

Supersector decompositions may also be applied to ribbon graphs, with or without a metric structure. In this case, take a ribbon graph $\Gamma$, and choose some subset of edges $E_c$ such that cutting each $e\in E_c$ into two half-edges at its midpoint decomposes $\Gamma$ into some number of simply connected components (i.e., trees), such that all leaves of such trees are half-edges and all non-leaf edges are full edges. We will also call these components supersectors. The \emph{supersector boundary points} of a given supersector are the endpoints of its leaves, and these endpoints are edge midpoints in $\Gamma$. For $X\in\M_{g,n}(L_1,...,L_n)$ admitting a supersector decomposition, set $\Gamma=\Phi(X)$. Then there is a natural bijection between intercostals on $X$ and edge midpoints on $\Gamma$, and decomposing $X$ into supersectors gives a corresponding supersector decomposition on $\Gamma$.

\begin{definition}
\label{goodframedef}

Let $X$ be $(\epsilon,\eta)$-good for $\epsilon$ a sufficiently large multiple of $\eta$ as above, and let $\gamma_1,...,\gamma_{6g-6+2n}$ be a set of simple (multi-)geodesics on $X$. We say $\gamma_1,...,\gamma_{6g-6+2n}$ is a \emph{good geodesic frame} if any two arcs $\xi_i\subset\gamma_i$ and $\xi_j\subset\gamma_j$ crossing the same supersector $P$ either do not intersect or intersect a common intercostal on the boundary of $P$.
	
\end{definition}

\begin{figure}
\centering

\includegraphics[scale=0.14]{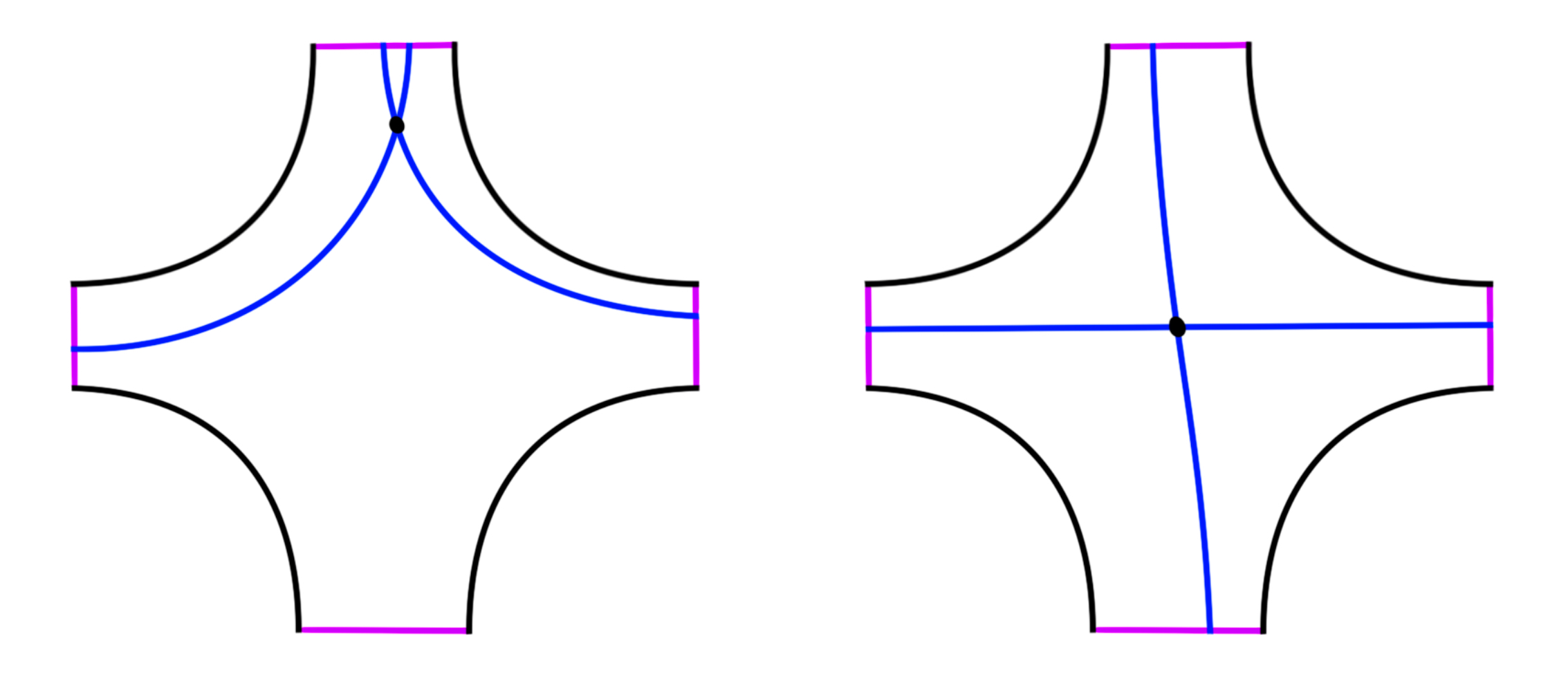}
\caption{Two illustrations showing a supersector and two intersecting geodesics. On the left, the geodesics do not form a bad intersection point. On the right, the geodesics form a bad intersection point.}

\label{badcross}
\end{figure}

Take two simple closed geodesics $\gamma_1$ and $\gamma_2$ intersecting at a point $p\in\gamma_1\cap\gamma_2$ on a surface $X$ admitting a supersector decomposition. Say $p$ lies in a supersector $P$. We will define $p$ to be a \emph{bad intersection point} if the respective arcs of $\gamma_1$ and $\gamma_2$ traversing $P$ and passing through $p$ do not intersect a common intercostal on the boundary of $P$ (see Figure \ref{badcross}). Likewise, given two geodesics $\tilde{\gamma}_1$ and $\tilde{\gamma}_2$ intersecting in a pair of arcs $\tilde{p}$, say $\tilde{p}$ is a \emph{bad intersection} if its body is contained in the interior of a single supersector.

\begin{lemma}
\label{badintersectionequiv}
Let $X\in\M_{g,n}(L_1,...,L_n)$ be a surface admitting a supersector decomposition, and fix simple closed geodesics $\gamma_1$ and $\gamma_2$ on $X$ intersecting at a point $p\in\gamma_1\cap\gamma_2$. Say $\Gamma=\Phi(X)$, and $\tilde{\gamma}_1$ corresponds to $\gamma_1$, $\tilde{\gamma}_2$ corresponds to $\gamma_2$, and $\tilde{p}\in\tilde{\gamma}_1\cap\tilde{\gamma}_2$ corresponds to $p$. Then $p$ is a bad intersection point if and only if $\tilde{p}$ is a bad intersection.
	
\end{lemma}

\textbf{Proof.} Assume $p$ is a bad intersection point occurring in a supersector $P$, and let $\tilde{P}$ be the corresponding supersector in $\Gamma$. Let $\gamma_i'$ be the arc of $\gamma_i$ traversing $P$ that contains $p$ for $i=1,2$. If the body of $\tilde{p}$ intersected one of the boundary points of $\tilde{P}$, then both $\gamma_1'$ and $\gamma_2'$ would intersect the corresponding intercostal, and $p$ would not be a bad intersection point. Likewise, if $\gamma_1'$ and $\gamma_2'$ were incident to a common intercostal, then $\gamma_1$ and $\gamma_2$ would not diverge until at least the sector on the other side of that intercostal, so the body of $\tilde{p}$ would pass into that sector. This contradicts that the body of $\tilde{p}$ is contained in the interior of $\tilde{P}$. $\blacksquare$

\

The reason for avoiding bad intersection points is as follows:

\begin{lemma}
\label{supersectorintersectioncontrol}

Let $X\in\M_{g,n}(L_1,...,L_n)$ be a $(B,A)$-good surface, and let $\gamma_1,...,\gamma_{6g-6+2n}$ be a good geodesic frame. Also assume that $B$ is a sufficiently large multiple of $A$. Fix two geodesics from the frame, say $\gamma_1$ and $\gamma_2$, and assume these geodesics intersect at a point $p$ with angle of incidence $\theta$ and crossing number $C$, measured from $\gamma_1$ to $\gamma_2$. Then
$$\cos(\theta)=C+O(e^{-A/2})$$
	
\end{lemma}

\textbf{Proof.} The proof is essentially identical to the proof of Lemma \ref{greatanglecontrol}, with the key property being that $\gamma_1$ and $\gamma_2$ can be followed from $p$ to a common intercostal attached to a long corridor. $\blacksquare$

\subsection{Dehn Twists and Winding}

We will use Dehn twists to obtain good geodesic frames. The key idea is that, given two simple closed geodesics $\gamma$ and $\eta$ on a surface, $\eta$ may be made to `wind around' $\gamma$ as many times as desired using repeated Dehn twists.

Recall the definition of the earthquake flow: given a surface $X\in\M_{g,n}(L_1,...,L_n)$ and a simple closed geodesic $\zeta$ on $X$, $\zeta$ may be identified with $\reals/l(\zeta)\intgr$. For any $t\in\reals$, we may cut $X$ along $\zeta$, obtaining two new boundary components we arbitrarily designate `left' and `right'. We then glue point $x$ on the right copy of $\zeta$ to point $x+t$ on the left copy, obtaining a new surface $X_t$. Given another simple closed geodesic $\eta$ on $X$, we may construct a new geodesic $\eta_t$ by, whenever $\eta$ approaches $\zeta$ from the `left', making a left turn at $\zeta$, traveling distance $t$ along $\zeta$, and turning right to continue along $\eta$. While this procedure produces a curve that is not a geodesic and may overlap itself, we can declare $\eta_t$ to be the unique geodesic in that curve's homotopy class. This $\eta_t$ will be simple. 

Furthermore, for any geodesic $\gamma$ either equal to $\zeta$ or disjoint from $\zeta$, $\gamma$ is fixed by the earthquake, and there is a natural bijection between points in $\gamma\cap \eta$ on $X$ and points in $\gamma\cap \eta_t$ on $X_t$. This bijection comes from the fact that points in $\gamma\cap\eta_t$ `slide' continuously along $\gamma$ as $t$ increases or decreases, and cannot be created or destroyed due to the assumption that $\gamma$ is equal to or disjoint from $\zeta$.

While the map $t\to X_t$ itself is $l(\gamma)$-periodic, $t\to\eta_t$ is crucially not periodic, and we may define the \emph{positive Dehn twist} of $\eta$ around $\gamma$ to be $\eta_{l(\gamma)}$ on $X_{l(\gamma)}=X$. 

The key technical lemma is as follows. The $\zeta=\gamma$ case was first explicitly stated by McMullen in \cite{CM98}, but was essentially proved by Kerckhoff in \cite{SK83}. The case when $\zeta$ is disjoint from $\gamma$ is not stated by either author, but follows near-immediately from Kerckhoff's proof.

\begin{lemma}
\label{twistangle}

Let $X$ be a hyperbolic surface, and let $\gamma$ and $\xi$ be two simple closed geodesics on $X$ intersecting at a point $p$. Let $\zeta$ be a third simple closed geodesic satisfying either $\zeta=\gamma$ or $\zeta\cap\gamma=\emptyset$. Let $\frac{\partial}{\partial \tau}$ be the infinitesimal twist derivative along $\zeta$ evaluated at $X$. Let $\eta(s)$ be a unit-speed parameterization of $\eta$ satisfying $\eta(0)=p$, and for a given $s\in\reals$ such that $\eta(s)\in\eta\cap\zeta$, let $\theta_{\eta(s)}$ be the angle between $\zeta$ and $\eta$ at $q$. Then
\begin{equation}
\label{etwistangle}
\frac{\partial \theta}{\partial \tau}=\sum_{\{s: \ \eta(s)\in\eta\cap\zeta\}}e^{-|s|}\sin(\theta_{\eta(s)})
\end{equation}
Furthermore, this quantity is positive.
\end{lemma}

While we will not include a full explanation of this formula, the basic idea is to first consider the scenario when $X=\H$ and $\gamma$ and $\zeta$ intersect at only a single point $q$. This case may be handled by hyperbolic trigonometry - in particular, the impact on $\theta$ of twisting along $\zeta$ decays exponentially in $d(p,q)$. The full case of Lemma \ref{twistangle} is then obtained by lifting to the universal cover and using linearity of derivatives.

Let $\gamma$ and $\eta$ be simple closed geodesics on a surface $X$, intersecting at a point $p$ (and possibly with other points of intersection). It is helpful to define a real-valued auxiliary function $f_{\gamma,\xi,p}(t)=f(t)$, which will, in informal terms, estimate the arclength on either side of $p$ for which $\xi$ `fellow travels' $\gamma$ after a twist along $\gamma$ of length $t$. 

Specifically, let $X_t$ be the surface obtained by performing a positive twist along $\gamma$ of distance $t$. As above, for any $t$ we have simple closed geodesics $\eta_t$ and $\zeta_t$ respectively homotopic to the images of $\eta$ and $\zeta$, as well as the image of $\gamma$, which we also denote as $\gamma$ in a slight abuse of notation. Furthermore, we have an intersection point $p_t\in\gamma\cap\eta_t$, uniquely defined by the properties that $p_0=p$ and that $p_t$ varies continuously with $t$.

Given any $t$, let $\gamma(s)$ and $\xi_t(s)$ be respective unit-speed parameterizations of $\gamma$ and $\xi_t$ with the property that $\gamma(0)=\eta_t(0)=p_t$. Then we wish to construct $f(t)$ such that the following properties are satisfied:
\begin{enumerate}
\item $f(t)$ is smooth for all $t$.
\item $f(t)$ is strictly increasing in $t$, and for all $t$ and $\delta$,
$$f(t+2\delta+4)\geq f(t)+\delta$$
\item For all $t$ satisfying $f(t)>0$, and all $s\in [-f(t),f(t)]$, we have
$$d_{X_t}(\gamma(s),\xi_t(s))\leq 4$$
\end{enumerate}
The last condition states that $\gamma$ and $\xi_t$ fellow travel for distance $f(t)$ in both directions from $p_t$.

\begin{lemma}
\label{shiftvalue}
There exists a function $f(t)$, dependent on $\gamma,\xi_t,p$, satisfying the three properties listed above.
\end{lemma}

\begin{figure}[h]

\centering
\includegraphics[scale=0.15]{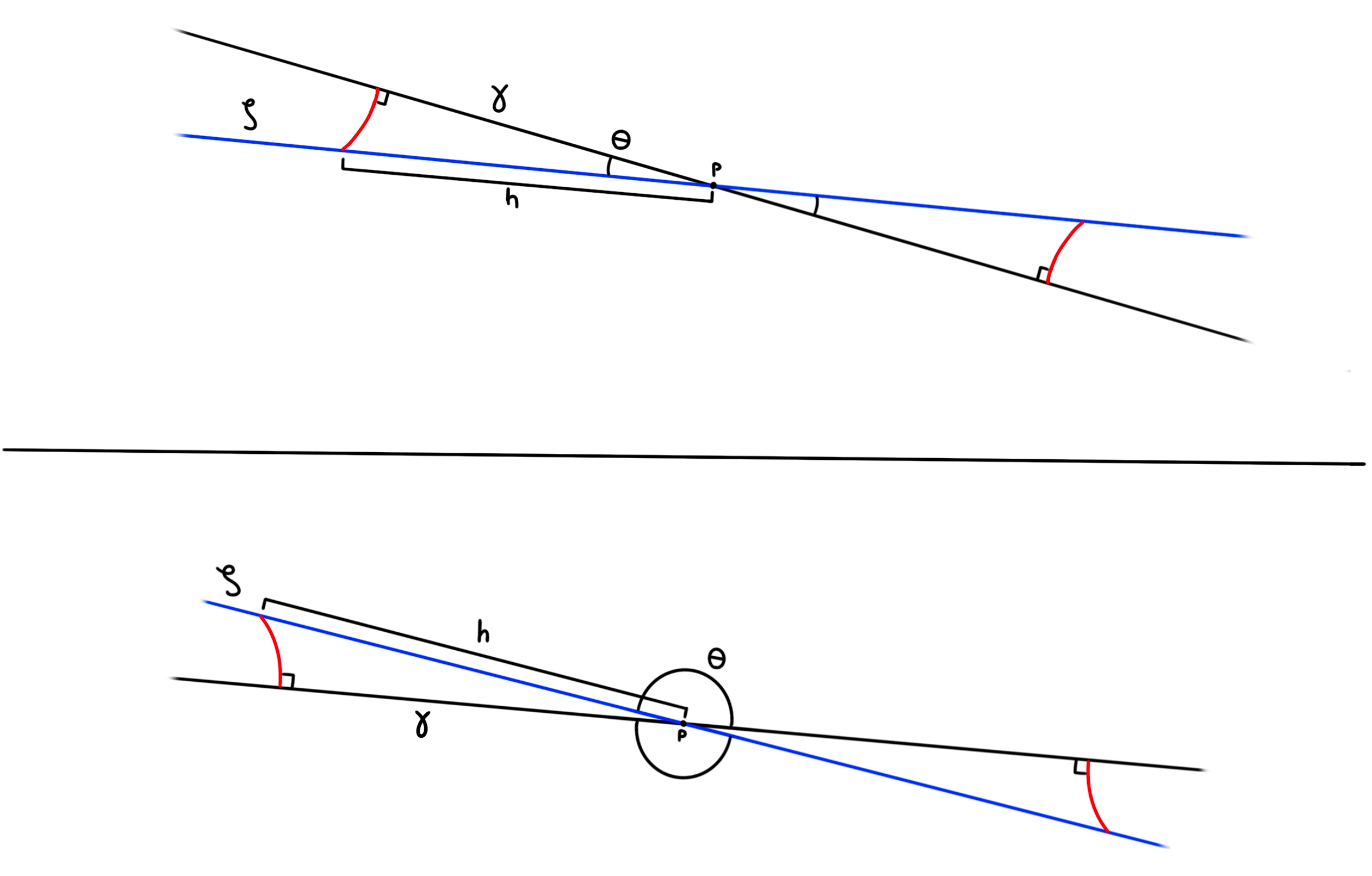}
\caption{The triangles used to define $f_{\gamma,\zeta,p}(t)$.}

\label{followingdef}
	
\end{figure}

\textbf{Proof.} The geodesics $\gamma$ and $\xi_t$ form four angles at $p_t$, two of which will be at most $\pi/4$. Choosing either of these small angles, we build a triangle such that the edge opposite that angle meets $\gamma$ at a right angle and has length $\sinh^{-1}(1)\approx 0.88$, as in Figure \ref{followingdef}. Let $h$ be the length of the subarc of $\xi$ forming the hypotenuse of this triangle.

If $0<\theta(t)\leq\pi/4$, we have the top case of Figure \ref{followingdef}, and we set $f(t)=-h$. If $3\pi/4\leq\theta(t)<\pi$, we have the bottom case of Figure \ref{followingdef} we set $f(t)=h$. Since $\theta(t)$ is a smooth increasing function, $\pi/4<\theta(t)<3\pi/4$ on at most a single interval of $\reals$, and we may define $f(t)=g(t)$ on this interval for some arbitrary interpolating function $g(t)$ such that the resulting $f(t)$ is increasing and smooth. By standard trigonometry,
$$f(t)=\begin{cases}\sinh^{-1}\left(\frac{1}{\sin(\theta(t))}\right)& 3\pi/4\leq \theta(t)<\pi\\ g(t)& \pi/4<\theta(t)<3\pi/4\\ -\sinh^{-1}\left(\frac{1}{\sin(\theta(t))}\right)& 0<\theta(t)\leq\pi/4\end{cases}$$
We will now verify the third condition. On the region where $f(t)=g(t)$ and $f(t)>0$, it is straightforward to verify from the definitions and the increasing condition that $f(t)\leq 2$, so for $s\in [-f(t),f(t)]$, clearly $d_{X_t}(\gamma(s),\xi_t(s))\leq 4$. Otherwise, we have from the definition of $f(t)=h$ that for $s\in [-f(t),f(t)]$, $d_{X_t}(\gamma(s),\xi_t(s))\leq \sinh^{-1}(1)<4$.

Lastly, we will verify the second condition, using the control on $\theta'(t)$ provided by Lemma \ref{twistangle}. If $\pi/4<\theta(t)<3\pi/4$, then by Equation \ref{etwistangle} we have $\theta'(t)\geq\sin(\theta(t))>\frac{1}{2}$, so the length of the interval where $f(t)=g(t)$ has length at most, say, $4$. For the second part, if $\zeta=\gamma$, then Equation \ref{etwistangle} yields $\theta'(t)\geq\sin(\theta(t))$. Furthermore, if $3\pi/4\leq\theta(t)<\pi$, then
$$f'(t)=-\frac{1}{\sqrt{1+\frac{1}{\sin^2(\theta(t))}}}*\frac{\cos(\theta(t))\theta'(t)}{\sin^2(\theta(t))}=\frac{-\cos(\theta(t))}{\sqrt{1+\sin^2(\theta(t))}}*\frac{\theta'(t)}{\sin(\theta(t))}$$
The first term is at least $\frac{1}{2}$ by the restriction on $\theta(t)$, while the second term is at least $1$, so $f'(t)\geq\frac{1}{2}$. An analogous computation holds for $0<\theta(t)\leq \pi/4$. $\blacksquare$

The next lemma essentially says that pairs of geodesic arcs originating from bad intersection points cannot fellow-travel for too long.

\begin{lemma}
\label{functionbadbound}

Let $X\in\M_{g,n}(L_1,...,L_n)$ be $(\epsilon,\eta)$-good. Let $\gamma$ and $\xi$ be two simple closed geodesics on $X$, and fix $p\in\gamma\cap\xi$. Choose either of the two angles formed by $\gamma$ and $\xi$ and $p$ with angle at most $\pi/4$, and define respective unit-speed parameterizations $\gamma(s)$ and $\xi(s)$ of $\gamma$ and $\xi$ traveling in the directions given by this angle for positive $s$.

Assume $p$ is a bad intersection point, and define $s_{\textup{close}}$ to be the largest value of $s$ such that $d_X(\gamma(s),\xi(s))\leq 4$. Then $s_{\textup{close}}\leq l(\gamma)$.
	
\end{lemma}

\textbf{Proof.} Take a supersector decomposition on $X$, and assume $p$ lies in a supersector $S$.  We may travel along $\gamma$ and $\xi$ in the directions given by the chosen angle until reaching points $p_\gamma$ and $p_\xi$, respectively lying on $\iota_\gamma$ and $\iota_\xi$. By the bad intersection assumption, $\iota_\gamma\neq\iota_\xi$.

\begin{figure}[h]

\centering
\includegraphics[scale=0.15]{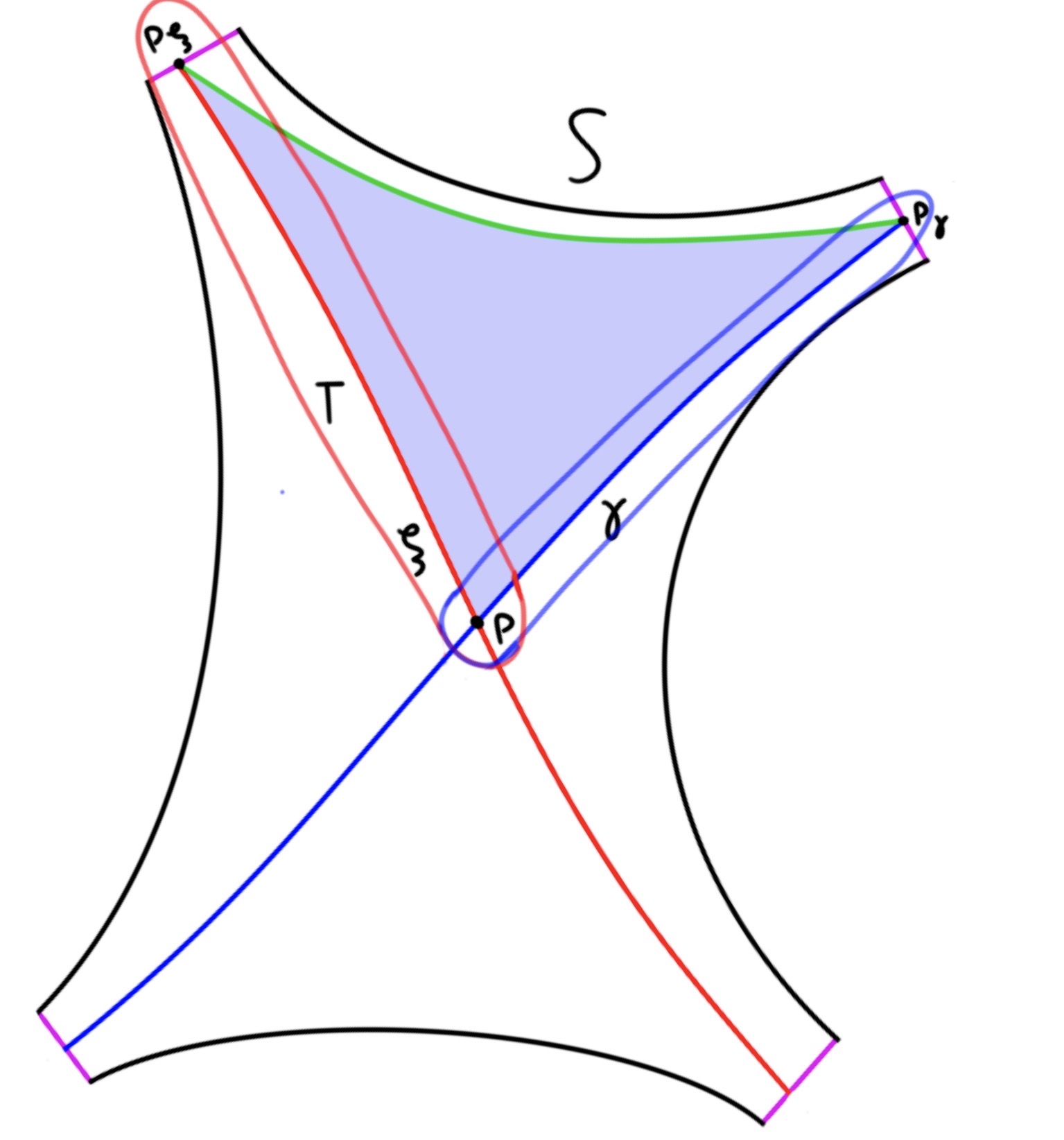}
\caption{An illustration of the construction used in Lemma \ref{functionbadbound}. The triangle $T$ is shaded, and the $4$-neighborhoods of the two relevant sides of $T$ are marked. The gist of the argument is that these neighborhoods do not possess an extremely long region where they overlap.}

\label{supersectortriangle}

\end{figure}

Consider the geodesic triangle $T$ in $S$ with vertices $p$, $p_\gamma$, and $p_\xi$ - see Figure \ref{supersectortriangle} for an illustration. Since all intercostals on the boundary of $S$ lie in corridors corresponding to spine edges with length at least $\eta$, one may straightforwardly verify that all side lengths of this triangle are $\Omega(\eta)$. As a consequence, elementary hyperbolic trigonometry yields that for $\eta$ sufficiently large, $s_{\text{close}}$ is bounded above by the minimum of the two side lengths of $T$ incident to $p$. These side lengths are $d_S(p,p_\gamma)$ and $d_S(p,p_{\xi})$. As $S$ is simply connected, these lengths are lengths of respective subarcs of $\gamma$ and $\xi$, and so are respective lower bounds for $l(\gamma)$ and $l(\xi)$. $\blacksquare$

\begin{lemma}
\label{functionlengthbound}

Let $X\in\M_{g,n}(L_1,...,L_n)$ be $(B,A)$-good. Let $\gamma$ and $\zeta$ be two distinct simple closed geodesics, fix $p\in\gamma\cap\zeta$, and let $\theta$ be the intersection angle of $\gamma$ and $\zeta$ at $p$, measured counterclockwise from $\gamma$ to $\zeta$. Assume that $\theta<\pi/4$, and let $h$ be the length of the hypotenuse of the triangle constructed as in Figure \ref{followingdef}. Recalling that $c(\cdot)$ measures the combinatorial length of a geodesic,
$$c(\zeta)\geq \left\lfloor\frac{h}{2l(\gamma)}-2\right\rfloor c(\gamma)$$

\end{lemma}

\textbf{Proof.} Let $U$ be the cylindrical neighborhood around $\gamma$ with diameter $\sinh^{-1}(1)$. If $k= \left\lfloor\frac{h}{2l(\gamma)}\right\rfloor$, then $\zeta$ lies within this cylinder for at least length $kl(\gamma)$. An argument as in the previous lemma shows that if $B$ is sufficiently large, then this arc of $\zeta$ intersects the same intercostals as some path following $\gamma$ with length at least $(k-1)l(\gamma)$. But any such path completes at least $k-2$ full loops around $\gamma$, and so intersects $(k-2)c(\gamma)$ intercostals when counted with multiplicity. $\blacksquare$

\begin{lemma}
\label{twistlengthincrease}

Let $X$ be a hyperbolic surface, and let $\gamma$ and $\zeta$ be two distinct simple closed geodesics. Let $\zeta_k$ be the result of applying $k$ positive Dehn twists to $\zeta$. Then
$$c(\zeta_k)\leq c(\zeta)+i(\gamma,\zeta)kc(\gamma)$$
	
\end{lemma}

\textbf{Proof.} One way to define a Dehn twist is to take $\zeta$, add $k$ revolutions around $\gamma$ at every point of $\gamma\cap\zeta$ with the correct directions, and homotope the result to a geodesic. The above formula holds with equality before the homotopy step, so it will become an upper bound after homotopy. $\blacksquare$

\subsection{Building a Good Geodesic Frame}

As an initial step towards building good geodesic frames, we will construct particular kinds of pants decompositions, following the work of Felikson and Natanzon in \cite{FN12}. First, we need an intermediate result. The constant $10$ appearing in this result is fairly arbitrary.

\begin{lemma}
\label{geodflip}

Fix $X\in\M_{g,n}(L_1,...,L_n)$ such that $\Gamma=\Phi(X)$ is trivalent. Let $\gamma_2$, $\gamma_3$, $\gamma_4$, and $\gamma_5$ be four simple closed geodesics bounding a four-holed sphere $S$ in $X$, and let $\gamma_1$ be a simple closed geodesic contained in $S$ such that $\gamma_1$ cuts $S$ into two pants $P_1$ and $P_2$. Then there exists a simple closed geodesic $\zeta_1$ that is also contained in $S$, also cuts $S$ into two pants, and that intersects $\gamma_1$ in exactly two points. Furthermore, $\zeta_1$ satisfies
$$c(\zeta_1)\leq 10(c(\gamma_1)+c(\gamma_2)+c(\gamma_3)+c(\gamma_4)+c(\gamma_5))$$
Alternately, let $\gamma_2$ be a simple closed geodesic bounding a one-holed torus $T$ in $X$, and let $\gamma_1$ be a simple closed geodesic contained in $T$ such that $\gamma_1$ cuts $T$ into a pants $P$. Then there exists a simple closed geodesic $\zeta_1$ that is also contained in $T$, als cuts $T$ into a pants, and that intersects $\gamma_1$ in exactly one point. Furthermore, $\zeta_1$ satisfies
$$c(\zeta_1)\leq 10(c(\gamma_1)+c(\gamma_2))$$
 
\end{lemma}

\textbf{Proof.} Note that if one drops the combinatorial length bound condition, this result is standard.

In the four-holed-sphere case, $\gamma_2$ through $\gamma_5$ correspond to geodesics $\tilde{\gamma}_2$ through $\tilde{\gamma}_5$ on $\Gamma$ bounding some trivalent sub-metric ribbon graph $\tilde{S}\subset\Gamma$ homotopic to a four-holed sphere, and containing a geodesic $\tilde{\gamma}_1$ corresponding to $\gamma_1$. There are now a finite number of topological ribbon graph isomorphism types for $\tilde{S}$ - this fact is analogous to the fact that all trivalent ribbon graphs homotopic to three-holed spheres are topologically dumbbell graphs or theta graphs, as shown in Figure \ref{figureeight}. For each possible topological type of $\tilde{S}$, straightforward casework shows that there exists a simple closed geodesic $\tilde{\zeta}_1$ such that $\tilde{\zeta}_1$ intersects $\tilde{\gamma}_1$ twice, and $\tilde{\zeta}_1$ is constructed from explicit subarcs of $\tilde{\gamma}_1$, $\tilde{\gamma}_2$, $\tilde{\gamma}_3$, $\tilde{\gamma}_4$, and $\tilde{\gamma}_5$. These properties together imply that the corresponding geodesic $\zeta_1$ on $S$ intersects $\gamma_1$ twice and satisfies the combinatorial length bound.

\begin{figure}[h]

\centering
\includegraphics[scale=0.18]{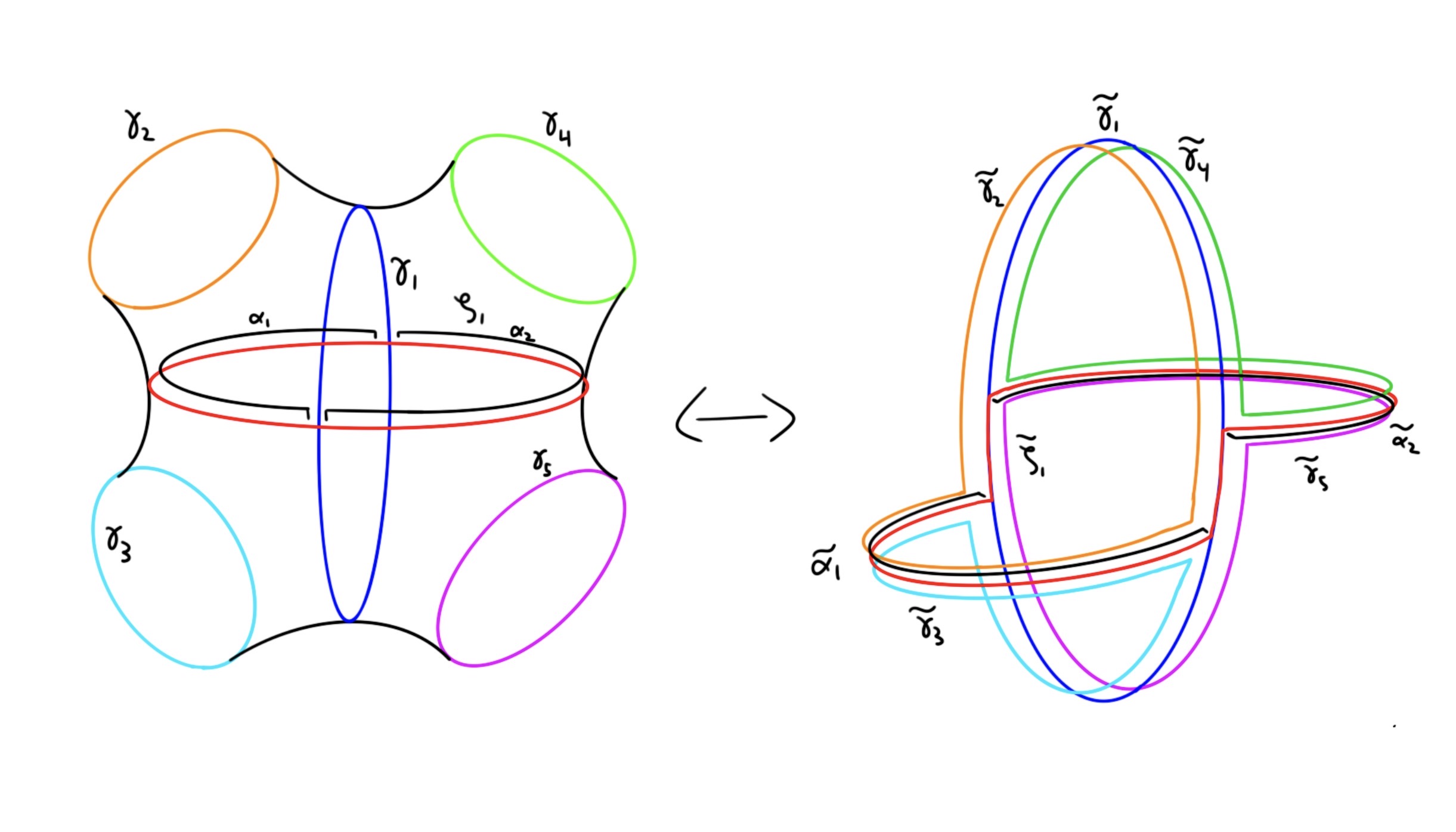}
\caption{An illustration of the argument used in Lemma \ref{geodflip} to argue that $c(\zeta_1)=O(1)$, in the case where both $P_1$ and $P_2$ correspond to theta graphs.}

\label{pantsflip}

\end{figure}

One possible case of this construction is illustrated in Figure \ref{pantsflip}.

The one-holed-torus case is analogous. $\blacksquare$

\

\begin{lemma}

\label{specialpants}

Fix $X\in\M_{g,n}(L_1,...,L_n)$ such that $\Gamma=\Phi(X)$ is trivalent, and assume the $L_i$ are $C$-controlled by $L$. Then there exist simple closed geodesics $\gamma_1,...,\gamma_{3g-3+n},\zeta_1,...,\zeta_{3g-3+n}$ on $X$ with the following properties.

\begin{enumerate}
\item All $\gamma_i$ and $\zeta_i$ respectively satisfy $l(\gamma_i)=O(L)$, $l(\zeta_i)=O(L)$.
\item All $\gamma_i$ and $\zeta_i$ respectively satisfy $c(\gamma_i)=O(1)$, $c(\zeta_i)=O(1)$.
\item The $\gamma_1,...,\gamma_{3g-3+n}$ form a pants decomposition.
\item The $\zeta_1,...,\zeta_{3g-3+n}$ form a pants decomposition.
\item For all $1\leq i,j\leq 3g-3+n$, $i(\gamma_i,\zeta_j)=O(1)$.
\item For all $1\leq i\leq 3g-3+n$, $i(\gamma_i,\zeta_i)=1$ or $i(\gamma_i,\zeta_i)=2$.
\item For all $1\leq j<i \leq 3g-3+n$, $i(\gamma_i,\zeta_j)=0$.
\end{enumerate}

\end{lemma}

\textbf{Proof.} First, note that there are only finitely many possible topological types for $\Gamma$, and in each type the combinatorial lengths of the boundary geodesics are fixed. Therefore, the combinatorial lengths of the boundary geodesics of $X$ are bounded in terms of $g$ and $n$.

Now, choose some edge $e$ of $\Gamma$, and find the corresponding pants as in Figure \ref{buildpants}. This pants has one or two boundary components that are geodesics in $\Gamma$, and we may cut along them to obtain a new metric ribbon graph with reduced topological complexity. Iterating this process yields some collection of geodesics $\tilde{\gamma}_1,...,\tilde{\gamma}_{3g-3+n}$ forming a pants decomposition on $\Gamma$, and transferring them to geodesics $\gamma_1,...,\gamma_{3g-3+n}$ on $X$ yields a pants decomposition on the surface. Furthermore, it is straightforward to verify via induction that this pants decomposition has bounded combinatorial lengths as in properties one and two above.

We now iteratively construct the $\zeta_i$. For $\zeta_1$, take $\gamma_1$ and examine the pants incident to it; we either obtain a single pants forming a one-holed torus or two pants forming a four-holed sphere. In either case, we may apply Lemma \ref{geodflip} to obtain $\zeta_1$. Notice that by the topological properties of $\zeta_1$, we have that $\zeta_1,\gamma_2,...,\gamma_{3g-3+n}$ forms a pants decomposition of $X$. Furthermore, the fact that all $\gamma_i$ and all boundary geodesics of $X$ have combinatorial length at most $O(1)$ implies that $c(\zeta_1)=O(1)$ as well.

We now iterate. At step $i$ of the iteration, we assume that we have geodesics $\zeta_1,...,\zeta_{i-1}$ such that $\zeta_1,...,\zeta_{i-1},\gamma_i,...,\gamma_{3g-3+n}$ form a pants decomposition and every one of these geodesics have combinatorial length at most $O(1)$. We then look at the pants incident to $\gamma_i$ in this decomposition and apply Lemma \ref{geodflip} to construct $\zeta_i$ as before, obtaining a new pants decomposition $\zeta_1,...,\zeta_i,\gamma_{i+1},...,\gamma_{3g-3+n}$. In particular, $\zeta_i$ is disjoint from $\gamma_{i+1}$ through $\gamma_{3g-3+n}$, and the property that $c(\zeta_i)=O(1)$ follows via the same argument used for $c(\zeta_1)$ above and the property that $c(\zeta_1)$ through $c(\zeta_{i-1})$ are $O(1)$.

Upon completion of the iteration, all $\gamma_i$ and $\zeta_i$ have combinatorial length $O(1)$, and properties three, four, six, and seven follow immediately from the construction. Since each edge of $\Gamma$ has length at most $O(L)$ by the assumption that the boundary lengths $L_i$ are $C$-controlled by $L$, the combinatorial length bound implies that all $\tilde{\gamma}_i$ and $\tilde{\zeta}_i$ have lengths at most $O(L)$. Corollary \ref{ribcontrol} then implies that all $\gamma_i$ and $\zeta_i$ have lengths at most $O(L)$.

Lastly, the fact that all $\gamma_i$ and $\zeta_i$ have combinatorial lengths at most $O(1)$ implies that they have pairwise intersection number at most $O(1)$. $\blacksquare$

\

\begin{figure}[h]

\centering
\includegraphics[scale=0.2]{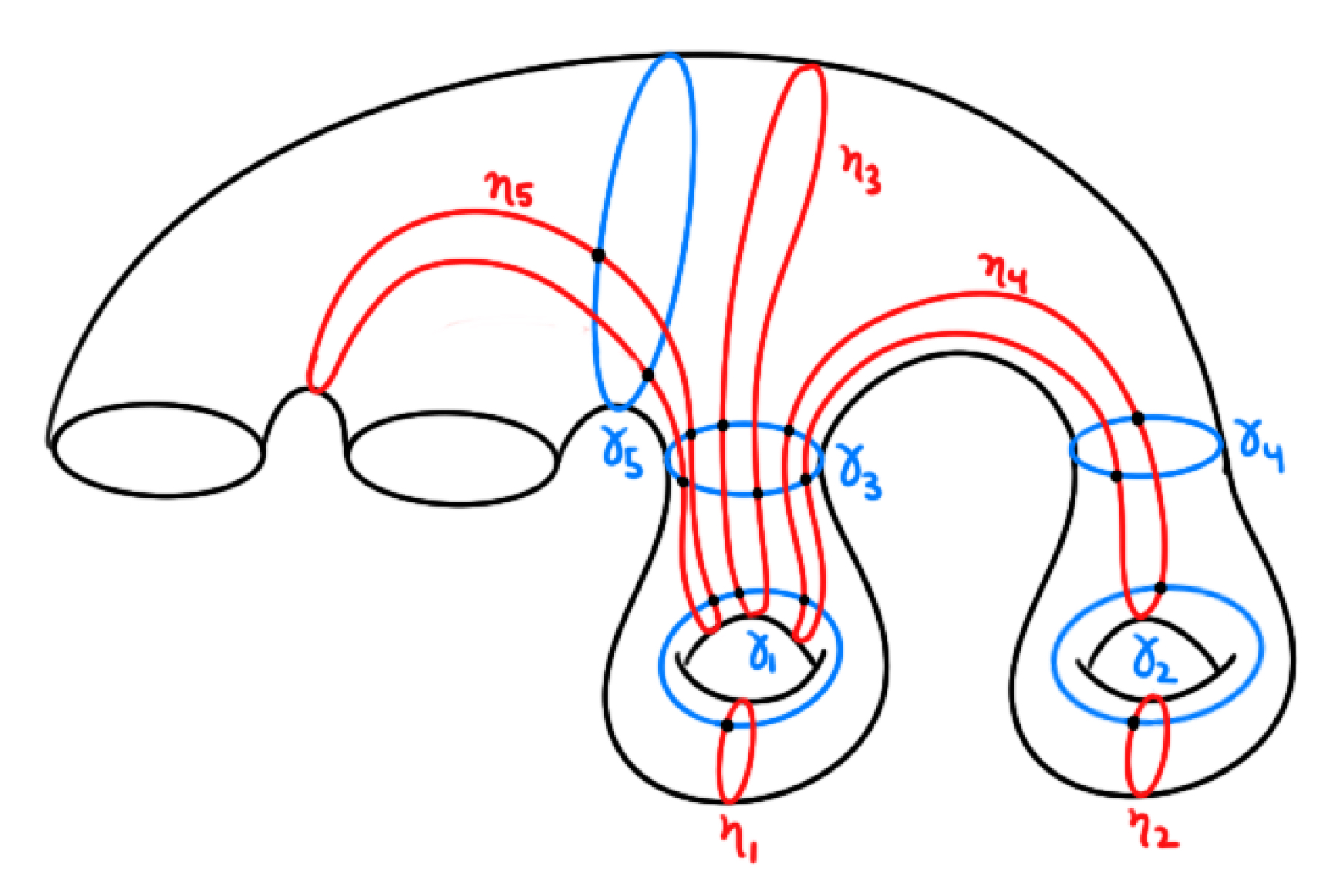}
\caption{An example of a special double pants decomposition on $S_{2,2}$.}

\label{specialdouble}

\end{figure}

Such double pants decompositions will be called \emph{special double pants decompositions}. See Figure \ref{specialdouble} for an example of a special double pants decomposition on a surface.

\begin{lemma}
\label{goodgeodframe}

Let $X\in \M_{g,n}(L_1,...,L_n)$ be a $(B,A)$-good surface. Let $\Gamma=\Phi(X)$. There exist simple (multi-)geodesics $\gamma_1,...,\gamma_{6g-6+2n}$ on $X$ and corresponding geodesics $\tilde{\gamma}_1,...,\tilde{\gamma}_{6g-6+2n}$ on $\Gamma$ with the following properties:
\begin{enumerate}
\item $\gamma_1,...,\gamma_{6g-6+2n}$ is a good geodesic frame.
\item We have $i(\gamma_i,\gamma_j)=O(1)$ for all $1\leq i\neq j\leq 6g-6+2n$.
\item We have $c(\gamma_i)=O(1)$ for all $1\leq i\leq 6g-6+2n$.
\item If $k_i=l(\tilde{\gamma}_i)$, then the $k_i$ form local coordinates in a neighborhood of $\Gamma$.
\end{enumerate}
	
\end{lemma}

\textbf{Proof.} Using Lemma \ref{specialpants}, let $\gamma_1,...,\gamma_{3g-3+n},\zeta_1,...,\zeta_{3g-3+n}$ be a special double pants decomposition on $X$; these curves have combinatorial lengths and pairwise intersection numbers at most $O(1)$.

We will iteratively modify the $\zeta_i$ using Dehn twists. For each $\gamma_i$, define
$$m_i=\min_{\substack{1\leq j\leq 3g-3+n\\ p\in\gamma_i\cap\zeta_j}}f_{\gamma_i,\zeta_j,p}(0)$$
using functions $f_{\gamma_i,\zeta_j,p}$ construction as in Lemma \ref{shiftvalue}.

By Lemma \ref{functionlengthbound} and properties one and two of Lemma \ref{specialpants}, we have, for all $\gamma_i$, $\zeta_j$, and $p\in\gamma_i\cap\zeta_j$,
$$|f_{\gamma_i,\zeta_j,p}(0)|=O(l(\gamma_i))$$

Therefore, all $m_i$ are bounded below by $-O(l(\gamma_i))$. We will now obtain new simple closed geodesics $\zeta_1',...,\zeta'_{3g-3+n}$ by applying $\frac{3m_i}{l(\gamma_i)}+2$ positive Dehn twists around each $\gamma_i$ to the collection $\zeta_1,...,\zeta_{3g-3+n}$. For each $\gamma_i$, $\zeta_j$, and $p\in\gamma_i\cap\zeta_j$, let $p'$ be the corresponding intersection point in $\gamma_i\cap \zeta_j'$. Fixing $\gamma_i$, $\zeta_j$, and $p\in\gamma_i\cap\zeta_j$, by Lemma \ref{shiftvalue}, 
$$f_{\gamma_i,\zeta_j',p'}(0)\geq f_{\gamma_i,\zeta_j,p}(0)+m_i+2l(\gamma_i)\geq 2l(\gamma_i)$$
By Lemma \ref{twistlengthincrease},
$$c(\zeta_j')\leq c(\zeta_j)+\sum_{i=1}^{3g-3+n}i(\gamma_i,\zeta_j)c(\gamma_i)\left(\frac{3m_i}{l(\gamma_i)}+2\right)=O(1)$$
Lastly, $i(\gamma_i,\zeta'_j)=i(\gamma_i,\zeta_j)=O(1)$ for all $1\leq i,j\leq 3g-3+n$, while the $\zeta_j'$ form a pants decomposition.

Set $\gamma_{3g-3+n+k}=\zeta'_k$ for $1\leq k\leq 3g-3+n$, obtaining $\gamma_1,...,\gamma_{6g-6+2n}$. From the work carried out above, the desired second and third properties are satisfied. Furthermore, since all $f_{\gamma_i,\zeta_j',p'}(0)$ are bounded below by $2l(\gamma_i)$, Lemma \ref{functionbadbound} implies there are no bad intersection points, so $\gamma_1,...,\gamma_{6g-6+2n}$ form a good geodesic frame. In particular, the fact that all $f_{\gamma_i,\zeta'_j,p'}$ are large and positive implies that all signed crossing numbers are $+1$, so for all $\gamma_i$ and $\zeta_j'$, $si(\gamma_i,\zeta_j')=i(\gamma_i,\zeta_j')$.

Transfer $\gamma_1,...,\gamma_{6g-6+2n}$ to geodesics $\tilde{\gamma}_1,...,\tilde{\gamma}_{6g-6+2n}$ on $\Gamma$, and define functions $k_1,...,k_{6g-6+2n}$ via $k_i=l(\tilde{\gamma}_i)$. Using the pants decomposition $\tilde{\gamma}_1,...,\tilde{\gamma}_{3g-3+n}$, construct combinatorial Fenchel-Nielsen coordinates $k_1,...,k_{3g-3+n},$ $\tilde{\tau}_1,...,\tilde{\tau}_{3g-3+n}$. We will analayze the Jacobian matrix $J$ obtained from the map
$$(k_1,...,k_{3g-3+n},\tilde{\tau}_1,...,\tilde{\tau}_{3g-3+n})\to(k_1,...,k_{3g-3+n},k_{3g-2+n},...,k_{6g-6+2n})$$
We will split $J$ into four $(3g-3+n)\times (3g-3+n)$ submatrices. Of course, $\frac{\partial k_i}{\partial k_j}=\delta_{ij}$ for $1\leq i,j\leq 3g-3+n$, so the upper left submatrix is $I_{3g-3+2n}$. Likewise, $\frac{\partial k_i}{\partial \tau_j}=0$ for $1\leq i,j\leq 3g-3+n$, so the upper right submatrix is $0$. By Lemma \ref{graphtwist}, 
$$\frac{\partial k_i}{\partial \tau_j}=si(\gamma_j,\zeta_{i-(3g-3+n)}')=i(\gamma_j,\zeta'_{i-(3g-3+n)})$$ for $1\leq i,j\leq 3g-3+n$, and by properties six and seven of Lemma \ref{specialpants}, we may conclude the lower right submatrix of $J$ is lower triangular with all nonzero entries on the diagonal. Therefore, the block structure of $J$ implies $\det(J)\neq 0$.

In \cite{ABCGLW21}, the authors establish that, for any choice of Fenchel-Nielsen coordinates and on any fixed maximal open cell of $MRG_{g,n}(L_1,...,L_n)$, the edge length functions are linear functions of the Fenchel-Nielsen coordinates (see the discussion below Theorem 2.41). In turn, $k_1,...,k_{6g-6+2n}$ are clearly linear functions of the edge lengths. So the prospective change of coordinates given above is linear in a neighborhood of $\Gamma$, and the nonsingularity of $J$ implies the $k_i$ form local coordinates in a neighborhood of $\Gamma$ that are analytic with respect to combinatorial Fenchel-Nielsen coordinates. $\blacksquare$

\subsection{Intersection Matrix Convergence}

The following lemma will allow us to choose geodesics (and hence choose local coordinates) with all the properties we will later need.

\begin{lemma}
\label{geodsetup}

Let $X\in\M_{g,n}(L_1,...,L_n)$ be a $(B,A)$-good surface with $B$ a sufficiently large multiple of $A$ in terms of $g$ and $n$, and let $\Gamma=\Phi(X)$. There exist geodesics $\gamma_1,...,\gamma_{6g-6+2n}$ on $X$ and corresponding geodesics $\tilde{\gamma}_1,...,\tilde{\gamma}_{6g-6+2n}$ on $\Gamma$ such that if the matrix $W$ is defined by $W_{ij}=w(\gamma_i,\gamma_j)$, the matrix $S$ is defined by $S_{ij}=iw(\gamma_i,\gamma_j)$, $l_i=l(\gamma_i)$, and $k_i=l(\tilde{\gamma}_i)$, then the following holds: $W$ and $S$ are invertible, the $l_i$ are a local coordinate system for $\M_{g,n}(L_1,...,L_n)$ around $X$, the $k_i$ are a local coordinate system for $MRG_{g,n}(L_1,...,L_n)$ around $\Gamma$, formulas \ref{ewpform} and \ref{ekontform} hold, the $\gamma_i$ have $O(1)$-bounded combinatorial length, and 
$$\det(W)=\det(S)+O(e^{-A})$$
$$\det(W^{-1})=\det(S^{-1})+O(e^{-A})$$
$$|\det(W^{-1})|,|\det(S^{-1})|=\Omega(1)$$
$$\max_{i,j}(|W_{ij}|),\max_{i,j}(|W_{ij}^{-1}|),\max_{i,j}(|S_{ij}|),\max_{i,j}(|S^{-1}_{ij}|)=O(1)$$
If $X$ is instead $\eta$-great, the same result holds, except with adjusted errors:
$$\det(W)=\det(S)+O(e^{-2\eta})$$
$$\det(W^{-1})=\det(S^{-1})+O(e^{-2\eta})$$

\end{lemma}

\textbf{Proof.} Assume $X$ is $\eta$-great. Choose $\tilde{\gamma}_1,...,\tilde{\gamma}_{6g-6+2n}$ via Lemma \ref{goodgeodframe}, and define $\gamma_1,...,\gamma_{6g-6+2n}$ by correspondence. By Lemma \ref{greatanglecontrol}, we have
$$W_{ij}=S_{ij}+O(i(\gamma_i,\gamma_j)e^{-2\eta})$$
and Lemma \ref{geodchoice} reduces this formula to
$$W_{ij}=S_{ij}+O(e^{-2\eta})$$

By Lemma \ref{geodchoice}, the entries of $S$ (and hence the entries of $W$) are bounded, while since the $k_i$ form a local coordinate system we know $S$ is invertible. Since furthermore the entries of $S$ are integers, $|\det(S^{-1})|=\Omega(1)$. By continuity of the determinant, as long as $A$ is sufficiently large in terms of $g$ and $n$, $\eta\geq A$ allows us to conclude $|\det(W^{-1})|=\Omega(1)$ as well. The adjoint formula for the determinant combined with the formulas given above now yields the desired result. In particular, $W$ is invertible, so Theorem \ref{cosineformula} applies. 

For the claim for good surfaces, we apply Lemma \ref{supersectorintersectioncontrol} and then follow essentially the same argument. In particular, since the $k_i$ form a local coordinate system at $\Gamma$, $S$ must be invertible. We obtain
$$W_{ij}=S_{ij}+O(i(\gamma_i,\gamma_j)e^{-2A})$$ 
and now need to choose $A$ sufficiently large that $e^{-2A}$ compensates for $i(\gamma_i,\gamma_j)$ and the overall error is small enough to conclude $|\det(W^{-1})|=\Omega(1)$. Since the bound on $i(\gamma_i,\gamma_j)$ depends only on $g$ and $n$ and $A$ is allowed to depend on these parameters, this can be done. $\blacksquare$

\section{Quantitative Convergence: Derivatives}

\subsection{Twist Derivatives}

Given a surface $X\in\M_{g,n}(L_1,...,L_n)$ such that $\Phi(X)=\Gamma$ is trivalent, choose two simple closed geodesics $\gamma$ and $\zeta$. Let $\tilde{\gamma}$ be the geodesic on $\Gamma$ corresponding to $\gamma$, let $k=l(\tilde{\gamma})$, and let $\frac{\partial}{\partial \tau}$ be the twist derivative corresponding to $\gamma$. Our goal is to compute $\frac{\partial \Phi^*k}{\partial \tau}$. We will generally drop the pullback notation, and so write $k$ for $\Phi^*k$ and $\frac{\partial k}{\partial \tau}$ for $\frac{\partial \Phi^*k}{\partial \tau}$.

\begin{figure}[h]

\centering
\includegraphics[scale=0.15]{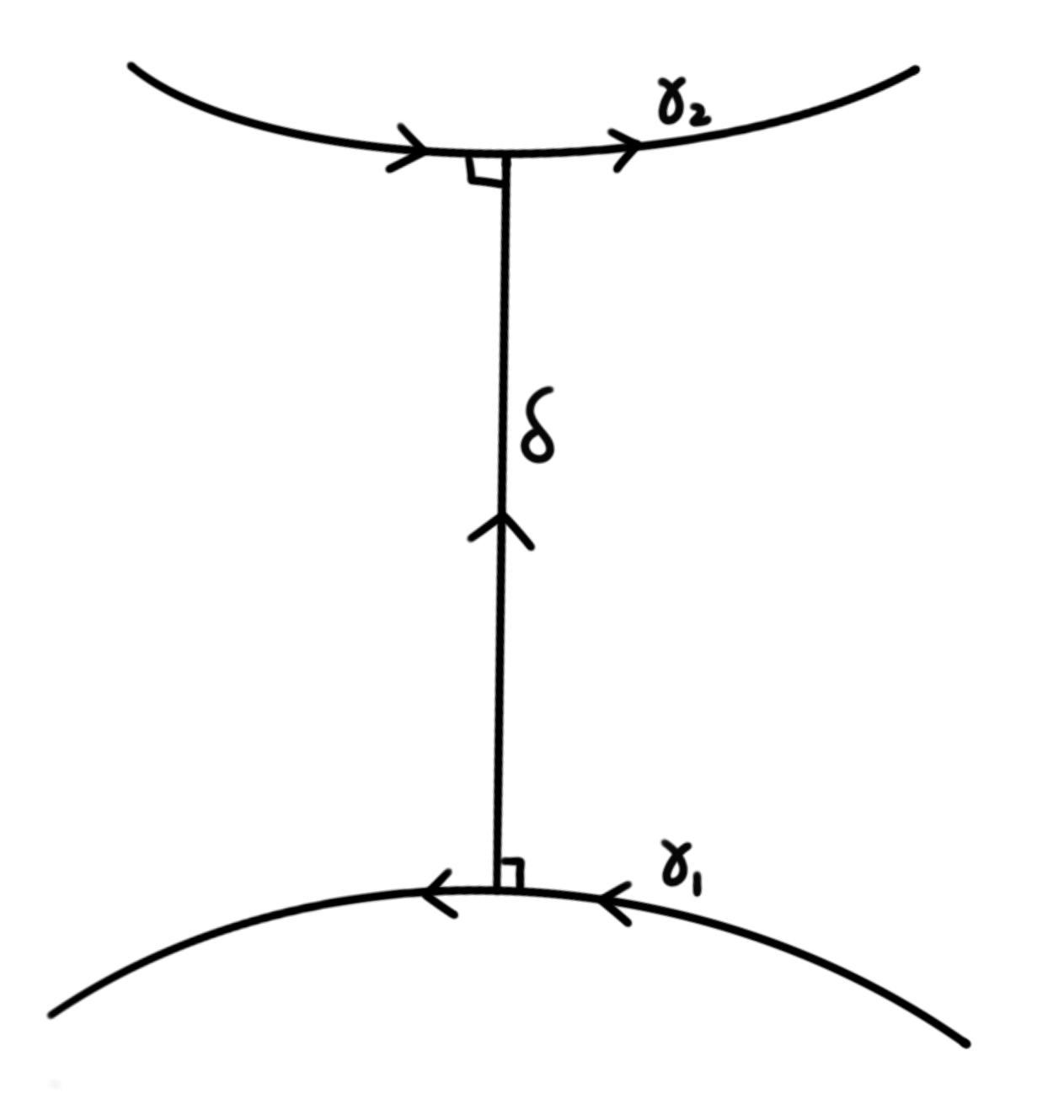}

\caption{The orientation scheme used to define signed distance in Theorem \ref{mondello}.}

\label{signeddistorient}

\end{figure}

\begin{theorem}[Mondello, \cite{GM09}, Thm. 3.7]
\label{mondello}

Let $X$ be a hyperbolic surface without boundary, and fix two simple closed geodesics $\gamma_1$ and $\gamma_2$ on $X$. Let $\delta$ be a simple geodesic arc from $y_1\in\gamma_1$ to $y_2\in\gamma_2$ realizing the distance between $\gamma_1$ and $\gamma_2$ in its homotopy class. Let $\xi$ be a third simple closed geodesic that does not intersect $\delta$, and let $\frac{\partial}{\partial\tau}$ be the (positive) twist derivative corresponding to $\xi$. Orient $\delta$ from $y_1$ to $y_2$, and use this orientation to orient $\gamma_1$ and $\gamma_2$ as in Figure \ref{signeddistorient}.

Then
$$\frac{\partial l(\delta)}{\partial\tau}=c_1+c_2$$
where $c_i=\sum_{x_j\in\xi\cap\gamma_i}c_i(x_j)$ and
$$c_i(x_j)=\frac{\sinh(l(\gamma_i)/2-d_s(y_i,x_j))}{2\sinh(l(\gamma_i)/2)}\sin(\theta_{x_j})$$
Here, $d_s(y_i,x_j)$ is the signed distance from $y_i$ to $x_j$, traveling around $\gamma_i$ with respect to its orientation. Furthermore, $\theta_{x_j}$ is the intersection angle of $\gamma_i$ and $\xi$ at $x_j$, measured counterclockwise from $\xi$ to $\gamma_i$. If one wishes to instead take a negative twist derivative, the formula holds with all signs reversed.
	
\end{theorem}

\begin{figure}[h]

\centering
\includegraphics[scale=0.15]{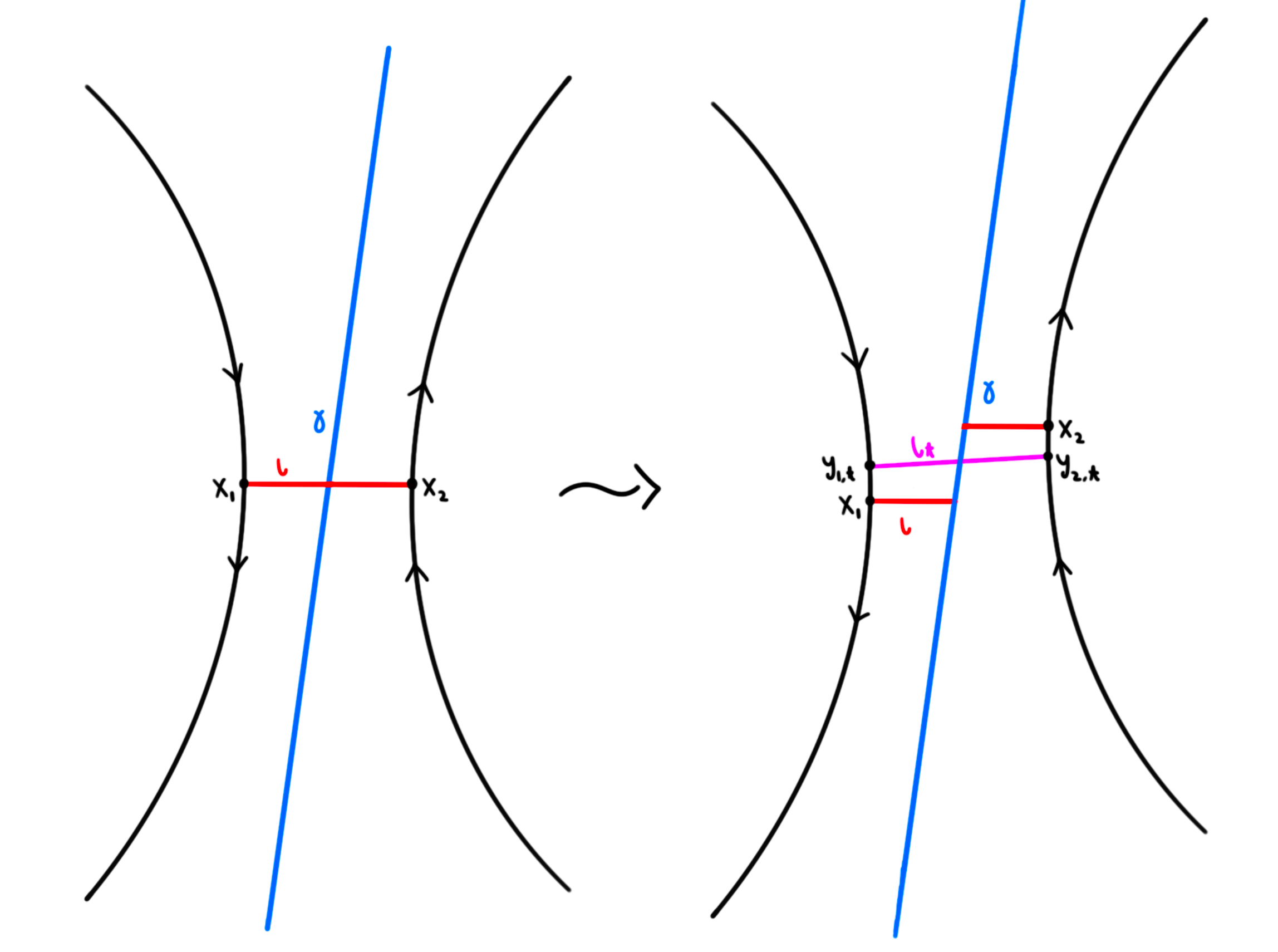}
\caption{An illustration of the setup of Corollary \ref{trianglecomp}, showing the new intercostal obtained after twisting and the distances being measured.}

\label{twistchange}

\end{figure}

We apply this theorem as follows:

\begin{corollary}
\label{trianglecomp}

Fix $X\in\M_{g,n}(L_1,...,L_n)$. Let $\iota\subset X$ be an intercostal, i.e. a simple path between two points on $\partial X$ that is length-minimizing in its homotopy class. Label the endpoints of $\iota$ as $x_1$ and $x_2$.

Let $\gamma$ be a simple closed geodesic on $X$, and let $\frac{\partial}{\partial \tau}$ be the twist derivative along $\gamma$. For $t$ in some sufficiently small neighborhood $(-\epsilon,\epsilon)$ let $X_t$ be the surface obtained by a twist along $\gamma$ of length $t$, and let $\iota_t$ be the intercostal homotopic to $\iota$ on $X_t$. Label the endpoints of $\iota_t$ as $y_{1,t}$ and $y_{2,t}$, so that $y_{i,t}$ lies on the same side as $x_{i,t}$. With this setup, $\frac{\partial}{\partial \tau}=\frac{\partial}{\partial t}$.

Using the orientation on $X$, we have an induced orientation on $\partial X$. Then $y_{i,t}$ lies `behind' $x_i$ with respect to this orientation in the sense of Figure \ref{twistchange}, and for $i=1,2$,
$$\frac{\partial d(x_i,y_{i,t})}{\partial t}\Big|_{t=0}=\sum_{p\in \gamma\cap\iota}\frac{\sinh(l(\iota)-d(x_i,p))}{\sinh(l(\iota))}\sin(\theta_p)$$
where $\theta_p$ is the angle between $\gamma$ and $\iota$ at $p$, measured counterclockwise from $\gamma$ to $\iota$.
	
\end{corollary}

\textbf{Proof.} We will prove the formula for $x_2$. Lift $\iota$ and $x_2$ to the universal cover $\tilde{X}$ of $X$, obtaining $\tilde{\iota}$ and $\tilde{x}_2$.\footnote{Notice that here tildes refer to objects on the universal cover $\tilde{X}$, while elsewhere in the paper tildes have generally referred to objects on a metric ribbon graph.} Let $\beta$ be the boundary geodesic of $X$ containing $x_2$, and let $\tilde{\beta}$ be its lift to $\tilde{X}$ containing $\tilde{x}_2$; it has an induced orientation.  Each $p\in\gamma\cap\iota$ lifts to a unique $\tilde{p}\in\tilde{\iota}$, and there is a unique lift $\tilde{\gamma}_p$ of $\gamma$ passing through this $\tilde{p}$. There is a finite number of these $\tilde{\gamma}_p$, and they are disjoint. As a consequence of this finiteness, there exists an intercostal $\tilde{\iota}'\subset \tilde{X}$ with an endpoint $\tilde{x}_2'$, satisfying the following properties:
\begin{enumerate}
\item 	The point $\tilde{x}'_2$ lies `ahead' of $\tilde{x}_2$ with respect to the orientation on $\tilde{\beta}$.
\item No $\tilde{\gamma}_p$ intersects $\tilde{\iota}'$.
\end{enumerate}
To obtain $\tilde{\iota}$, first note that $\tilde{\beta}$ contains infinitely many lifts of $x_2$ incident to lifts of $\iota$, spaced out along $\tilde{\beta}$ at interval $l(\beta)$. As $\tilde{\beta}$ is disjoint from any $\tilde{\gamma}_p$, the ideal endpoints of $\tilde{\beta}$ are distinct from the set of ideal endpoints of the $\tilde{\gamma}_p$. In particular, we can find lifts of $\iota$ incident to $\tilde{\beta}$ that are arbitrarily far from any $\tilde{\gamma}_p$, and can find them in either direction from $\tilde{x}_2$. As the lifts of $\iota$ have finite length, we can certainly find some $\tilde{\iota}'$ ahead of $\tilde{x}_2$ that is disjoint from all $\tilde{\gamma}_p$.

Now, double $\tilde{X}$ to create a new boundaryless surface $\tilde{Y}$. The intercostals $\tilde{\iota}$ and $\tilde{\iota}'$ extend to respective closed geodesics $\tilde{\zeta}$ and $\tilde{\zeta'}$ on $Y$. Also, we may mirror each $\tilde{\gamma}_p$ to a new geodesic $\tilde{\gamma}_p^Y$ on the other `side' of $Y$, intersecting $\tilde{\zeta}$ at some point $\tilde{p}'$.  Letting $\delta$ be the segment between $\tilde{x}_2$ and $\tilde{x}_2'$ along $\tilde{\beta}$, $\delta$ is a geodesic arc on $Y$ realizing the distance between $\tilde{\zeta}$ and $\tilde{\zeta}'$.

Consider a $t$-length positive twist along each $\tilde{\gamma}_p$ and an $t$-length negative twist along each $\tilde{\gamma}^Y_p$ for some sufficiently small $t>0$. This will produce a new geodesic $\tilde{\zeta}_t$ homotopic to $\tilde{\zeta}$, which restricts on $X$ to a new intercostal $\tilde{\iota}_t$ homotopic to $\tilde{\iota}$ by symmetry of the twists. By choice of $\tilde{\iota}'$, $\tilde{\zeta}'$ and $\tilde{\iota}'$ are unchanged under twisting. Furthermore, the endpoint of $\tilde{\iota}_t$ lying on $\beta_t$ projects to $y_{2,t}$, and we have some new $\delta_t$ between $\tilde{\iota}_t$ and $\tilde{\iota}'$. Define $l_\delta(t)$ to be the length of $\delta_t$.

We are now in a position to apply Theorem \ref{mondello} to calculate $\frac{\partial l_\delta}{\partial t}|_{t=0}$. While we are twisting along multiple geodesics simultaneously, linearity of tangents means we may consider their individual contributions and sum the results. By construction, each $\tilde{\gamma}_p$ or $\tilde{\gamma}_p^Y$ intersects $\tilde{\zeta}$ once, and does not intersect $\tilde{\zeta}'$. While the $\tilde{\gamma}_p^Y$ are twisted negatively, the resulting negative sign ends up cancelling when
$$\sinh(l(\tilde{\zeta})/2-d_s(y_i,\tilde{p}'))=-\sinh(l(\tilde{\zeta})/2-d(y_i,\tilde{p}))$$
is taken into account.

Applying Theorem \ref{mondello}, we obtain
$$\frac{\partial l_\delta}{\partial t}|_{t=0}=\sum_{p\in\gamma\cap\iota}\frac{\sinh(l(\tilde{\zeta})/2-d(\tilde{x}_2,\tilde{p}))}{\sinh(l(\tilde{\zeta})/2)}\sin(\theta_p)$$
We now observe that $l_\delta(t)$ is the distance from $y_{2,t}$ to $\tilde{x}_2$ plus the distance from $\tilde{x}_2$ to $\tilde{x}_2'$, and this second distance is invariant under twisting. After projecting back down, we obtain the desired formula. $\blacksquare$

\

\begin{lemma}
\label{derivcalc}

Fix $X\in\M_{g,n}(L_1,...,L_n)$ such that $X$ is $(B,A)$-good, and consider two simple closed geodesics $\gamma$ and $\zeta$. Let $\tilde{\zeta}$ be the geodesic on $\Gamma=\Phi(X)$ corresponding to $\zeta$, and define $k=l(\tilde{\zeta})$. Let $\frac{\partial}{\partial\tau}$ be the twist derivative corresponding to $\gamma$. Then
$$\frac{\partial k}{\partial \tau}=O(c(\gamma)c(\zeta))$$
Furthermore, if $X$ is in fact $\eta$-great, then
$$\frac{\partial k}{\partial\tau}=iw(\zeta,\gamma)+O(e^{-\eta/2}c(\gamma)c(\zeta)))$$

\end{lemma}

\textbf{Proof.} Decompose $X$ into sectors, and cut $\zeta$ along all intercostals it intersects. This splits $\zeta$ into a collection of arcs $\zeta_i$ indexed by some set $I$, such that each $\zeta_i$ traverses a single sector. Within that sector, the two intercostals incident to the endpoints of $\zeta_i$ are in turn incident to a single boundary arc, which we denote $b_i$. We then have
$$k:=l(\tilde{\zeta})=\sum_{i\in I}l(b_i)$$
This is straightforward if all corridors of $X$ contain their intercostals, since then each $l(b_i)$ is the length of two adjacent half-edges on $\Gamma$. Some additional casework handles the cases where some corridors of $X$ do not contain their intercostals.

\begin{figure}

\centering
\includegraphics[scale=0.2]{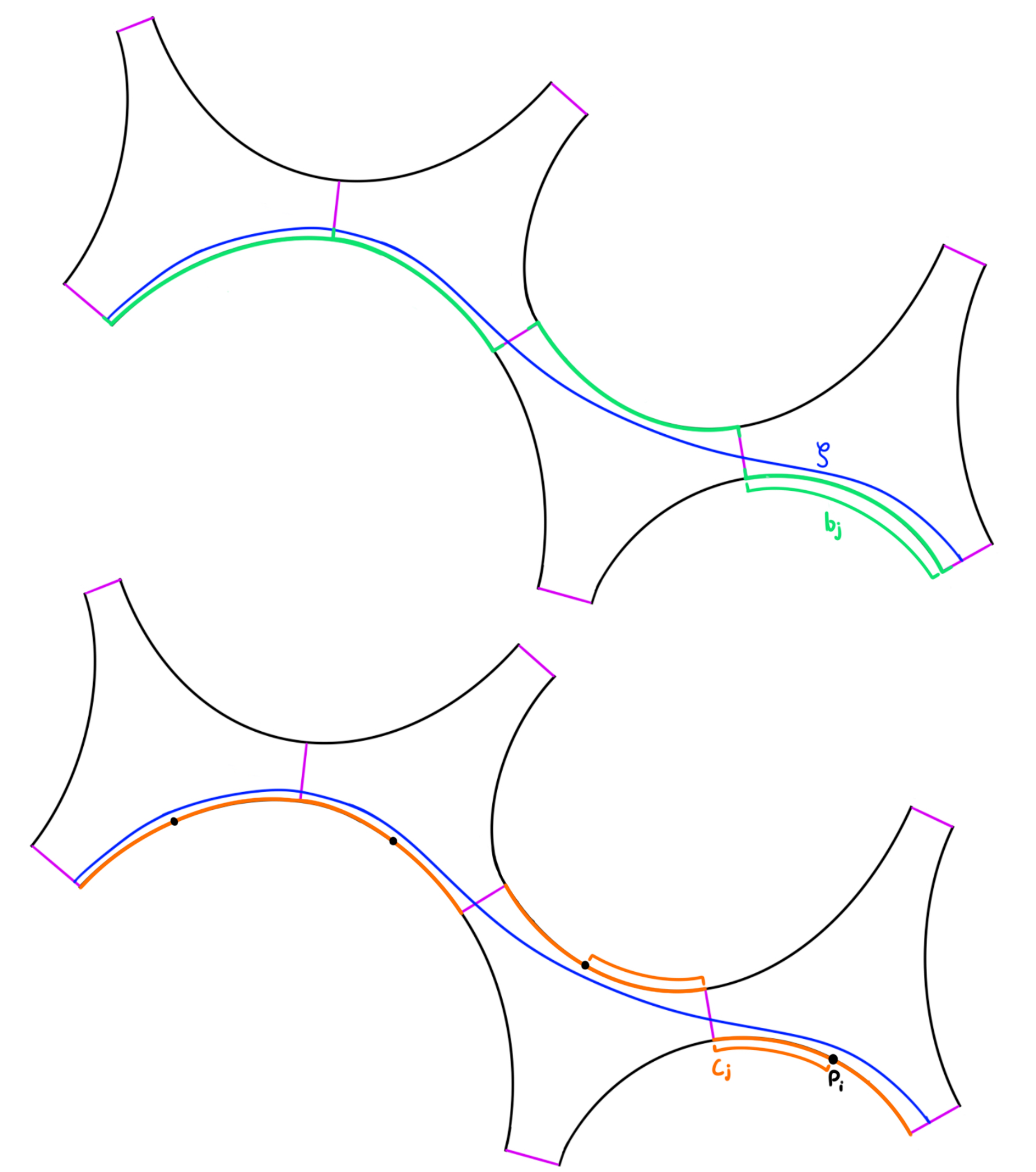}
\caption{An illustration of the decomposition procedure used to obtain the $b_i$ and $c_j$ arcs used in Lemma \ref{derivcalc}.}

\label{cutandglue}

\end{figure}

Let $p_i$ be the midpoint of $b_i$ (the choice of midpoint is arbitrary - any interior point would work). Now, each $\zeta_i$ is incident to two other arcs $\zeta_i'$ and $\zeta_i''$ at its endpoints. Cutting all $b_i$ at the chosen $p_i$, we may use this incidence to reglue halves of the $b_i$ together, resulting in a collection of $c_j$ indexed by some set $J$ with $|J|=|I|$. Each $c_j$ is either a single arc between two $p_i$, or a union of two arcs that are incident to a common intercostal. See Figure \ref{cutandglue} for an illustration of this procedure. We observe that the first case represents a `no hop' as in Figure \ref{hopdef}, while the second case represents a left hop or a right hop.

Since all we have done is redistribute lengths,
$$l(\tilde{\zeta})=\sum_{j\in J}l(c_j)$$
It therefore suffices to calculate $\frac{\partial l(c_j)}{\partial \tau}$ for each $c_j$. If $c_j$ represents a no hop, then since twisting does not move the $p_i$, we have
$$\frac{\partial l(c_j)}{\partial\tau}=0$$

\begin{figure}[h]

\centering
\includegraphics[scale=0.18]{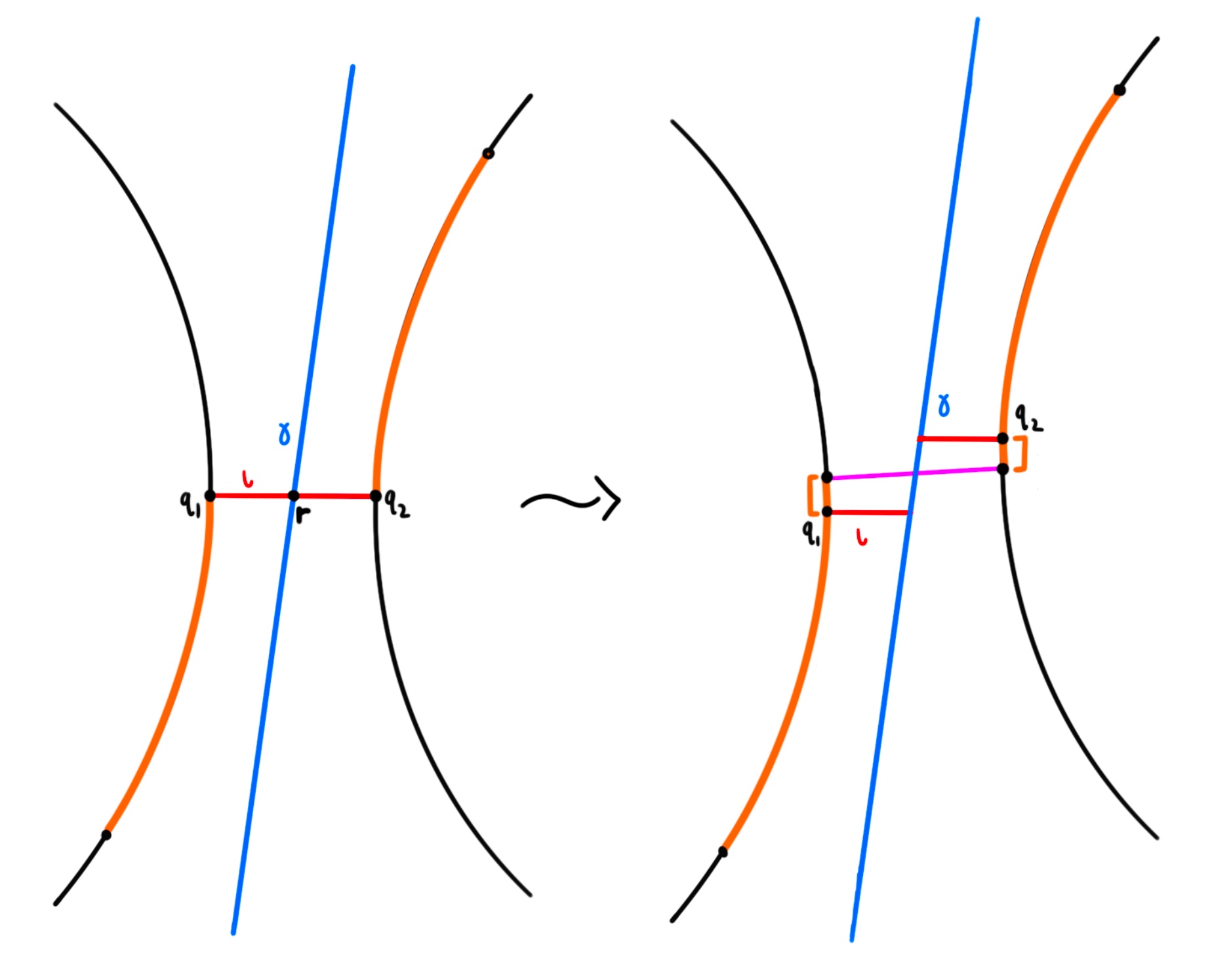}
\caption{An illustration of the additional length created for a $c_j$ representing a `right hop' under twisting $\gamma$. The additional length is marked in brackets; we wish to understand the infinitesimal increase in this length.}

\label{righthoptwist}

\end{figure}

On the other hand, assume $c_j$ represents a right hop. Then an infinitesimal twist along $\gamma$ causes an infinitesimal increase in the lengths of both components of $c_j$. Let $\iota_j$ be the intercostal incident to both components of $c_j$, and let $q_1$ and $q_2$ be its endpoints. We may recognize that Corollary \ref{trianglecomp} applies to this setup - each component of $c_j$ has an endpoint on $\iota_j$ and an endpoint that is fixed under twisting, so the change in its length comes solely from the new length created `behind' the endpoint on $\iota_j$. See Figure \ref{righthoptwist} for an illustration. 

We obtain the following formula:
$$\frac{\partial l(c_j)}{\partial \tau}=\sum_{r\in \gamma\cap\iota_j}\left(\frac{\sinh(l(\iota_j)-d(q_1,r))}{\sinh(l(\iota_j))}\sin(\theta_r)+\frac{\sinh(l(\iota_j)-d(q_2,r))}{\sinh(l(\iota_j))}\sin(\theta_r)\right)$$
All terms in this sum are bounded in absolute value by $2$, and the same holds for the case of a left hop. Furthermore, $|\gamma\cap\iota_j|$ is bounded above by $c(\gamma)$. In turn $|J|=|I|=c(\zeta)$, so we obtain
$$\frac{\partial k}{\partial \tau}=O(c(\gamma)c(\zeta))$$
We now assume that $X$ is $\eta$-great. By Lemma \ref{intercostalintersectioncontrol}, we may replace all $\sin(\theta_r)$ terms with $1+O(e^{-\eta/2})$. Distributing the errors gives a bounded number of error terms that are each at most $O(e^{-\eta/2}c(\gamma)c(\zeta)))$, along with a main term
$$\sum_{r\in\gamma\cap\iota_j}\left(\frac{\sinh(l(\iota_j)-d(q_1,r))}{\sinh(l(\iota_j))}+\frac{\sinh(l(\iota_j)-d(q_2,r))}{\sinh(l(\iota_j))}\right)$$
By Lemma \ref{intercostalcontrol}, $l(\iota_j)=O(e^{-\eta/2})$. Combined with the power-series estimate $\sinh(x)=x+O(x^3)$ and the fact that $d(q_1,r)+d(q_2,r)=l(\iota_j)$, we obtain
$$\frac{\sinh(l(\iota_j)-d(q_1,r))}{\sinh(l(\iota_j))}+\frac{\sinh(l(\iota_j)-d(q_2,r))}{\sinh(l(\iota_j))}=1+O(e^{-\eta/2})$$
Therefore, if $c_j$ represents a right hop, it contributes $|\gamma\cap\iota_j|+O(e^{-\eta/2}c(\gamma))$. Likewise, if $c_j$ represents a left hop, it contributes $-|\gamma\cap\iota_j|+O(e^{-\eta/2}c(\gamma))$. Letting $h_j=0$ if $c_j$ represents no hop, $h_j=1$ if $c_j$ represents a right hop, and $h_j=-1$ if $c_j$ represents a left hop, we end up with the formula
$$\frac{\partial l(\tilde{\zeta})}{\partial \tau}=\sum_{j\in J}\left(h_j|\gamma\cap\iota_j|+O(e^{-\eta/2}c(\gamma))\right)$$
$$=\sum_{j\in J}h_j|\gamma\cap \iota_j|+O(e^{-\eta/2}c(\gamma)c(\zeta))$$
Lastly, by Lemma \ref{intersectiondef}, 
$$\sum_{j\in J}h_j|\gamma\cap\iota_j|=iw(\zeta,\gamma)$$
This proves the desired formula. $\blacksquare$

\subsection{Length Derivatives from Twist Derivatives}

Our goal is to quantitatively compare the Kontsevich and Weil-Petersson forms. The following lemma is a reformulation of Lemma \ref{earlyformratio}.

\begin{lemma}
\label{ratiocontrol}

Fix a $(B,A)$-good surface $X\in \M_{g,n}(L_1,...,L_n)$, and let $\Gamma\in MRG_{g,n}(L_1,...,L_n)$. Let $dV_{WP}$ be the Weil-Petersson form, and let $dV_K$ be the Kontsevich form. Recall that $\Phi: \M_{g,n}(L_1,...,L_n)\to MRG_{g,n}(L_1,...,L_n)$ is the spine graph map. Then
\begin{equation}\label{radon1}\frac{\Phi^*dV_{K}}{dV_{WP}}|_X=O(1)\end{equation}
Moreover, if $X$ is in fact $\eta$-great, then
\begin{equation}\label{radon2}\frac{\Phi^*dV_{K}}{dV_{WP}}|_X=1+O(e^{-\eta/2})\end{equation}	
	
\end{lemma}

\textbf{Proof.} We first apply Lemma \ref{geodsetup} to find geodesics $\gamma_1,...,\gamma_{6g-6+2n}$ on $X$ and corresponding geodesics $\tilde{\gamma}_1,...,\tilde{\gamma}_{6g-6+2n}$, and define $l_i=l(\gamma_i)$, $k_i=l(\tilde{\gamma}_i)$. We also define matrices $W$ and $S$ by $W_{ij}=w(\gamma_i,\gamma_j)$ and $S_{ij}=iw(\gamma_i,\gamma_j)$, such that the $l_i$ are local coordinates around $X$, and the $k_i$ are local coordinates around $\Gamma$, and both $W$ and $S$ are invertible. Moreover, we have formulas
$$dV_{WP}=\sqrt{\det(W^{-1})}dl_1\wedge\cdots\wedge dl_{6g-6+2n}$$
and
$$dV_K=\sqrt{\det(S^{-1})}dk_1\wedge\cdots\wedge dk_{6g-6+2n}$$
From now on, we will slightly stretch notation by writing $\Phi^*dV_K$ as $dV_K$ and $\Phi^*k_i$ as $k_i$, since $\Phi$ is a local diffeomorphism in a neighborhood of $X$. We will also leave off the $|_X$ notation. We want to control
\begin{equation}\frac{dV_K}{dV_{WP}}=\label{eratio}\sqrt{\frac{\det(S^{-1})}{\det(W^{-1})}}\frac{dk_1\wedge\cdots\wedge dk_{6g-6+2n}}{dl_1\wedge\cdots\wedge dl_{6g-6+2n}}\end{equation}

If $X$ is $(B,A)$-good, then by Lemma \ref{geodsetup},
$$\sqrt{\frac{\det(S^{-1})}{\det(W^{-1})}}=O(1)$$
Furthermore, if $X$ is $\eta$-great, we instead have
$$\sqrt{\frac{\det(S^{-1})}{\det(W^{-1})}}=1+O(e^{-2\eta})$$
We now wish to use the coordinate system $l_1,...,l_{6g-6+2n}$ to apply the standard formula
$$dk_i=\sum_{j=1}^{6g-6+2n}\frac{\partial k_i}{\partial l_j}dl_j$$
However, we only know how to take twist derivatives, not length derivatives. To rectify this problem, we will use the fact that if $\tau_j$ is the twist derivative corresponding to $\gamma_j$, then $\frac{dl_i}{d\tau_j}=w(\gamma_i,\gamma_j)=W_{ij}$. As in the proof of Lemma \ref{symplectic}, we take
$$\frac{\partial}{\partial \tau_j}=\sum_{i=1}^{6g-6+2n}\frac{\partial l_i}{\partial \tau_j}\frac{\partial}{\partial l_i}=\sum_{i=1}^{6g-6+2n}W_{ij}\frac{\partial}{\partial l_i}$$
which implies that
$$\frac{\partial}{\partial l_i}=\sum_{j=1}^{6g-6+2n}W^{-1}_{ij}\frac{\partial}{\partial\tau_j}$$
Therefore
\begin{equation}\label{elengthderiv}\sum_{j=1}^{6g-6+2n}\frac{\partial k_i}{\partial l_j}dl_j=\sum_{j=1}^{6g-6+2n}\sum_{k=1}^{6g-6+2n}W_{jk}^{-1}\frac{\partial k_i}{\partial\tau_k}dl_j\end{equation}
Again applying Lemma \ref{derivcalc} and using the combinatorial length estimates in \ref{geodsetup}, this is
$$dk_i=\sum_{j=1}^{6g-6+2n}O(1)dl_j$$
and plugging this bound into $dk_1\wedge\cdots\wedge dk_{6g-6+2n}$ yields $O(1)dl_1\wedge\cdots\wedge dl_{6g-6+2n}$. We've now shown equation \ref{radon1}. If $X$ is $\eta$-great, then recall that Lemma \ref{geodsetup} guarantees the $\gamma_i$ have $O(1)$-bounded combinatorial length, and that Lemma \ref{derivcalc} then allows us to write equation \ref{elengthderiv} as
$$\sum_{j=1}^{6g-6+2n}\sum_{k=1}^{6g-6+2n}W_{jk}^{-1}(si(\gamma_i,\gamma_k)+O(e^{-\eta/2}))dl_j= \sum_{j=1}^{6g-6+2n}\sum_{k=1}^{6g-6+2n}W_{jk}^{-1}(S_{ik}+O(e^{-\eta/2}))dl_j $$
Applying the bounds on $W_{ij}^{-1}$ and $S_{ik}$ from Lemma \ref{geodsetup}, this is
$$\sum_{j=1}^{6g-6+2n}(W^{-1}S)_{ji}dl_j+O(e^{-\eta/2})\sum_{j=1}^{6g-6+2n}dl_j$$
Now,
$$dk_1\wedge\cdots\wedge dk_{6g-6+2n}=\bigwedge_{i=1}^{6g-6+2n}\left(\sum_{j=1}^{6g-6+2n}(W^{-1}S)_{ji}dl_j\right)+O(e^{-\eta/2})\left(\sum_{j=1}^{6g-6+2n}dl_j\right)$$
The errors may be handled by expansion, leading to
$$dk_1\wedge\cdots\wedge dk_{6g-6+2n}=\left(\frac{\det(S)}{\det(W)}+O(e^{-\eta/2})\right)dl_1\wedge\cdots\wedge dl_{6g-6+2n}$$
Plugging this into equation \ref{eratio} and applying Lemma \ref{geodsetup} yet again yields the desired bound. $\blacksquare$

\section{Convergence of Distributions}

\subsection{Critical Exponent Upper Bounds for Graphs}

Before proving our main theorems, we need to establish some coarse critical exponent bounds that will be useful for surfaces and graphs that are not good. For surfaces, the bound is immediate: as mentioned in the introduction, if $X$ is a hyperbolic surface with or without boundary and $\delta_X$ is its critical exponent, $0<\delta_X\leq 1$. For graphs, the situation is more complicated, as graph critical exponents may be arbitrarily large.

\begin{lemma}
\label{critexpupperbound}

Let $\Gamma\in MRG_{g,n}(L_1,...,L_n)$ be a metric ribbon graph such that all geodesics of $\Gamma$ have length at least $\epsilon$. Then if $\delta_\Gamma$ is the critical exponent of $\Gamma$, $\delta_\Gamma=O(1/\epsilon)$.
	
\end{lemma}

\textbf{Proof.} We may ignore the ribbon graph structure on $\Gamma$, and treat it solely as a metric graph. Contracting an edge of $\Gamma$ can only increase the critical exponent; we will iteratively perform such contractions.

Assuming $\Gamma$ has more than one vertex, take a minimal-length edge $e$ of $\Gamma$ between two distinct vertices. Any geodesic on $\Gamma$ either does not contain $e$, or contains $e$ and at least one other edge, which is is at least as long as $e$. Therefore, contracting $e$ reduces the length of each geodesic on $\Gamma$ by at most a factor of $\frac{1}{2}$. Moreover, the number of such contractions that can be performed until $\Gamma$ has only one vertex is bounded by $g$ and $n$. So we can find a one-vertex graph $\Gamma'$ that is obtained from contracting $\Gamma$, such that all geodesics (and hence all loops) on $\Gamma'$ have length at least $\Omega(\epsilon)$. From here, replace $\Gamma'$ with the graph obtained by shortening each loop on $\Gamma'$ so they all have length equal to the length of the minimal-length loop. Say this length is $\epsilon'$.

For a graph with one vertex and $k$ loops of length $1$, the number of geodesics of length $\leq L$ is $O((2k-1)^L)$, so the critical exponent of such a graph is $O(\ln(2k-1))$. The number of loops of $\Gamma'$ is bounded in terms of $g$ and $n$, and so by scaling the critical exponent of $\Gamma'$ is $O(1/\epsilon)$. In conclusion, $\delta_\Gamma\leq\delta_{\Gamma'}=O(1/\epsilon)$. $\blacksquare$

\subsection{Distribution Convergence on Moduli Space}

\begin{lemma}
\label{volumecontrol}

Fix $g$, $n$, and $C$, let $V_{g,n}(L_1,...,L_n)$ be the Weil-Petersson volume of $\M_{g,n}(L_1,...,L_n)$, and let $V^K_{g,n}(L_1,...,L_n)$ be the Kontsevich volume of $MRG_{g,n}(L_1,...,L_n)$. Assume that the $L_i$ are $C$-controlled by a parameter $L$. Then
$$V_{g,n}^K(L_1,...,L_n)=V_{g,n}(L_1,...,L_n)\left(1+O\left(\frac{1}{L}\right)\right)$$
	
\end{lemma}

\textbf{Proof.} By definition,
$$V_{g,n}(L_1,...,L_n)-V_{g,n}^K(L_1,...,L_n)$$
$$=\int_{\M_{g,n}(L_1,...,L_n)}dV_{WP}-\int_{MRG_{g,n}(L_1,...,L_n)}dV_K$$
We can now use the map $\Phi:\M_{g,n}(L_1,...,L_n)\to MRG_{g,n}(L_1,...,L_n)$ to perform a pullback. Key to this operation is the fact that $\Phi$ is a diffeomorphism between a full-$dV_{WP}$-measure locus of $\M_{g,n}(L_1,...,L_n)$ and a full-$dV_K$-measure locus of $MRG_{g,n}(L_1,...,L_n)$. From for the rest of the proof, we will write $\Phi^*dV_K=dV_{K}$. Pulling back yields.
$$\int_{\M_{g,n}(L_1,...,L_n)}dV_{WP}-dV_K$$
The strategy is to decompose $\M_{g,n}(L_1,..,L_n)$ into $\M_I$, $\M_{NG}(B)$, $\M_{G}(B,A)$, and $\M_{GG}(B,A)$ as in definition \ref{surfaceclassification}. First, the irrelevant locus has measure zero with respect to both $dV_{WP}$ and $dV_K$ by Corollary \ref{trivalentvol}, so we may ignore it. For $\M_{NG}(B)$, we undo the pullback and bound its contribution by
\begin{equation}
\label{efullvol}
\vol_{\M_{g,n}(L_1,...,L_n)}(\M_{NG}(B))+\vol_{MRG_{g,n}(L_1,...,L_n)}(\Phi(\M_{NG}(B)))
\end{equation}
By Lemma \ref{ribcontrol}, we may choose $B$ sufficiently large with respect to $n$ and $g$ that assuming $\Phi(X)$ has no geodesics of length $\leq B$ implies $X$ has no geodesics with length $\leq \frac{B}{2}$. By Lemmas \ref{shortgeodcontrol} and \ref{shortloopgraphvol}, Equation \ref{efullvol} is bounded by $O\left(\frac{B^2}{L^2}\right)$, which since $B$ depends only on $g$, $n$, and $C$, may be written as $O\left(\frac{1}{L^2}\right)$. Next, we apply Radon-Nikodym to the good and great loci:
$$\int_{\M_G(B,A)\cup\M_{GG}(B,A)}dV_{WP}-dV_K$$
\begin{equation}
\label{epredyad}
=\underbrace{\int_{\M_G(B,A)}\left(1-\frac{dV_K}{dV_{WP}}\right)dV_{WP}}_{\ref{epredyad}-A}+\underbrace{\int_{\M_{GG}(B,A)}\left(1-\frac{dV_K}{dV_{WP}}\right)dV_{WP}}_{\ref{epredyad}-B}
\end{equation}
For the integral \ref{epredyad}-$A$, we will apply Lemma \ref{ratiocontrol} to replace the integrand with $O(1)$, and then Lemma \ref{shortedgecontrol} to bound the volume of $\M_G(B,A)$ by $O\left(\frac{1}{L}\right)V_{g,n}(L_1,...,L_n)$. For the integral \ref{epredyad}-$B$, we will apply a dyadic decomposition. Specifically, define
$$\M_{GG}^k(B,A)\coloneq\M_{GG}(B,2^kA)\setminus\M_{GG}(B,2^{k+1}A)$$
so that
$$\M_{GG}(B,A)=\bigcup_{k=0}^\infty \M_{GG}^k(B,A),$$
although for any fixed $L$ all but finitely many of these sets will be empty. Now, applying Lemma \ref{shortedgecontrol},
$$\vol_{\M_{g,n}(L_1,...,L_n)}(\M_{GG}^k(B,A))\leq \vol_{\M_{g,n}(L_1,...,L_n)}(\M_{GG}(B,2^{k+1}A)^c)$$
$$=O\left(\frac{2^{k+1}A}{L}\right)V_{g,n}(L_1,...,L_n) $$
Furthermore, on $\M_{GG}^k(B,A)\subset \M_{GG}(B,2^kA)$, Lemma \ref{ratiocontrol} yields
$$1-\frac{dV_K}{dV_{WP}}=O(e^{-{2^{k-1}A}})$$
Therefore
$$\int_{\M_{GG}(B,A)}\left(1-\frac{dV_K}{dV_{WP}}\right)dV_{WP}\leq O\left(\sum_{k=0}^\infty \frac{2^{k+1}Ae^{-2^{k-1}A}}{L}\right)V_{g,n}(L_1,...,L_n)$$
The sum $\sum_{k=0}^\infty 2^{k+1}Ae^{-2^{k-1}A}$ converges to a value uniformly bounded in terms of $A$, so we can conclude that
$$V_{g,n}^K(L_1,...,L_n)-V_{g,n}(L_1,...,L_n)$$
$$=O\left(\frac{1}{L}\right)V_{g,n}(L_1,...,L_n)\text{ $\blacksquare$}$$

\

We may now prove our main result.

\begin{theorem}
\label{mainresult2}

Fix $g$, $n$, and $C$, and consider the probability space $(\M_{g,n}(L_1,...,L_n),d\mu_{WP})$, with the $L_i$ $C$-controlled by a parameter $L$. Also let $f:\reals\to\reals$ be an arbitrary $1$-Lipschitz function satisfying $f(0)=0$. Then for $L$ sufficiently large,
$$\int_{\M_{g,n}(L_1,...,L_n)}\left|f(\delta_X)-f(\Phi^*\delta_\Gamma)\frac{\Phi^*d\mu_K}{d\mu_{WP}}\right|d\mu_{WP}=O\left(\frac{1}{L^{13/12}}\right)$$
	
\end{theorem}

\textbf{Proof.} As before, we will write $\delta_\Gamma$ for $\Phi^*\delta_\Gamma$ and $d\mu_K$ for $\Phi^*d\mu_K$. We now subdivide the domain on slightly different scales than in Lemma \ref{volumecontrol}, employing $\M_I$, $\M_{NG}(L^{1/3})$, $\M_G(L^{1/3},L^{1/4})$, and $\M_{GG}(L^{1/3},L^{1/4})$. This $L$ will be taken sufficiently large so that $L^{1/4}\geq B$ and $L^{1/3}\geq A$, and that $L^{1/3}$ is a sufficiently large multiple of $L^{1/4}$ in the sense of Lemma \ref{goodgeodframe}.

For $\M_{NG}(L^{1/3})$, we bound its contribution by
\begin{equation}
\label{esplitcontribution}
\underbrace{\int_{\M_{NG}(L^{1/3})}|f(\delta_X)|d\mu_{WP}}_{\ref{esplitcontribution}-A}+\underbrace{\int_{\M_{NG}(L^{1/3})}|f(\delta_\Gamma)|\frac{d\mu_K}{d\mu_{WP}}d\mu_{WP}}_{\ref{esplitcontribution}-B}
\end{equation}
Since $|f(\delta_X)|\leq 1$, the integral \ref{esplitcontribution}-$A$ is at most $\vol_{\M_{g,n}(L_1,...,L_n)}(\M_{NG}(L^{1/3}))$, which by Lemma \ref{shortgeodcontrol} is $O\left(\frac{L^{2/3}}{L^2}\right)=O\left(\frac{1}{L^{4/3}}\right)$. For the integral \ref{esplitcontribution}-$B$, we use $|f(\delta_\Gamma)|\leq \delta_\Gamma$, undo the Radon-Nikodym derivative, and push forward along $\Phi$:
$$\leq\int_{\Phi(\M_{NG}(L^{1/3}))}\delta_\Gamma d\mu_K$$
We will now define
$$\M_{NG}^k(L^{1/3})=\Phi(\M_{NG}(2^{-k}L^{1/3}))\setminus \Phi(\M_{NG}(2^{-(k+1)}L^{1/3}))$$
Then
$$\Phi(\M_{NG}(L^{1/3}))=\bigcup_{k=0}^\infty\M_{NG}^k(L^{1/3})$$
Each $\M^k_{NG}(L^{1/3})$ has volume less than $\M_{NG}(2^{-k}L^{1/3})$, which is at most $O\left(\frac{2^{-2k}}{L^{4/3}}\right)$ by Lemma \ref{shortgeodcontrol}. Furthermore, on $\M_{NG}^k(L^{1/3})$, $|\delta_\Gamma|\leq O\left(\frac{2^{k+1}}{L^{1/3}}\right)$ by Lemma \ref{critexpupperbound}. Therefore
$$\int_{\Phi(\M_{NG}(L^{1/3}))}\delta_\Gamma d\mu_K=O\left(\sum_{k=0}^\infty\frac{2^{-k+1}}{L^{5/3}}\right)=O\left(\frac{1}{L^{5/3}}\right)$$
To handle good surfaces, we note that if $X\in\M_G(L^{1/3},L^{1/4})$, then by Corollary \ref{goodbridge} and Lemma \ref{critexpupperbound}, we have $\delta_X\leq\delta_{\Phi(X)}\leq \frac{1}{L^{1/3}}$, so $|f(\delta_X)|,|f(\delta_\Gamma)|\leq \frac{1}{L^{1/3}}$. Furthermore, by Lemma \ref{volumecontrol}, $\frac{V_{g,n}^K(L_1,...,L_n)}{V_{g,n}(L_1,...,L_n)}=1+O\left(\frac{1}{L}\right)$, so we can apply Lemma \ref{ratiocontrol} to conclude $\frac{d\mu_K}{d\mu_{WP}}=O(1)$ on $\M_{G}(L^{1/3},L^{1/4})$. So by Lemma \ref{shortedgecontrol},
$$\int_{\M_G (L^{1/3},L^{1/4})}\left|f(\delta_X)-f(\delta_\Gamma)\frac{d\mu_K}{d\mu_{WP}}\right|d\mu_{WP}=O\left(\frac{1}{L^{1/3}}\int_{\M_G(L^{1/3},L^{1/4})}d\mu_{WP}\right)$$
$$=O\left(\frac{1}{L^{1/3}}\frac{L^{1/4}}{L}\right)=O\left(\frac{1}{L^{13/12}}\right)$$
We now define
$$\M_{GG}^k(L^{1/3},L^{1/4})=\M_{GG}^k(L^{1/3},2^kL^{1/4})\setminus\M_{GG}(L^{1/3},2^{k+1}L^{1/4})$$
By Lemma \ref{critexpcontrol}, $0<\delta_X\leq 1$, and the Lipschitz condition, on $\M_{GG}^k(L^{1/3},L^{1/4})$, 
$$f(\delta_\Gamma)=f(\delta_X)+O\left(\delta_X\frac{e^{-2^{k-1}L^{1/4}}}{2^kL^{1/4}}\right)=f(\delta_X)+ O\left(\frac{e^{-2^{k-1}L^{1/4}}}{2^kL^{1/4}}\right) $$
Furthermore, by Lemmas \ref{ratiocontrol} and \ref{volumecontrol}, on $\M^k_{GG}(L^{1/3},L^{1/4})$,
$$\frac{d\mu_K}{d\mu_{WP}}=\frac{V_{g,n}^K(L_1,...,L_n)}{V_{g,n}(L_1,...,L_n)}\frac{dV_K}{dV_{WP}}$$
$$=\left(1+O\left(\frac{1}{L}\right)\right)\left(1+O\left(e^{-2^{k-1}L^{1/4}}\right)\right)$$
For sufficiently large $L$, this is $1+O\left(\frac{1}{L}\right)$. Lastly, a graph with all edges of length at least $2^kL^{1/4}$ clearly has no geodesics of length at most $ 2^kL^{1/4}$, so Lemma \ref{critexpupperbound} and $\delta_X\leq\delta_{\Phi(X)}$ tell us that $|f(\delta_X)|\leq O\left(\frac{1}{2^kL^{1/4}}\right)\leq O\left(\frac{1}{L^{1/4}}\right)$ on $\M^k_{GG}(L^{1/3},L^{1/4})$. Now,
$$\left|f(\delta_X)-f(\delta_\Gamma)\frac{d\mu_K}{d\mu_{WP}}\right|=\left|f(\delta_X)-f(\delta_\Gamma)\left(1+O\left(\frac{1}{L}\right)\right)\right|$$
$$=\left|f(\delta_X)-\left(f(\delta_X)+ O\left(\frac{e^{-2^{k-1}L^{1/4}}}{2^kL^{1/4}}\right)\right) \left(1+O\left(\frac{1}{L}\right)\right)\right|$$
Expanding and using $|f(\delta_X)|\leq O\left(\frac{1}{L^{1/4}}\right)$ yields
$$O\left(\frac{e^{-2^{k-1}L^{1/4}}}{2^kL^{1/4}}\right)+O\left(\frac{1}{L^{5/4}}\right)$$
Lastly, by Lemma \ref{shortedgecontrol}, 
$$\prob_{\M_{g,n}(L_1,...,L_n)}(\M_{GG}^k(L^{1/3},L^{1/4}))=O\left(\frac{2^k}{L^{3/4}}\right)$$
Therefore
$$\int_{\M_{GG}(L^{1/3},L^{1/4})}\left|f(\delta_X)-f(\delta_\Gamma)\frac{d\mu_K}{d\mu_{WP}}\right|d\mu_{WP}$$
$$=\sum_{k=0}^\infty\int_{\M_{GG}^k (L^{1/3},L^{1/4})} \left|f(\delta_X)-f(\delta_\Gamma)\frac{d\mu_K}{d\mu_{WP}}\right|d\mu_{WP}$$
\begin{equation}\label{etwointegrand}=\sum_{k=0}^\infty \int_{\M^k_{GG}(L^{1/3},L^{1/4})}O\left(\frac{e^{-2^{k-1}L^{1/4}}}{2^kL^{1/4}}\right)+O\left(\frac{1}{L^{5/4}}\right)d\mu_{WP}\end{equation}
We may now integrate the two integrand terms of Equation \ref{etwointegrand} separately. For the first, we apply the aforementioned volume bound.
$$\sum_{k=0}^\infty \int_{\M^k_{GG}(L^{1/3},L^{1/4})}O\left(\frac{e^{-2^{k-1}L^{1/4}}}{2^kL^{1/4}}\right)d\mu_{WP}=\sum_{k=0}^\infty O\left(\frac{e^{-2^{k-1}L^{1/4}}}{L}\right)$$
It is now straightforward to check that this sum is at most $O\left(\frac{e^{-L^{1/4}}}{L}\right)$, which in turn is certainly less than, say, $O\left(\frac{1}{L^{5/4}}\right)$ for sufficiently large $L$. For the second integral-sum, we observe that the integrand does not depend on $k$, so we may undo the dyadic decomposition:
$$\int_{\M_{GG}(L^{1/3},L^{1/4})}O\left(\frac{1}{L^{5/4}}\right)d\mu_{WP}$$
We now use the fact that $\prob_{\M_{g,n}(L_1,...,L_n)}(\M_{GG} (L^{1/3},L^{1/4}))\leq 1$ to obtain a bound of $O\left(\frac{1}{L^{5/4}}\right)$. $\blacksquare$

\subsection{Corollaries}

From here, we can easily obtain Corollaries \ref{cor1} and \ref{cor2}. For Corollary \ref{cor1}, set $f(x)=x$ and recall that rescaling preserves the Kontsevich form while multiplying the graph critical exponent by $\alpha$. Corollary \ref{cor2} is slightly more involved. We need the following dual representation of the Wasserstein metric:

\begin{theorem}[Kantorovich Duality, Wasserstein case]
Fix two probability measures $\mu$ and $\nu$ on $\reals$. Then
$$d_W(\mu,\nu)=\sup\left\{\int_\reals f(x)d(\mu-\nu)(x): \emph{Lip}(f)\leq 1\right\}$$
where $\emph{Lip}(f)$ is the Lipschitz constant of $f$.
\end{theorem}

A proof may be found in \cite{CV03}.

\

\textbf{Proof of Corollary \ref{cor2}.}
First, observe that we may always replace the test function $f(x)$ used in the dual representation of the Weierstrass norm with $f(x)-f(0)$ without changing the value of the difference of integrals. We can therefore assume without loss of generality that $f(0)=0$. Recall that we wish to bound
$$\left|\int_\reals f(x)(\alpha\delta_X)_*(d\mu_{WP})-\int_\reals f(x)(\delta_\Gamma)_*(d\mu_{K})\right|$$
Undo the pushforward and rescale:
$$=\left|\int_{\M_{g,n}(\alpha L_1,...,\alpha L_n)}f(\alpha\delta_X)d\mu_{WP}-\int_{MRG_{g,n}(L_1,...,L_n)}f(\delta_\Gamma)d\mu_K\right|$$
$$=\left|\int_{\M_{g,n}(\alpha L_1,...,\alpha L_n)}f(\alpha\delta_X)d\mu_{WP}-\int_{MRG_{g,n}(\alpha L_1,...,\alpha L_n)}f(\alpha \delta_\Gamma)d\mu_K\right|$$
$$=\left|\int_{\M_{g,n}(\alpha L_1,...,\alpha L_n)}\left(f(\alpha\delta_X)-f(\alpha\delta_\Gamma)\frac{d\mu_K}{d\mu_{WP}}\right)d\mu_{WP}\right|$$
$$\leq \alpha \int_{\M_{g,n}(\alpha L_1,...,\alpha L_n)}\left|1/\alpha f(\alpha\delta_X)-1/\alpha f(\alpha\delta_\Gamma)\frac{d\mu_K}{d\mu_{WP}}\right|d\mu_{WP}$$
Since $f(x)$ is $1$-Lipschitz, so is $1/\alpha f(\alpha x)$, and we may apply Theorem \ref{maintheorem} to obtain the desired bound. $\blacksquare$

\pagebreak

\appendix

\section{Trigonometric Formulas}

\begin{figure}[h]

\centering
\includegraphics[scale=0.15]{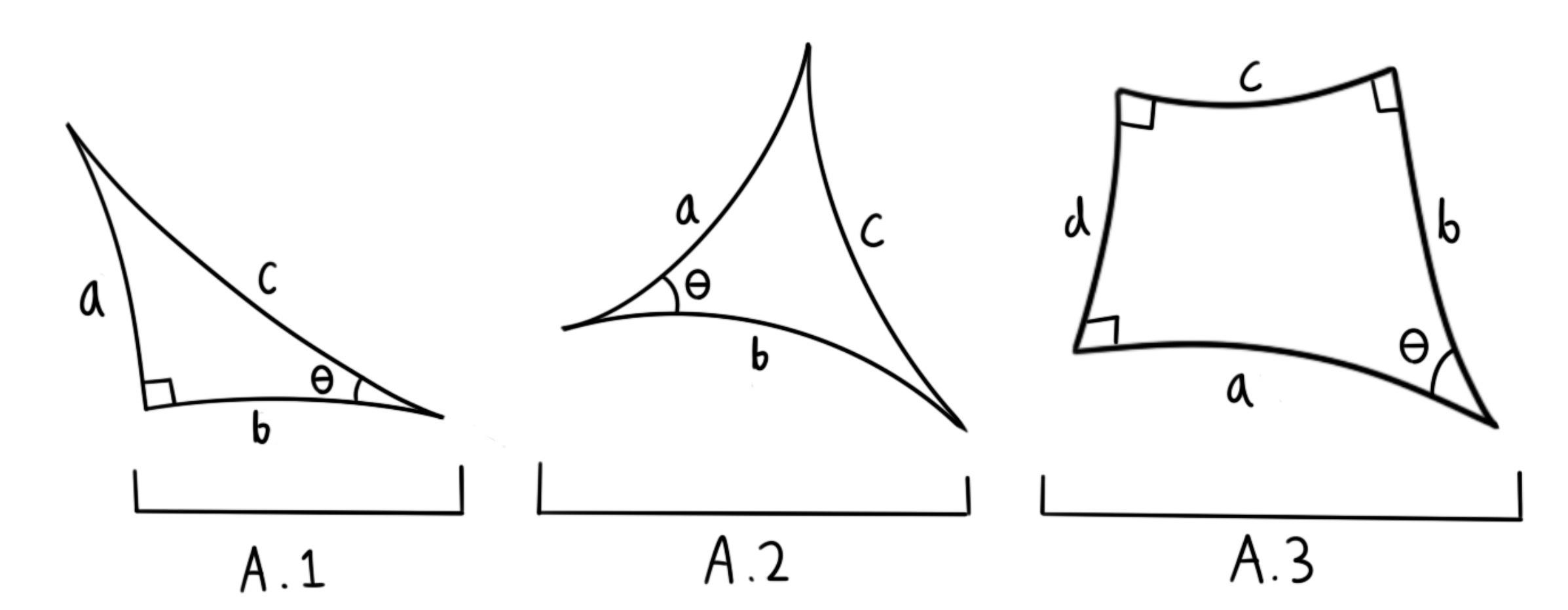}
\caption{The geometric setups used in Lemmas \ref{righttriangle}, \ref{triangle}, and \ref{quadrilateral}.}

\label{trig}

\end{figure}

The following formulas are all standard. For a concise reference, see \url{https://www.maths.gla.ac.uk/wws/cabripages/hyperbolic/hyperbolic0.html}. Furthermore, the geometric objects under consideration are depicted in Figure \ref{trig}. 

\begin{lemma}
\label{righttriangle}

Let $T$ be a hyperbolic triangle with sides $a$, $b$, and $c$, and assume that $T$ has a right angle between sides $a$ and $b$. Let $\theta$ be the angle between sides $b$ and $c$. Then
$$\sin(\theta)=\frac{\sinh(l(b))}{\sinh(l(c))}$$
	
\end{lemma}

\begin{lemma}
\label{triangle}

Let $T$ be a hyperbolic triangle with sides $a$, $b$, and $c$, and let $\theta$ be the angle between sides $a$ and $b$. Then
$$\cos(\theta)=\coth(l(a))\coth(l(b))-\frac{\cosh(l(c))}{\sinh(l(a))\sinh(l(b))}$$
	
\end{lemma}

\begin{lemma} 
\label{quadrilateral}

Let $Q$ be a hyperbolic quadrilateral with sides $a$, $b$, $c$, and $d$, arranged in that order. Let $\theta$ be the angle between $a$ and $b$, and assume all other angles are right. Then:

\begin{enumerate}

\item $\cos(\theta)=\sinh(c)\sinh(d)$, and in particular $|\sinh(c)\sinh(d)|\leq 1$.
\item $\tanh(b)=\cosh(c)\tanh(d)$, and in particular $|\cosh(c)\tanh(d)|\leq 1$.
	
\end{enumerate}
	
\end{lemma}

\section{Non-Trivalency}
\label{nontrivalent}

We now deal with non-trivalent ribbon graphs. Let $\M_I$ be the locus of surfaces $X\in\M_{g,n}(L_1,...,L_n)$ satisfying the property that $\Phi(X)$ is not trivalent. Equivalently, $\Phi(\M_I)$ is the locus of non-trivalent graphs in $MRG_{g,n}(L_1,...,L_n)$. We would like to know that the Weil-Petersson volume of $\M_I$ and the Konsevich volume of $\Phi(\M_I)$ are both $0$. For the Kontsevich form, this is true by definition; see Definition 3.1 in \cite{ABCGLW21}.

Calculating the Weil-Petersson volume presents more difficulties. First, we may sort graphs on $MRG_{g,n}(L_1,...,L_n)$ according to their topological isomorphism classes (i.e. forgetting the edge length information but keeping the ribbon structure); this induces a cell decomposition on $MRG_{g,n}(L_1,...,L_n)$ that lifts to the Teichm{\"u}ller space $TRG_{g,n}(L_1,...,L_n)$. Since $\Phi$ and its corresponding Teichm{\"u}ller variant are homeomorphisms, we obtain topological cellularizations of $\M_{g,n}(L_1,...,L_n)$ and $\T_{g,n}(L_1,...,L_n)$.

Now, a \emph{maximal arc system} $\alpha$ on a topological surface $S_{g,n}$ is a choice of simple geodesic arcs $\iota_1,...,\iota_{k}$ on $S_{g,n}$ such that each arc travels between two boundary arcs, and such that cutting along the $\iota_j$ decomposes $X$ into hexagons. By a genus calculation, such a system must contain exactly $6g-6+3n$ arcs. For any marked surface $X\in \T_{g,n}(L_1,...,L_n)$, the marking gives us a way to realize $\alpha$ on $X$ such that all arcs are length-minimizing in their homotopy classes (relative to the boundary).

The following result is due to Luo. Any maximal arc system $\alpha$ on a surface $S_{g,n}$ generates a corresponding \emph{ideal graph} $\Gamma_\alpha$, where the vertices correspond to hexagons obtained by cutting along the arc system and each arc generates an edge between vertices corresponding to hexagons on the two sides of that arc. A \emph{fundamental edge cycle} is an edge cycle in $\Gamma_\alpha$ such that each edge  appears at most twice, and each boundary component of $X$ generates a corresponding \emph{boundary edge cycle} by tracking which arcs are incident to the boundary (with multiplicity). Then:

\begin{theorem}[Luo, \cite{FL07}, Thm. 1.2]
\label{widthcoord}

Fix a maximal arc system $\alpha$ on a surface $S_{g,n}$, and choose positive real numbers $L_1,...,L_n$. For each arc $\iota\in\alpha$, one may define a \emph{width coordinate} $w_\iota:\T_{g,n}(L_1,...,L_n)\to\reals$, and the map $W_\alpha$:
$$W_\alpha:\T_{g,n}(L_1,...,L_n)\to \reals^{6g-6+3n}$$
$$W_\alpha(X)=(w_\iota)_{\iota\in\alpha}$$
is a diffeomorphism with respect to any choice of Fenchel-Nielsen coordinates from $\T_{g,n}(L_1,...,L_n)$ onto the convex polytope $P_\alpha$ defined by the following inequalities. First, for every fundamental cycle $\iota_1,...,\iota_k$, 
\begin{equation}\label{efundcycle}w_{\iota_1}+\cdots+w_{\iota_k}>0\end{equation}
Second, for every boundary edge cycle $\iota_1,...,\iota_k$ coming from the $j$th boundary component, 
\begin{equation}\label{eboundarycycle}w_{\iota_1}+\cdots+w_{\iota_k}=L_j\end{equation}
	
\end{theorem}

Note that since $P_\alpha$ is a convex polytope inside an affine linear subspace, some subset of the $w_\iota$ serve as global coordinates on $P_\alpha$.

While we will not give the full definition of the functions $w_\iota$, their key property is as follows:

\begin{theorem}[Guo-Luo, \cite{GL11}]
\label{widthisedge}

Choose a cell $C$ of the cellularization of $\T_{g,n}(L_1,...,L_n)$ such that $C$ corresponds to a trivalent ribbon graph. For any $X$ in that cell, define a maximal arc system $\alpha$ by pulling back the intercostals of the spine $S_X$ to $S_{g,n}$ via the marking map. Then $\alpha$ is independent of the $X$ chosen. Furthermore, each arc $\iota\in\alpha$ is associated to an edge $e\in S_X$ via edge-intercostal duality, and $w_\iota(X)$ is the length of that edge in $\Phi(X)$.

Lastly, there exists an open neighborhood $U$ of the closure of $C$ such that for all $X\in U$, $X$ fails to be trivalent if and only if $w_a(X)=0$ for some $\iota\in\alpha$.
	
\end{theorem}

In other words, on each cell we may find a maximal arc system so that the width coordinates become edge-length coordinates for the corresponding cell in $TRG_{g,n}(L_1,...,L_n)$, and with these coordinates the walls of the cell are described by simple linear equations. The first claim of this lemma is essentially Theorem 2 in \cite{GL11}, although one needs to unravel how Guo and Luo define the cell decomposition. The second claim is a straightforward consequence of Lemma 5 in \cite{GL11} (see in particular the discussion in section 5.2 of this paper). 

\begin{corollary}
\label{analytic}

Fix $X\in\M_{g,n}(L_1,...,L_n)$, let $\Gamma=\Phi(X)$, and assume $\Gamma$ is trivalent. Assume we have simple closed geodesics $\tilde{\eta}_1,...,\tilde{\eta}_{6g-6+2n}$ on $\Gamma$ such that the length functions $k_i=l(\tilde{\eta}_i)$ form local coordinates around $\Gamma$. Then there is a neighborhood $U$ of $X$ and a neighborhood $V$ of $\Gamma$ such that $\Phi(U)=V$ and $\Phi|_U$ is a diffeomorphism from a choice of Fenchel-Nielsen coordinates to the $(k_1,...,k_{6g-6+2n})$ coordinates.
	
\end{corollary}

\textbf{Proof.} Theorem \ref{widthcoord} tells us that any choice of $6g-6+2n$ width coordinates with full rank in $P_\alpha$ are locally diffeomorphic to any choice of Fenchel-Nielsen coordinates, and Theorem \ref{widthisedge} tells us that the map $W_\alpha$ precisely gives the edge lengths of $\Phi(X)$. The $k_i$ coordinates are therefore linear with respect to some set of width coordinates.

\

Furthermore, let $K(L_1,...,L_n)$ be the affine linear subspace of $\reals^{6g-6+3n}$ defined by the equations \ref{efundcycle} and \ref{eboundarycycle}. Any top-dimensional cell $C$ lies inside of $K(L_1,...,L_n)$. If any hyperplane defined by $w_a=0$ did not intersect $K(L_1,...,L_n)$ transversely, $K(L_1,...,L_n)$ would have to lie inside this hyperplane, implying no $X\in C$ would be trivalent. But this contradicts the fact that the map in Theorem 2 of \cite{FL07} is surjective.

In conclusion, cell boundaries of the cellularization of $\T_{g,n}(L_1,...,L_n)$ are contained in lower-dimensional submanifolds with respect to coordinates that are diffeomorphic to Fenchel-Nielsen coordinates. Since the Weil-Petersson form becomes the Euclidean volume form in Fenchel-Nielsen coordinates by Wolpert's magic formula, and $\T_{g,n}(L_1,...,L_n)$ covers $\M_{g,n}(L_1,...,L_n)$, we may conclude:

\begin{corollary}
\label{trivalentvol}

$\vol_{\M_{g,n}(L_1,...,L_n)}(\M_I)=0$.
	
\end{corollary}

\printbibliography

@misc{AHC22,
      title={The shapes of complementary subsurfaces to simple closed hyperbolic multi-geodesics}, 
      author={Francisco Arana-Herrera and Aaron Calderon},
      year={2022},
      eprint={2208.04339},
      archivePrefix={arXiv},
      primaryClass={math.GT},
      url={https://arxiv.org/abs/2208.04339}, 
}

@book{AT10,
  title     = "Zeta Functions of Graphs: A Stroll through the Garden",
  author    = "Audrey Terras",
  year      = 2010,
  publisher = "Cambridge University Press",
  address   = "Cambridge"
}

@article{GL11,
  author  = "Ren Guo and Feng Luo",
  title   = "Cell decompositions of Teichm{\"u}ller spaces of surfaces with boundary",
  journal = "Pacific Journal of Mathematics",
  year    = 2011,
  volume  = "253",
  number  = "2",
  pages   = "423--438"
}

@misc{ABCGLW21,
      title={On the Kontsevich geometry of the combinatorial Teichm\"uller space}, 
      author={Jørgen Ellegaard Andersen and Ga{\"e}tan Borot and Séverin Charbonnier and Alessandro Giacchetto and Danilo Lewa{\'n}ski and Campbell Wheeler},
      year={2021},
      eprint={2010.11806},
      archivePrefix={arXiv},
      primaryClass={math.DG},
      url={https://arxiv.org/abs/2010.11806}, 
}

@misc{ND10,
      title={The asymptotic Weil-Petersson form and intersection theory on moduli space}, 
      author={Norman Do},
      year={2010},
      eprint={1010.4126},
      archivePrefix={arXiv},
      primaryClass={math.GT},
      url={https://arxiv.org/abs/1010.4126}, 
}

@article{FL07,
  author  = "Feng Luo",
  title   = "On Teichm{\"u}ller spaces of surfaces with boundary",
  journal = "Duke Math J.",
  year    = 2007,
  volume  = "139",
  number  = "3",
  pages   = "463--482"
}

@article{MM07,
  author  = "Maryam Mirzakhani",
  title   = "Simple geodesics and Weil-Petersson volumes of moduli spaces of bordered Riemann surfaces",
  journal = "Invent. Math.",
  year    = 2007,
  volume  = "167",
  number  = "1",
  pages   = "179--222"
}

@article{MM08a,
  author  = "Maryam Mirzakhani",
  title   = "Ergodic theory of the earthquake flow",
  journal = "Int. Math. Res. Not.",
  year    = 2008,
  volume  = "2008",
  number  = "9",
  pages   = "1--39"
}

@article{MM08b,
  author  = "Maryam Mirzakhani",
  title   = "Growth of the number of simple closed geodesics on hyperbolic surfaces",
  journal = "Ann. of Math.",
  year    = 2008,
  volume  = "168",
  number  = "1",
  pages   = "97--125"
}

@article{MM13,
  author  = "Maryam Mirzakhani",
  title   = "Growth of Weil-Petersson volumes and random hyperbolic surfaces of large genus",
  journal = "J. Differential Geom.",
  year    = 2013,
  volume  = "94",
  number  = "2",
  pages   = "267--300"
}

@article{GM09,
  author  = "Gabriele Mondello",
  title   = "Triangulated Riemann surfaces with boundary and the Weil-Petersson Poisson Structure",
  journal = "J. Differential Geom.",
  year    = 2009,
  volume  = "81",
  pages   = "391--436"
}

@article{SL89,
  author  = "Steven P. Lalley",
  title   = "Renewal theorems in symbolic dynamics, with applications to geodesic flows, noneuclidean tessellations and their fractal limits",
  journal = "Acta Math.",
  year    = 1989,
  volume  = "163",
  pages   = "1--55"
}

@article{BE88,
  author  = "Brian Bowditch and David Epstein",
  title   = "Natural triangulations associated to a surface",
  journal = "Topology",
  year    = 1988,
  volume  = "27",
  number  = "1",
  pages   = "91--117"
}

@article{JH86,
  author  = "John L. Harer",
  title   = "The virtual cohomological dimension of the mapping class group of an orientable surface",
  journal = "Invent. Math.",
  year    = 1986,
  volume  = "84",
  number  = "1",
  pages   = "157--176"
}

@article{MK92,
  author  = "Maxim Kontsevich",
  title   = "Intersection theory on the moduli space of curves and the matrix Airy function",
  journal = "Comm. Math. Phys.",
  year    = 1992,
  volume  = "147",
  number  = "1",
  pages   = "1--23"
}

@article{RP87,
  author  = "R. C. Penner",
  title   = "The decorated Teichm{\"u}ller space of punctured surfaces",
  journal = "Comm. Math. Phys.",
  year    = 1987,
  volume  = "113",
  number  = "2",
  pages   = "299--339"
}

@article{SW83,
  author  = "Scott Wolpert",
  title   = "On the symplectic geometry of deformations of a hyperbolic surface",
  journal = "Ann. of Math.",
  year    = 1983,
  volume  = "117",
  number  = "2",
  pages   = "207--234"
}

@article{SWP83,
  author  = "Scott Wolpert",
  title   = "On the K{\"a}hler form of the moduli space of once punctured tori",
  journal = "Comm. Math. Helvetici",
  year    = 1983,
  volume  = "58",
  pages   = "246--256"
}

@article{AS56,
    AUTHOR = {Alte Selberg},
     TITLE = {Harmonic analysis and discontinuous groups in weakly symmetric
              {R}iemannian spaces with applications to {D}irichlet series},
   JOURNAL = {J. Indian Math. Soc. (N.S.)},
  FJOURNAL = {The Journal of the Indian Mathematical Society. New Series},
    VOLUME = {20},
      YEAR = {1956}}

@book{CV03,
  title     = "Topics in Optimal Transportation",
  author    = "C\'edric Villani",
  year      = 2003,
  publisher = "American Mathematical Society",
  address   = "Providence"
}

@article{FN12,
AUTHOR = {Anna Felikson and Sergey Natanzon},
     TITLE = {Moduli via double pants decompositions},
   JOURNAL = {Diff. Geo. and its Applications},
    VOLUME = {30},
    ISSUE = {5},
      YEAR = {2012}
}

@article{SK83,
AUTHOR = {Steven Kerckhoff},
     TITLE = {The Neilsen realization problem},
   JOURNAL = {Ann. of Math.},
    VOLUME = {177},
      YEAR = {1983},
    PAGES = {235-265}
}

@article{CM98,
AUTHOR = {Curtis McMullen},
     TITLE = {Complex earthquakes and Teichm{\"u}ller Theory},
   JOURNAL = {Jour. of the Amer. Math. Soc.},
    VOLUME = {11},
    ISSUE = {2},
      YEAR = {1998},
    PAGES = {283-320},
}

@article{DGZZ21,
  author  = "Vincent Delecroix and {\'E}lise Goujard and Peter Zograf and Anton Zorich",
  title   = "Masur-Veech volumes, frequencies of simple closed geodesics, and intersection numbers of moduli spaces of curves",
  journal = "Duke Math. J.",
  year    = 2021,
  volume  = "170",
  number  = "12",
  pages   = "2633--2718"
}

@article{ABCDGLW19,
  author  = "Jørgen Ellegaard Andersen and Ga{\"e}tan Borot and Séverin Charbonnier and Vincent Delecroix and Alessandro Giacchetto and Danilo Lewa{\'n}ski and Campbell Wheeler",
  title   = "Topological recursion for Masur-Veech volumes",
  journal = "Jour. of the London Math. Soc.",
  year    = 2022,
  volume  = "107",
  number  = "1",
  pages   = "254-332"
}

@phdthesis{AG21,
 type = {Doctoral Dissertation},
  author       = {Alessandro Giacchetto}, 
  title        = {Geometric and topological recursion and invariants of the moduli space of curves},
  school       = {University of Bonn},
  year         = 2021,
  url = {https://bonndoc.ulb.uni-bonn.de/xmlui/bitstream/handle/20.500.11811/9385/6402.pdf?sequence=1&isAllowed=y#page305}
}

@article{MP19,
AUTHOR = {Maryam Mirzakhani and Bram Petri},
     TITLE = {Lengths of closed geodesics on random surfaces of large genus},
   JOURNAL = {Comm. Math. Helvetici},
    VOLUME = {94},
    ISSUE = {4},
      YEAR = {2019},
    PAGES = {869--889},
}

@misc{BGL23,
      title={Length spectrum of large genus random metric maps}, 
      author={Simon Barazer and Alessandro Giacchetto and Mingkun Liu},
      year={2023},
      eprint={2312.10517},
      archivePrefix={arXiv},
      primaryClass={math.PR},
      url={https://arxiv.org/abs/2312.10517}, 
}

@article{JL23,
AUTHOR = {Svante Janson and Baptiste Louf},
     TITLE = {Unicellular maps vs. hyperbolic surfaces in large genus: simple closed curves},
   JOURNAL = {Ann. Prob.},
    VOLUME = {51},
    ISSUE = {3},
      YEAR = {2023},
    PAGES = {899-929},
}

\end{document}